%% file: rapid-pdes-main.tex
\documentclass{article}
\usepackage{amssymb}
\usepackage{amsthm}
\usepackage{amsmath}
\usepackage[cp1250]{inputenc}
\usepackage[a4paper, left=3cm, right=3cm, top=3cm, bottom=3cm]{geometry}

\newtheorem{theorem}{Theorem}[section]
\newtheorem{lemma}[theorem]{Lemma}
\newtheorem{rem}[theorem]{Remark}
\newtheorem{definition}[theorem]{Definition}

\theoremstyle{definition}

\newcommand{\inte }{{\rm int}\,}
\newcommand{\dom }{{\rm dom}\,}
\newcommand{\dist }{\,{\rm dist}\,}

\providecommand{\subspaceH}{H^\prime}

\def\comment#1{{}}

\def\qed{{\hfill{\vrule height5pt width3pt depth0pt}\medskip}}

\begin{document}
\begin{center}
{\Large \bf  Stabilizing effect of large average initial velocity in forced dissipative PDEs invariant with respect to Galilean transformations}

 \vskip 0.5cm
 {\large
Jacek~Cyranka$^{\ddag,*}$,
Piotr~Zgliczy\'nski$^{*}$}

 \vskip 0.5cm

{\small$^*$ Institute of Computer Science and Computational Mathematics,
Jagiellonian University}\\
{\small  S. {\L}ojasiewicza 6, 30-348 Krak\'ow, Poland}

 \vskip 0.5cm

{\small$^{\ddag}$ Institute of Applied Mathematics and Mechanics, University of Warsaw}\\
{\small  Banacha 2, 02-097 Warszawa, Poland}

 \vskip 0.5cm

 jacek.cyranka@ii.uj.edu.pl, piotr.zgliczynski@ii.uj.edu.pl

\vskip 0.5cm

\today
\end{center}

\begin{abstract}
 We describe a topological method to study the dynamics of dissipative PDEs on a torus with rapidly oscillating forcing terms.
We show that a dissipative PDE, which is invariant with respect to  Galilean transformations,  with a large average initial velocity can be reduced 
to a problem with
rapidly oscillating forcing terms. We apply the technique to the Burgers equation, and the incompressible 2D Navier-Stokes equations
with a time-dependent forcing. We prove that for  a large initial average speed  the equation admits
a bounded eternal solution, which attracts all other solutions forward in time.
For the incompressible 3D Navier-Stokes equations  we establish existence of a locally attracting solution.

\end{abstract}


\input intro.tex

\input linNonAuto.tex

\input infdim.tex

\input rapid-osc-infdim.tex


\input burgers.tex

\input fast-mov-burgers.tex

\input ns.tex


\section{Acknowledgments}
Research has been supported by Polish
National Science Centre grant 2011/03B/ST1/04780.
The presented work has been done while the first author (JC) held a post-doctoral position at
Warsaw Center of Mathematics and Computer Science.

\input ref.tex
\end{document}

%% file: intro.tex
\section{Introduction}


Let us consider a \emph{dissipative partial differential equation} (PDE) with periodic boundary conditions  in $\mathbb{R}^d$ (i.e. on the $d$-dimensional torus $\mathbb{T}_d$) of the following form
\begin{equation}
  \frac{d u}{dt} = \nu \Delta u + N(\nabla u,u) + f(t,x) \label{eq:pde-abstract}
\end{equation}
where $\nu >0$, $u \in \mathbb{R}^d$ and  $f \in \mathbb{R}^d$ is the forcing term.

Additionally, we assume that (\ref{eq:pde-abstract}) is invariant with respect the transformation $(t,x,u) \mapsto (t,x+at,u+a)$, where $a \in \mathbb{R}^d$
(this is \emph{the Galilean transformation} to a coordinate frame moving with the velocity $a$) and
we have the following conservation law
\begin{equation}
  \frac{d}{dt} \int_{\mathbb{T}^d} u(t,x)=  \int_{\mathbb{T}^d} f(t,x).
\end{equation}
If we assume that
\begin{equation}
  \int_{\mathbb{T}^d} f(t,x)=0,
\end{equation}
then
\begin{equation}
 \frac{1}{(2\pi)^d} \int_{\mathbb{T}^d} u(t,x)= u_0.  \label{eq:u-ave-const}
\end{equation}

For the problem (\ref{eq:pde-abstract},\ref{eq:u-ave-const}), for large $\|u_0\|$ we  are interested in the existence of a bounded solution,
which looks like
$u(t,x)=u_0 + \mathcal{O}\left(\frac{1}{\|u_0\|}\right)$ defined for all $t \in \mathbb{R}$ (termed an \emph{eternal solution}), which attracts all other solutions forward in time.
The size of $\mathcal{O}\left(\frac{1}{\|u_0\|}\right)$ term depends on the magnitude of forcing $f$, its time derivatives and the viscosity $\nu$.
The important point is that we do not assume smallness of $f$.

The basic idea in our approach can be described as follows:
for generic choices of $u_0$  the transformation moving (\ref{eq:pde-abstract}) into
a coordinate frame in which the average of $u$ vanishes, leads to a new form of (\ref{eq:pde-abstract}) with a rapidly oscillating forcing term. We prove that in the transformed system
this very rapid oscillation is effectively equivalent to small forcing term of the size $\mathcal{O}\left(\frac{1}{\|u_0\|}\right)$ (for a fixed forcing $f$ and $\nu$). As the result we obtain
an absorbing set very close to \emph{zero} and then by topological reasoning we show that it contains a attracting orbit bounded by $\mathcal{O}(1/\|u_0\|)$, defined for all
$t \in \mathbb{R}$.
By reversing  the initial coordinate change we obtain an attracting orbit of the form $u_0 + \mathcal{O}\left(\frac{1}{\|u_0\|}\right)$.
We exploit the fact that the rapid oscillations could be effectively 'integrated out' (or averaged) without assuming the smallness of oscillating
term. This idea is known for some time in the numerical analysis, where it is used to obtain effective quadratures and ODE solvers for system with rapid oscillations
(see \cite{I} and the references cited there), and in the averaging theory developed by Bogolyubov, Neishstadt and others (see for example \cite{BM,BZ,FW,Nei84,Wa}).

As an application of our approach we chose two models: one dimensional Burgers equation and the Navier-Stokes system in dimension two and three. For the Burgers equation in dimension
one we established the result described above. For the Navier-Stokes equations in 2D, the result is true for a generic direction of $u_0$ and
in 3D for a generic direction $u_0$ we establish the existence of small locally attracting solution. Each of the mentioned equations required slightly different
approach.

Similar results to ours can be found in  \cite{JKM}. There
for any $\nu>0$ the authors established the existence of a
globally attracting solution for 1D viscous Burgers equation with   periodic boundary conditions, under assumption that forcing is periodic in time.
Without any proof they state that this solution scales like $u_0 + \mathcal{O}(1/\|u_0\|)$ (for fixed forcing $f$ and $\nu$).
Our result applies to  more general forcings, moreover, we
are able to establish the
exponential convergence rate to the attracting solution, while in
\cite{JKM} the authors clearly indicated that they cannot make
such claim and they asked for the convergence rate in one of the
stated problems \cite[Problem 3(i)]{JKM}. The method in \cite{JKM}
appears to be restricted to the scalar equation on one-dimensional
domains, partially due to the use of the maximum principles.

As it was already mentioned we use some variant of averaging in our approach. There exists a large literature on the application of averaging to PDEs
(see for example \cite{Bam03a,Bam05,Cw12,HVL90,He,Mat01,MS03,Mat08,Pr05}). In the available literature there are  basically two approaches to averaging. The first one, in the spirit
of Bogolyubov and his coworkers \cite{BM,BZ} asks for coordinate change, which 'absorbs' the leading term of perturbation. The second one, which we will attribute to Henry \cite{He},
shows in a direct way that the influence of rapid oscillations is small. Our approach belongs to the second category.
Since we are after global results, we  cannot stay in the framework provided by Henry \cite{He}, but we need to construct a priori bounds, to which
the averaging principle is applied. To handle this we use the method of \emph{self-consistent bounds} developed in \cite{ZM,ZKS2,ZKS3}, which was proposed for computer assisted proofs
in dynamics of Kuramoto-Sivashinski PDE. Recently, this method was further applied to prove the existence of a globally attracting solution for 1D
viscous Burgers equation with some particular choices of nonautonomous forcings \cite{Cy}, \cite{CyZ}.

Let us briefly describe our approach, and layout the structure of the paper.
To establish the existence of a globally attracting solution for (\ref{eq:pde-abstract},\ref{eq:u-ave-const}) with large $\|u_0\|$, we proceed as follows.
First, we establish an averaging lemma  for a finite dimensional non-autonomous non-homogenous linear ODE.
Second, we generalize the result for finite dimensional nonlinear ODE problems.
We describe our results for finite dimensional problems in Section~\ref{sec:basic-estm}.
Third, using the self-consistent bounds approach
in the dissipative PDE setting, which we recall in Section~\ref{sec:self-consist-bnds}, the result is established for
each Galerkin approximation of the PDE and then using high compactness of trapping regions
we pass to the limit with the approximation dimension, see Section~\ref{sec:lem-rapid-osc}.
This third step is where our method may fail, because it requires the existence of the compact absorbing set.
Sections~\ref{sec:def-burgers} and ~\ref{sec:burgers}
contain the application of our method to the viscous Burgers equation, with Theorem~\ref{thm:burgers-main} being the main result.
In Section~\ref{sec:NSE} we study the Navier-Stokes equation on $2$D and $3$D tori, the main result is Theorem~\ref{thm:ns-main}.

\subsection{Notation}
Consider   nonautonomous ODE
\begin{equation}
  x'(t)=f(t,x(t)), \label{eq:notation-Non_auto}
\end{equation}
 where $x \in \mathbb{R}^n$ and $f: \mathbb{R} \times \mathbb{R}^n \to \mathbb{R}^n$ is regular enough to guarantee
 the uniqueness of the initial value problem $x(t_0)=x_0$.  We set
 $\varphi(t_0,t,x_0)=x(t_0+t)$, where $x(t)$ is a solution of (\ref{eq:notation-Non_auto}) with  initial condition  $x(t_0)=x_0$.
 Obviously in each context it will be clearly stated what is the ordinary differential equation generating $\varphi$. We will
 sometimes refer to $\varphi$ as to the local process generated by  (\ref{eq:notation-Non_auto}).

 If a solution $x(t)$ of (\ref{eq:notation-Non_auto}) is defined for all $t \in \mathbb{R}$, then we will call it an \emph{eternal solution} (or occasionally an \emph{eternal orbit}).

 For  function $f(t,z)$ we will often use $D_tf$ and $D_z f$
to denote the partial derivatives. For example, $D_z
f(t,z)=\frac{\partial f}{\partial z}(t,z)$.

For the partial derivatives of $f: \mathbb{R}^n \to \mathbb{R}$ will sometimes use the following notation $D_\alpha f$, where $\alpha \in \{1,\dots,n\}^l$ is a multiindex of length $|\alpha|=l$, to denote $\frac{\partial^{|\alpha|}f }{\partial x_{\alpha_1}\cdot \partial x_{\alpha_l}}$.

 For  matrix $U$ by $U^t$ we will denote its transpose. For a square matrix $U$ we will denote its spectrum by $\mbox{Sp}(U)=\{ \lambda \in \mathbb{C} \ | \ \lambda \ \mbox{is an eigenvalue of $U$}\}$. If $A \in \mathbb{R}^{n\times n}$ by $\mu(A)$ we will denote its logarithmic norm \cite{D,HNW,L,KZ}, which is defined
 by
 \begin{equation}
   \mu(A) = \lim_{h \to 0^+} \frac{\|I + hA\|- 1}{h}.
 \end{equation}
 For the properties of logarithmic norm and its relation with the Lipschitz constant for the flow induced by ODEs see \cite{KZ} and
the literature cited there.
 The logarithmic norm depends on the norm used. For the Euclidean norm the logarithmic norm of a matrix is given as follows
 \begin{equation}
   \mu(A)=\max \{ \lambda, \quad \lambda \in \mbox{Sp}((A+A^t)/2) \}. \label{eq:log-norm-eucl}
 \end{equation}

 For a vector field  $f(t,x)$ as in (\ref{eq:notation-Non_auto}) of class $C^1$  and a set $J \times W \subset \mathbb{R} \times \mathbb{R}^n$ we define
 \begin{equation}
   \mu(D_x f,J \times W ) = \sup_{(t,x) \in J \times W} \mu(D_xf(t,x)).
 \end{equation}
 For an autonomous system $x'=f(x)$ and  $W \subset \mathbb{R}^n$ we set analogously
 \begin{equation}
   \mu(D f, W ) = \sup_{x \in  W} \mu(D f(x)).
 \end{equation}

\comment{
Quite often the following expression will appear
$h_l(t)=\frac{e^{lt}-1}{l}$, where $l$ is a fixed parameter. For
$l=0$ we will understand $h_0$ as follows $h_0(t)=\lim_{l \to
0}\frac{e^{lt}-1}{l}=t$.  \textbf{PZ: nigdzie w tekscie nie znalazlem $h_l$, a sam to wprowadzilem}
}

We will use the following notation, for $d \in \mathbb{N}$ we set
\begin{equation}
  S_d (p) = 1+ \sum_{k \in \mathbb{Z}^d \setminus \{0\}} \frac{1}{|k|^p}.
\end{equation}
The sum above converges for $p > d$. The one added the sum is responsible for term with $k=0$.

 (
By $\mathbb{T}_d$ we will denote the $d$-dimensional torus, i.e. $\mathbb{T}_d=\left(\mathbb{R}/2\pi\right)^d$.

\begin{definition}
      Let $\nmid\cdot\nmid\colon\mathbb{R}\to\mathbb{R}$ be given by
      \begin{equation*}
        \nmid a\nmid:=\left\{\begin{array}{ll}|a|&\text{ if }a\neq 0,\\1&\text{ if }a=0.\end{array}\right.
      \end{equation*}
    \end{definition}

%% file: linNonAuto.tex

\section{Basic estimates for oscillations in finite dimension}
\label{sec:basic-estm}

\subsection{Linear nonautonomous equations}
\label{subsec:lin-non-auto}

Assume that $A: \mathbb{R} \to \mathbb{R}^{n \times n}$ is continuous and
for $k=1,\dots,m$ $v_k:\mathbb{R}  \to \mathbb{R}$
are $C^1$, $g_k:\mathbb{R} \to \mathbb{R}$ are continuous.

Let us consider the following non-autonomous non-homogenous linear ODE
\begin{equation}
  x'(t)=A(t)x(t) + \sum_{k \in J_V}  g_k(\omega_k t) v_k(t),  \quad x \in \mathbb{R}^n.  \label{eq:linNonHomoNonAuto}
\end{equation}
The set $J_V$ in the sum in (\ref{eq:linNonHomoNonAuto}) might be  finite or infinite, or the sum might be an integral over some measure on $J_V$.

For each $k \in J_V$ let $G_k(t)$ be a primitive of $g_k$, so
\begin{equation}
 G_k'(t)=g_k(t).
\end{equation}
We will assume later that $|G_k(t)|$ are  bounded. This is the reflection of the oscillating nature of $g_k$.

Let $M(t,t_0)$ be a fundamental matrix of solutions of the homogenous version of (\ref{eq:linNonHomoNonAuto})
\begin{equation}
   x'(t)=A(t)x(t). \label{eq:linNonAuto}
\end{equation}
This means that for any $t_0 \in \mathbb{R}$ and $x_0 \in \mathbb{R}^n$ the function $x(t)=M(t,t_0)x_0$ solves
(\ref{eq:linNonAuto}) with the initial condition $x(t_0)=x_0$.

It is well known that $M$ has the following properties
\begin{eqnarray}
 M(t_0,t_0)&=&I, \\
 M(t,t_0)^{-1}&=&M(t_0,t), \\
 \frac{\partial}{\partial t} M(t,t_0)&=&A(t) M(t,t_0), \\
  \frac{\partial}{\partial t_0} M(t,t_0)&=&- M(t,t_0) A(t_0). \label{eq:Mfund-4}
\end{eqnarray}

 The general solution of (\ref{eq:linNonHomoNonAuto}) is given by
\begin{eqnarray}
 \varphi(t_0,t,x_0)&=&M(t_0+t,t_0)x_0 +   \label{eq:int-sol-fund-matrix} \\
   & & \int_{0}^{t} M(t_0+t,t_0+s)  \sum_{k \in J_V} g_k(\omega_k (t_0+s)) v_k(t_0+s)
   ds. \nonumber
\end{eqnarray}

We compute the integral in the above formula as follows. Using the integration by parts and (\ref{eq:Mfund-4}) we obtain for $k\in J_V$
\begin{eqnarray}
 I_k(t+t_0):=\int_{0}^{t}  g_k(\omega_k(t_0+s)) M(t_0+t,t_0+s) v_k(t_0+s)ds=  \label{eq:Dk-lin} \\
 \left. \frac{G(\omega_k(t_0+s))}{\omega_k} M(t_0+t,t_0+s)v_k(t_0+s)\right|^{s=t}_{s=0} +  \nonumber \\
-\frac{1}{\omega_k} \int_{0}^t G_k(\omega_k(t_0+s)) \frac{\partial }{\partial s}\left(M(t_0+t,t_0+s) v_k(t_0+s)\right)ds= \nonumber \\
 \frac{1}{\omega_k}(G_k(\omega_k(t_0+t)) v_k(t_0+t)- G_k(\omega_k t_0)M(t_0+t,t_0)v_k(t_0)) + \nonumber \\
 \frac{1}{\omega_k} \int_{0}^t G_k(\omega_k(t_0+s)) M(t_0+t,t_0+s)A(t_0+s)v_k(t_0+s)ds + \nonumber \\
- \frac{1}{\omega_k} \int_{0}^t G_k(\omega_k(t_0+s))
M(t_0+t,t_0+s)u'_k(t_0+s)ds \nonumber
\end{eqnarray}

For Galerkin projections of dissipative PDEs while $\|A\|$ will not have any uniform bound independent of the projection dimension, we expect $\|Av_k\|$ to be uniformly bounded.

Therefore we have proved that for the process generated by (\ref{eq:linNonHomoNonAuto}) it holds
\begin{equation}
   \varphi(t_0,t,x_0)=M(t_0+t,t_0)x_0 + \sum_{k \in J_V} I_k(t+t_0)
\end{equation}


\subsection{Estimates for a nonlinear problem}
\label{subsec:non-prob1}

Assume that $F: \mathbb{R} \times \mathbb{R}^{n} \to \mathbb{R}^n$ is $C^1$ function and
for $k=1,\dots,m$ $v_k:\mathbb{R} \times \mathbb{R}^n \to \mathbb{R}$
are $C^1$, $g_k:\mathbb{R} \to \mathbb{R}$ are continuous.

Consider problem
\begin{equation}
  x' = \tilde{F}(t,x):= F(t,x) + \sum_{k \in J_V}  g_k(\omega_k t)  v_k(t,x)\label{eq:prob1-loc}
\end{equation}
and its oscillation-free version
\begin{equation}
  y' = F(t,y).  \label{eq:auto-prob1-loc}
\end{equation}


\begin{lemma}
\label{lem:prob1-estm-norm}
Let  $x:[t_0,t_0 + h] \to \mathbb{R}^n$ and $y:[t_0,t_0+h] \to \mathbb{R}^n$ be solutions to (\ref{eq:prob1-loc}) and (\ref{eq:auto-prob1-loc}),
respectively, such that $x(t_0)=y(t_0)$.

Let $W$ be a  compact set, such that for any  $t \in [0,h]$ the segment joining $x(t_0+t)$ and $y(t_0+t)$ is contained in $W$.

Assume that for $k \in J_V$ $G_k'(t)=g_k(t)$.

Assume that there exist constants $l$, $C(\dots)$  such that for all $k \in J_V$ holds
\begin{eqnarray}
   \sup_{t \in \mathbb{R}}\|G_k(t)\| &=& C(G_k), \\
   \sup_{t \in \mathbb{R}}|g_k(t)| &=& C(g_k), \\
  \sup_{z \in W, s \in [t_0,t_0+h]} \mu(D_z F(s,z))&=&l  \\
  \sup_{z \in W, s \in [t_0,t_0+h]} \|v_k(s,z)\|&=&C(v_k)  \\
  \sup_{z,z_1 \in W, s \in [t_0,t_0+h]} \|(D_z F(s,z)) v_k(s,z_1) \| &=& C\left(D_z F v_k\right) \\
  \sup_{z \in W, s \in [t_0,t_0+h]} \left\|\frac{\partial v_k}{\partial t} (s,z) \right\| &=&  C\left(\frac{\partial v_k}{\partial t}\right) \\
 \sup_{z \in W, s \in [t_0,t_0+h]} \left\|(D_z v_k (s,z)) \tilde{F}(s,z)  \right\| &=&  C(D_z v_k  \tilde{F} )
\end{eqnarray}

Assume that
\begin{eqnarray}
  \sum_{k\in J_V} C(G_k) C(v_k) < \infty  \label{eq:rosc-conv1} \\
    \sum_{k\in J_V} C(G_k) C\left(D_z F v_k\right) < \infty \\
     \sum_{k\in J_V} C(G_k) C\left(\frac{\partial v_k}{\partial t}\right) < \infty \\
      \sum_{k\in J_V} C(G_k) C\left(D_z v_k \tilde{F}\right) < \infty. \label{eq:rosc-conv4}
\end{eqnarray}

Then for $t \in [0,h]$ holds
\begin{equation}
 \|x(t_0+t)-y(t_0+t)\| \leq \sum_{k\in J_V} \frac{1}{|\omega_k|} b_k(t) 
   \label{eq:estm-diff-osc}
 \end{equation}
where  continuous functions $b_k:[0,h] \to \mathbb{R}_+$ depend on constants  $l$, $C(g_i)$, $C(G_i)$, $C(v_i)$, $C\left(D_z F v_i\right)$, $C\left(\frac{\partial v_i}{\partial t}\right)$ and
 $C\left(D_z v_i\tilde{F}\right)$ as follows
 \begin{eqnarray}
    b_k(t)&=&C(v_k) C(G_k)(1+e^{lt}) + C\left(D_z F v_k\right) C(G_k)(e^{lt}-1)/l +  \label{eq:bk-expression} \\
    & & C(G_k) \left(C\left(\frac{\partial v_k}{\partial t}\right) + C\left(D_z v_k \tilde{F}\right) \right)
(e^{lt}-1)/l \nonumber
 \end{eqnarray}


If $\inf_{k \in J_V} |\omega_k| >0$, the sum in (\ref{eq:estm-diff-osc}) is convergent.

\end{lemma}

\noindent
\textbf{Proof:}
Let $z(t)=x(t)-y(t)$. We have
\begin{eqnarray*}
 z'(t)=F(x(t)) - F(y(t)) +  \sum_{k\in J_V}  g_k(\omega_k t) v_k(t,x) =   \label{eq:prob1-diff}\\
   \left(\int_0^1 D_x F(t,s(x(t)-y(t)) + y(t))ds\right) \cdot z(t) + \sum_{k\in J_V}  g_k(\omega_k t) v_k(t,x(t)).
\end{eqnarray*}
Therefore
\begin{equation}
  z'(t) =  A(t)z(t) + \sum_{k\in J_V}  g_k(\omega_k t) v_k(t,x(t)),
\end{equation}
where
\begin{equation}
  A(t)=  \left(\int_0^1 D_x F(t,s(x(t)-y(t)) + y(t))ds\right).
\end{equation}
Let $M(t_1,t_0)$ is the fundamental matrix of solutions for $x'=A(t)x$.

Since $z(t_0)=0$, then from (\ref{eq:int-sol-fund-matrix}) and (\ref{eq:Dk-lin}) it follows that
\begin{equation}
  z(t_0+t)= \sum_{k\in J_V} I_k(t+t_0)
\end{equation}
where  
\begin{eqnarray*}
  I_k(t+t_0)=\frac{1}{\omega_k}(G_k(\omega_k(t_0+t)) v_k(t_0+t,x(t+t_0))- G_k(\omega_k t_0) M(t_0+t,t_0)v_k(t_0,x(t_0)) + \\
 \frac{1}{\omega_k} \int_{0}^t G_k(\omega_k(t_0+s)) M(t_0+t,t_0+s)A(t_0+s)v_k(t_0+s,x(t_0+s))ds + \\
- \frac{1}{\omega_k} \int_{0}^t G_k(\omega_k(t_0+s)) M(t_0+t,t_0+s)\left(\frac{d}{ds}v_k(t_0+s,x(t_0+s))\right)ds
\end{eqnarray*}

From the standard estimate for the
logarithmic norms (see for example Lemma 4.1 in \cite{KZ}) we know that
for $t \geq 0$ holds
\begin{equation*}
  \|M(t+t_0,t_0)\| \leq \exp(l t).
\end{equation*}
Hence we obtain the following estimate of $I_k(t)$ for $t \in [0,h]$ and $k\in J_V$
\begin{eqnarray*}
  |\omega_k| \cdot \|I_k(t+t_0)\| \leq C(v_k) C(G_k)(1+e^{lt}) + C(D_z F v_k ) C(G_k)  \int_0^t e^{l(t-s)}ds + \\
     C(G_k) \left(C\left(\frac{\partial v_k}{\partial t}\right) + \sup_{s \in [0,h]} \left\|\frac{\partial v_k}{\partial z}(t_0+s,x(t_0+s)) x'(t_0+s)\right\|  \right) \int_0^t e^{l(t-s)}ds \leq \\
  C(v_k) C(G_k)(1+e^{lt}) + C(D_z F v_k) C(G_k)(e^{lt}-1)/l + \\
    C(G_k) \left(C\left(\frac{\partial v_k}{\partial t}\right) + C\left(D_z v_k \tilde{F}\right) \right)
(e^{lt}-1)/l
\end{eqnarray*}

This proves (\ref{eq:bk-expression}).

To finish the proof observe that if $|\omega_k| > \epsilon>0$ for all $k \in J_V$, then assumptions (\ref{eq:rosc-conv1}--\ref{eq:rosc-conv4}) together with the formula (\ref{eq:bk-expression}) imply the convergence of the sum in (\ref{eq:estm-diff-osc}).
\qed

%% file: infdim.tex
\section{Self-consistent bounds for non-autonomous dissipative PDEs}
\label{sec:self-consist-bnds}

The goal of this section is recall from \cite{ZKS3} the technique of self-consistent bounds in order to show that the results from Section~\ref{sec:basic-estm}
can be carried over also in the context of the dissipative PDEs with periodic boundary conditions. In \cite{ZKS3} the problem was autonomous, hence the need
to extend some definitions, lemmas and theorems to the present content. The changes turn out to be minor.

First we recall some definitions and results from Section 2 and 3 in \cite{ZKS3}.

\subsection{The problem}
\label{subsec:problem}
 We consider PDEs of the
following type
\begin{equation}
  u_t = L u + N(t,u,Du,\dots,D^ru) + V(t,x), \label{eq:genpde}
\end{equation}
where $u \in \mathbb{R}^n$, $x \in \mathbb{T}_d$, is an
$d$-dimensional torus), $L$ is a linear operator,  $N$ - a real
polynomial of $u,Du,\dots,D^ru$ with bounded time-dependent coefficients, here by $D^s u$ we denote $s$-th order derivative of
$u$, i.e. the collection of all partial derivatives of $u$ of
order $s$. $V$ is a smooth ($C^\infty$)  real function of both variables.

Due to the periodic boundary conditions we will use the Fourier basis $\{e^{ikx}\}_{k \in \mathbb{Z}^d}$ and we use
the notation
\begin{equation*}
e_k=\exp(ikx).
\end{equation*}


 We require that $L$ is diagonal in the Fourier basis, namely
\begin{equation}
  L e_k= \lambda_k e_k, \label{eq:Ldiag}
\end{equation}
and the eigenvalues $\lambda_k$ satisfy
\begin{eqnarray}
  \lambda_k &=& - \beta(|k|) |k|^p  \label{eq:lambdak} \\
  0 &<& \beta_0 \leq \beta(|k|)  \leq \beta_1, \qquad \mbox{for $|k| > K_-$} \label{eq:lambdak2} \\
  p &>& r.  \label{eq:p>r}
\end{eqnarray}
 The fact that we are considering  functions on the torus  means that we impose periodic
 boundary conditions.
We may also seek odd or even solutions or impose some other
conditions.

 If $u(t,x)$ is a sufficiently regular solution of
(\ref{eq:genpde}), then  we can expand it in Fourier series
$u(t,x)=\sum_{k \in \mathbb{Z}^d} u_k(t)e^{ik\cdot x}$ to obtain
an infinite ladder of ordinary differential equations for the
coefficients $u_k$
\begin{equation}
  \frac{d u_k}{dt}=F_k(t,u):=\lambda_k u_k + N_k(t,u) + V_k(t), \quad k \in
  \mathbb{Z}^d,  \label{eq:fugenpde}
\end{equation}
where $N_k(t,u)$ is $k$-th Fourier coefficient of function
$N(t,u,Du,\dots,D^ru)$.

Observe that  $u_k \in \mathbb{C}^n$ and equations in (\ref{eq:fugenpde}) are
not independent, because the reality of $u$ imposes the following
condition
\begin{equation}
  u_{-k}=\overline{u}_{k}.  \label{eq:reality}
\end{equation}
Observe that we will also have a reality condition for $V$ and $N(u)$.



\subsection{The method of self-consistent bounds}
\label{subsec:method}

We  begin with an abstract nonlinear evolution equation in a real
Hilbert space $H$ ($L^2$ or some its subspaces in our treatment of
dissipative PDEs) of the form
\begin{equation}
\label{eq:pde} \frac{du}{dt} = F(t,u),
\end{equation}
where the domain of $F$ is  dense in $\mathbb{R} \times H$. By a solution of
(\ref{eq:pde}) we understand a function $u:[t_0,t_{max}) \to
H$, such that $(t,u(t)) \in \dom F $ for $t \in [t_0,t_{max}]$ such that $u$ is differentiable and (\ref{eq:pde}) is
satisfied for all $t \in [t_0,t_{max})$.

The scalar product in $H$  will be denoted by $(u | v)$.
Throughout the paper we assume that there is a set $I \subset
\mathbb{Z}^d$ and a sequence of subspaces $H_k \subset H$ for $k
\in I$, such that $\dim H_k=d_1 < \infty$ and $H_k$ and $H_{k'}$
are mutually orthogonal for $k \neq k'$. Let $A_k: H \to H_k$ be
the orthogonal projection onto $H_k$. We assume that for each $u
\in H$ holds
\begin{equation}
   u=\sum_{k \in I} u_k=\sum_{k \in I} A_k u. \label{eq:H-decmp}
\end{equation}
The above equality  for a given $u \in H$ and $k \in I$ defines
$u_k$. Analogously if $B$ is a function with the range in $H$,
then $B_k(u)=A_k B(u)$. Equation (\ref{eq:H-decmp}) implies that
$H=\overline{\bigoplus_{k \in I} H_k}$.

For $k \in \mathbb{Z}^d$ we define
\begin{displaymath}
  |k|=\sqrt{\sum_{i=1}^d k_i^2}
\end{displaymath}

 For $n > 0$ we set
\begin{eqnarray*}
  X_n = \bigoplus_{|k| \leq n, k \in I } H_k \\
  Y_n = X_n^\bot,
\end{eqnarray*}
by $P_n:H \to X_n$ and $Q_n:H\to Y_n$ we will denote the
orthogonal projections onto $X_n$ and onto $Y_n$, respectively.

\begin{definition}
\label{defn:Faddmissible} We say that $F:\mathbb{R} \times H \supset \dom (F) \to H$
is admissible  if the following conditions are satisfied for any
$i \in \mathbb{R}$, such that $\dim X_i >0$
\begin{itemize}
\item $ \mathbb{R} \times X_i \subset \dom(F)$
\item $P_i F : \mathbb{R} \times X_i \to X_i$ is a $C^1$ function
\end{itemize}

For an admissible map $F$ and $i \in \mathbb{R}$, such that $\dim X_i >0$ we define a projection of $F$, as a map $\mathcal{P}_i F: \mathbb{R} \times X_i \to X_i $ by $(\mathcal{P}_i F)(t,x)=P_i F(t,x) $ for $x \in X_i$.
\end{definition}

\begin{definition}
Assume $F$ is admissible. For a given number $n>0$ the ordinary
differential equation
\begin{equation}
  x'=P_n F(t,x), \qquad x \in X_n  \label{eq:galproj}
\end{equation}
will be called \emph{the $n$-th Galerkin projection} of
(\ref{eq:pde}).

By $\varphi^n(t_0,t,x)$ we denote the local process on $X_n$ induced by
(\ref{eq:galproj}).
\end{definition}

\begin{definition}
\label{defn:selfconsistent} Assume $F$ is an admissible function.
Let $m,M\in \mathbb{R}$ with $m\leq M$. Let $S \subset \mathbb{R}$ be a connected set.   Consider an object
consisting of:  a compact set $W\subset X_m$ and a sequence of
compact sets $B_k \subset H_k$ for $|k|
> m$, $k \in I$. We define the conditions
\textbf{C1}, \textbf{C2}, \textbf{C3}, \textbf{C4a} as follows:
\begin{description}
\item[C1] For $|k|>M$, $k \in I$ holds $0 \in B_k$ .
\item[C2] Let $\hat{a}_k := \max_{a \in B_k} \|a\|$ for $|k| > m$, $k \in I$ and then
$\sum_{|k| > m, k \in I} \hat{a}_k^2 < \infty$. In particular
\begin{equation}
   W \oplus \Pi_{|k| > m} B_k \subset H
\end{equation}
and for every $u \in  W \oplus \Pi_{k \in I,|k| > m} B_k $ holds,
$\| Q_n u \| \leq \sum_{|k|>n, k \in I}\hat{a}_k^2$.
\item[C3] The function $(t,u) \mapsto F(t,u)$ is continuous on
$S \times W\oplus\prod_{k \in I, |k|>m}B_k \subset H$.

 Moreover, if we define for $k
\in I$, $F_k=\max_{(t,u) \in S \times W \oplus \prod_{k\in I, |k| >m} B_k}
|F_k(t,u)| $, then $\sum F_k^2 < \infty$.

\item[C4a] For $|k|>m$, $k \in I$  $B_k$ is given by (\ref{eq:C4aball}) or (\ref{eq:C4aint})
\begin{eqnarray}
  B_k&=&\overline{B(c_k,r_k)},  \quad r_k > 0  \label{eq:C4aball} \\
  B_k&=&\Pi_{s=1}^d[a_{s}^-,a_s^+], \qquad a_s^- < a_s^+, \:
  s=1,\dots,d_1 \label{eq:C4aint}
\end{eqnarray}
 Let $(t,u) \in S \times W \oplus \Pi_{|k| > m}B_k$.
 Then for $|k| > m$ holds:
\begin{itemize}
\item if $B_k$ is given by (\ref{eq:C4aball}) then
   \begin{eqnarray}
     u_k \in  \partial_{H_k} B_k  & \Rightarrow  & (u_k -c_k | F_k(t,u)) <
     0. \label{eq:C4ain}
 \end{eqnarray}
 \item if $B_k$ is given by (\ref{eq:C4aint}) then
  \begin{eqnarray}
     u_{k,s}=a_{k,s}^-  & \Rightarrow  & F_{k,s}(t,u) >   0,  \label{eq:C4ain-} \\
     u_{k,s}=a_{k,s}^+  & \Rightarrow  & F_{k,s}(t,u) <   0.  \label{eq:C4ain+}
 \end{eqnarray}
 \end{itemize}
\end{description}
\end{definition}

In the sequel we will refer to equations (\ref{eq:C4ain}) and
(\ref{eq:C4ain-}--\ref{eq:C4ain+}) as  \emph{isolation equations}.

\begin{definition}
Assume $F$ is an admissible function. Let $m,M\in \mathbb{R}$ with
$m\leq M$. Let $S \subset \mathbb{R}$ be a connected set.  Consider an object consisting of: a compact set
$W\subset X_m$ and a sequence of compacts $B_k \subset H_k$ for
$|k| > m, k \in I $. We say that set $W \oplus \Pi_{k \in I, |k|
>m}B_k$ forms \emph{self-consistent  bounds for
$F$ over the time interval $S$} if conditions C1, C2, C3 are satisfied.

If additionally condition  C4a holds, then  we say that $W \oplus
\Pi_{k \in I,|k| > m} B_k$ forms \emph{topologically
self-consistent bounds for $F$ over the time interval $S$}
\end{definition}
If $F$ and $S$ is clear from the context, then we will often drop $F$ and $S$,
 and we will speak simply  about \emph{self-consistent bounds}
 or \emph{topologically self-consistent bounds}.

 Given self-consistent  bounds $W$ and $\{B_k\}_{k \in I, |k|>m}$,
 by $T$ (the tail) we will denote
\begin{equation}
T:=\prod_{|k| > m} B_k \subset Y_m.
\end{equation}

Here are some useful lemmas from~\cite{ZKS3} illustrating the implications of
conditions C1, C2, C3.

 From condition C2 it follows immediately that
\begin{lemma}
\label{lem:compact} If $W \oplus T$ forms self-consistent bounds,
then $W\oplus T$ is a compact subset of $H$.
\end{lemma}

The following  lemma is an immediate consequence of conditions C2
and C3.

\begin{lemma}
\label{lem:compact2} Given self-consistent  bounds $W \oplus T$,
then
\begin{displaymath}
  \lim_{n\to \infty} P_n(F(t,u)) =F(t,u), \quad \mbox{uniformly for
   $(t,u) \in S \times W \oplus T$}
\end{displaymath}
\end{lemma}

\begin{lemma}
\label{lem:convsubseq2} Let $W_i \oplus T_i$, $i=1,\dots,k$ forms
self-consistent bounds for (\ref{eq:pde}). Let  $\{d_n\}_{n \in
\mathbb{N}} \subset \mathbb{R}$ be a sequence, such that $\lim_{n
\to \infty} d_n = \infty$. Assume that, for all $n$,
$x_n:[t_1,t_2] \to \bigcup_{i=1}^k W_i \oplus T_i$ is a solution
of
\begin{equation}
   \frac{dp}{dt}(t)=P_{d_n}(F(t,p(t))), \qquad p(t) \in X_{d_n}.
\end{equation}
Then there exists a convergent subsequence $\{d_{n_l}\}_{l \in
\mathbb{N}}$ such that,
\newline $\lim_{l \to \infty} x_{n_l} = x^*$, where
$x^*:[t_1,t_2]\to \bigcup_{i=1}^k W_i \oplus T_i$ and the
convergence is uniform on $[t_1,t_2]$. Moreover, $x^*$ satisfies
(\ref{eq:pde}).
\end{lemma}

\providecommand{\gp}{Galerkin projection of \eqref{eq:fugenpde} }
\providecommand{\gps}{Galerkin projections of \eqref{eq:fugenpde} }

For a class of dissipative partial differential equations, in our approach, we will often construct \emph{trapping regions}, and \emph{absorbing sets}.
In the sequel we construct those sets for viscous Burgers' equation and the Navier-Stokes equations.
\begin{definition}
  \label{def:trappingReg}
Let $N_0\geq 0$, $\varphi^n$ be a local process induced by the $n$-th Galerkin projection of
    \eqref{eq:fugenpde}.
    A set $W\subset H$ is called \emph{the  trapping region} for large \gps if it is \emph{forward invariant},
namely if for any pair $(t_0,u_0)\in  \mathbb{R} \times W$, for all $n>N_0$ and all $t\geq t_0$ holds  $\varphi^n\left(t_0,t, P_nu_0\right)\in P_n W$.
\end{definition}

\begin{definition}
    \label{def:absorbingSet}
    Let $N_0\geq 0$, $\varphi^n$ be a local process induced by the $n$-th Galerkin projection of
    \eqref{eq:fugenpde}.
    A set $\mathcal{A}\subset H$ is called \emph{the absorbing set}
    for large \gps, if for any pair $(t_0,u_0)\in  \mathbb{R} \times H$ there exists
    $t_1(u_0)\geq 0$ such that
    for all $n>N_0$ and all $t \in \mathbb{R}$, $t\geq t(u_0)$ holds  $\varphi^n\left(t_0,t, P_nu_0\right)\in P_n
    \mathcal{A}$.
     Moreover, $P_n\mathcal{A}$ is forward invariant for $\varphi^n$.
\end{definition}

\subsection{Estimates}

In the sequel to make some notations and statements shorter and more transparent that whenever we have $\frac{1}{|k|^s}$, then for $k=0$ we mean that
$\frac{1}{|0|^s}$ should be replaced by $1$. With this convention we have the following lemma.

\begin{lemma}[Lemma 3.4 in \cite{ZKS3}]
\label{lem:convolution}
If $\gamma >d$, then there exists constant $C_2(d,\gamma)$ such that for any $k \in \mathbb{Z}^d$
the following holds true
\begin{equation}
  \sum_{k_1,k_2 \in \mathbb{Z}^d, k_1+k_2=k } \frac{1}{|k_1|^\gamma |k_2|^\gamma} \leq \frac{C_2(d,\gamma)}{|k|^\gamma}.
\end{equation}
\end{lemma}

Observe that $N_k$ is a finite sum of the  terms of the following form
\begin{equation}
 g(t,a)=  w(t)  \sum_{j_1+j_2 + \dots +j_l=k} (D^{\alpha_1}a_{p_1})_{j_1} (D^{\alpha_2}a_{p_2})_{j_2} \cdots (D^{\alpha_l}a_{p_l})_{j_{l}}, \label{eq:terms-in-N}
\end{equation}
where $w(t)$ is a coefficient in $N$, $|w(t)| \leq W$, $j_1,\dots,j_l \in \mathbb{Z}^d$,  multi-indices $\alpha_1,\dots,\alpha_l$ such that $|\alpha_i| \leq r$ and $(D^{\alpha}a_{p})_{j}$ is $p$-th component (coordinate)   of the $j$-th Fourier coefficient of $(D^\alpha a)$.

Therefore each term in $N$ give rise to some convolution, which is an infinite sum. We will refer to these sums as the sums defining $N_k$ or $\frac{\partial N_k}{\partial a_j}$.

The following lemma is an obvious generalization of Lemma 3.1 in \cite{ZKS3}.

\providecommand{\D}{\widehat{C}}

\begin{lemma}
\label{lem:Dgen}
Assume that $N(t,a,Da,\dots,D^ra)$ is a polynomial in variables $a,Da,\dots,D^ra$ with time dependent and bounded coefficients (i.e. there exists a constant
$w$, such that for all $t \in \mathbb{R}$ the coefficients in $N$ are less than or equal to $W$).

 Let $s > s_0=d+r$.
 If $|a_k| \leq C/|k|^s$, $|a_0| \leq C$, then
there exists $\D$, which depends on $C$, $s$ and $W$, such that
\begin{equation}
  |N_k| \leq \frac{\D}{|k|^{s-r}}, \qquad |N_0| \leq \D.
\end{equation}
Moreover, the series defining $N_k$ is absolutely uniformly converging.
\end{lemma}
The proof is almost the same as that of Lemma 3.1 in \cite{ZKS3}. The only difference is the following. Where we had previously constant coefficients,
we now just insert the upper bound for the time dependent coefficient.

\begin{lemma}
\label{lem:DNestm}
Assume that $N(t,a,Da,\dots,D^ra)$ is a polynomial in variables $a,Da,\dots,D^ra$ with time dependent and bounded coefficients (i.e. there exists a constant
$W$, such that for all $t \in \mathbb{R}$ the coefficients in $N$ are less than or equal to $W$). Assume that $N$ does not contain constant terms or degree one terms.

 Let $s > s_0=d+r$.
 If $|a_k| \leq C/|k|^s$, $|a_0| \leq C$, then
there exists $\D$, which depends on $C$, $s$ and $W$, such that
 \begin{eqnarray*}
  \left|\frac{\partial N_k}{\partial a_j}(t,a)\right| &\leq& \frac{\D|j|^r}{|k-j|^{s-r}}, \qquad   k \neq j  \\
   \left|\frac{\partial N_k}{\partial a_k}(t,a)\right| &\leq& \D|k|^r.
 \end{eqnarray*}
 Moreover, the series defining $\frac{\partial N_k}{\partial a_j}(t,a)$ is absolutely uniformly converging.
\end{lemma}

\noindent
\textbf{Proof:}
We will use the following notation $a=(a_1,\dots,a_n)$ for the components of function $a$.

 $N_k$ is a finite sum of the  terms of the following form $g(t,a)$'s given by (\ref{eq:terms-in-N}).
Observe that
\begin{equation}
  (D^{\alpha} a_p)_{k}=  i^{|\alpha|} k_{\alpha_1} \cdots k_{\alpha_{|\alpha|}} a_{p,k}. \label{eq:exp-D-alpha-ap}
\end{equation}

The partial derivative of $g(t,a)$ with respect to $a_{p,j}$
\begin{equation}
  \frac{\partial g(t,a)}{\partial a_{p,j}}= w(t)  \sum_{j_1+j_2 + \dots + j_l=k} \sum_{e=1}^l \delta_e(D^{\alpha_1}a_{p_1})_{j_1}  (D^{\alpha_2}a_{p_2})_{j_2} \cdots (D^{\alpha_l}a_{p_l})_{j_{l}} / a_{p,j},
\end{equation}
where $\delta_e=1$, when the $e$-th factor is $(D^\alpha_p)_j$ for some $\alpha$ and $0$, otherwise.
The $e$-th term in the above sum is $(D^{\alpha_1}a_{p_1})_{j_1}  (D^{\alpha_2}a_{p_2})_{j_2} \cdots  \frac{\partial  (D^{\alpha_e}a_{p_e})_{j_e} }{\partial a_{p,j}} \cdots (D^{\alpha_l}a_{p_l})_{j_{l}}$.

Taking into account the symmetry of $j$'s, equation (\ref{eq:exp-D-alpha-ap}) and our assumption about $|a_j|$'s we obtain the following upper bound for $\left|\frac{\partial g(t,a)}{\partial a_{p,j}}\right|$.
\begin{eqnarray*}
  \left|\frac{\partial g(t,a)}{\partial a_{p,j}}\right| \leq W l \sum_{j_2 + \dots +j_l=k-j} |j|^r \cdot |j_2|^r|a_{j_2}| \cdot \dots \cdot |j_l|^r |a_{j_l}| \leq \\
  W l |j|^r C^{l-1} \sum_{j_2 + \dots +j_l=k-j} \frac{1}{|j_2|^{s-r}} \cdot \cdots \cdot  \frac{1}{|j_l|^{s-r}}
\end{eqnarray*}

Therefore by Lemma~\ref{lem:convolution}  we obtain
\begin{eqnarray*}
   \left|\frac{\partial g(t,a)}{\partial a_{p,j}}\right| \leq W l  \frac{|j|^r C^{l-1} C_2^{l-2}(d,s-r)}{|k-j|^{s-r}}.
\end{eqnarray*}

This concludes the proof.
\qed

For the computation of the logarithmic norms we will need the following lemma.
\begin{lemma}
\label{lem:logNestm}
The same assumptions as in Lemma~\ref{lem:DNestm}. Assume additionally that $s > s_0'=d+2r$. Then
\begin{equation}
  \sum_j  \left| \frac{\partial N_k}{\partial a_j}(t,a)\right| +  \sum_j \left|\frac{\partial N_j}{\partial a_k}(t,a)\right| \leq |k|^r G
\end{equation}
Moreover, the above bound holds if we replace  $\left| \frac{\partial N_k}{\partial a_j}(t,a)\right|$  and  $\sum_j \left|\frac{\partial N_j}{\partial a_k}(t,a)\right|$ by absolute values of the terms defining $\frac{\partial N_k}{\partial a_j}(t,a)$ and $\left|\frac{\partial N_j}{\partial a_k}(t,a)\right|$.
\end{lemma}
\noindent
\textbf{Proof:}
From Lemma~\ref{lem:DNestm} it follows that  for some constant $G=G(d,C,s-r)$ holds
\begin{eqnarray*}
  \sum_j \left|\frac{\partial N_j}{\partial a_k}\right| \leq |k|^r \sum_j \frac{\D}{|k-j|^{s-r}} \leq |k|^r G_1.
\end{eqnarray*}
For the other sum we reason as follows. From Lemma~\ref{lem:DNestm}
\begin{eqnarray*}
  \left| \frac{\partial N_k}{\partial a_j}\right|  \leq \frac{\D |j|^r}{|k-j|^{s-r}} \leq  \frac{\D}{|k-j|^{s-2r}} \left(\frac{|j|}{|k-j|} \right)^r \leq \\
  \frac{\D}{|k-j|^{s-2r}} \left(\frac{|k|+|k-j|}{|k-j|} \right)^r= \frac{\D}{|k-j|^{s-2r}} \left(\frac{|k|}{|k-j|} +1 \right)^r =\\
   \sum_{p=0}^r \binom{r}{p}\frac{\D}{|k-j|^{s-2r}}\left(\frac{|k|}{|k-j|} \right)^p \leq |k|^r  \sum_{p=0}^r \binom{r}{p}\frac{\D}{|k-j|^{s-2r+p}}.
\end{eqnarray*}
Now we finish as with the first sum.
\qed

The next lemma shows the logarithmic norm in suitable neighborhood of the origin is negative, if all eigenvalues $\lambda_k$ are negative.

\begin{lemma}
\label{lem:logNormNegCloseToZero}
Consider (\ref{eq:fugenpde}). Assume that
conditions (\ref{eq:lambdak}), (\ref{eq:lambdak2}) and
(\ref{eq:p>r}) hold. Assume that for all $k$ holds $\lambda_k <0$ i.e. $K_- <0$
and  that $N(t,a,Da,\dots,D^ra)$ is a polynomial in variables $a,Da,\dots,D^ra$ with time dependent and bounded coefficients (i.e. there exists a constant
$C_N$, such that for all $t \in \mathbb{R}$ the coefficients in $N$ are less than or equal to $C_N$). Assume that $N$ does not contain constant terms or degree one terms.

 Let $s>s'_0=d+2r$, $C>0$, $E>0$. We set
\begin{equation}
  W(E,C,s)=\left\{ \{a_k\} \ | \  |a_k| \leq E, \ |a_k| \leq \frac{C}{|k|^s} \right\}.
\end{equation}

Then for any $s > s_0'$ and $C>0$, there exists $E_0$, such that for any $E < E_0$ and $t \in \mathbb{R}$ holds
\begin{equation}
  \mu(D_a F(t,a),\mathbb{R} \times W(E,C,s)) < 0.
\end{equation}
\end{lemma}
\noindent \textbf{Proof:}
Observe first that the perturbation term $V_k(t)$ does not influence  the logarithmic norm of $D_a F(t,a)$.

Let us fix $s>s_0$ and $C>0$.

We use (\ref{eq:log-norm-eucl}), Lemma~\ref{lem:logNestm} and the Gershgorin  theorem \cite{G} to bound the eigenvalues of
 $Sym(D_zF(t,z))=(D_z F(t,z)+D_zF(t,z)^t)/2$. We have for any $E$ and $z \in W(E,C,s)$
\begin{equation*}
  \mbox{Sp}(Sym(DF_z(t,z))) \subset \bigcup_{k} \overline{B}\left(\lambda_k, G |k|^r \right)
\end{equation*}
for some constant $G$ independent from $E$.
 From (\ref{eq:lambdak},\ref{eq:lambdak2},\ref{eq:p>r}) it follows that there exists $N>0$, such that
 \begin{equation}
    \bigcup_{|k| > N} \overline{B}\left(\lambda_k, G |k|^r \right) \cap \mathbb{R}  \subset \mathbb{R}_-.
 \end{equation}

 Now we will argue that if we take $E$ small enough, then also the other part (with $|k| < N$) will be negative.

 We want to show that for $|k|\leq N$ holds
 \begin{equation}
   \lambda_k +   \sum_j  \left| \frac{\partial N_k}{\partial a_j}(t,a)\right| +  \sum_j \left|\frac{\partial N_j}{\partial a_k}(t,a)\right| < 0.
 \end{equation}

Observe that we have only a finite number of inequalities to satisfy and $\lambda_k < 0$. Hence it is enough to show that for any $k$ and $\epsilon >0$  by taking $E$ sufficiently small have
\begin{equation}
  S= \sum_j  \left| \frac{\partial N_k}{\partial a_j}(t,a)\right| +  \sum_j \left|\frac{\partial N_j}{\partial a_k}(t,a)\right| < \epsilon.
\end{equation}
From Lemma~\ref{lem:logNestm} it follows that the sums in $S$ are uniformly converging with respect to $(t,a) \in \mathbb{R} \times W(C,C,s)$.
 and the particular terms contributing to $\frac{\partial N_k}{\partial a_j}(t,a)$ are estimated by a uniformly absolutely converging series.

 Therefore we can bound $S$ as follows
 \begin{equation}
   S\leq S'= \sum_{l \in \mathbb{N}} t_l
 \end{equation}
 where $t_l$ are  upper bounds for  absolute values of any term entering into $\frac{\partial N_k}{\partial a_j}$ or  $\frac{\partial N_j}{\partial a_k}$.
 Observe that since all terms in $N$ were at least of degree two, each of $t_l$ contains as a factor some $a_u$ for some $u \in \mathbb{Z}^d$.

 Therefore there exist $N_1$ such that
 \begin{equation}
   \sum_{l\geq N_1}  t_l  < \epsilon/2.
 \end{equation}

 Now observe that by taking $E$ small enough (remember that $t_l$ are at least degree $1$ in $a$)
  \begin{equation}
   \sum_{l < N_1}  t_l  < \epsilon/2.
 \end{equation}

 Hence $S<\epsilon$ for $E$ small enough.

\qed

\subsection{Existence of self-consistent bounds and solutions for short time step}
The main result in this section is Theorem~\ref{thm:selfexists},
which states that equation (\ref{eq:fugenpde}) satisfying
conditions (\ref{eq:lambdak}), (\ref{eq:lambdak2}), (\ref{eq:p>r})
has solutions within self-consistent bounds for a sufficiently
short time.

\begin{theorem}
\label{thm:existsCi} Consider (\ref{eq:fugenpde}). Assume that
conditions (\ref{eq:lambdak}), (\ref{eq:lambdak2}) and
(\ref{eq:p>r}) hold. Let  $s_0=p+d+1$ and $m \in \mathbb{R}$.

Consider compact set $W\subset X_m$ and a sequence of compact sets
$B_k \subset H_k$ for $|k| > m$, such that there exist $s \geq
s_0$ and $C \in \mathbb{R}$ and  the following condition is
satisfied
\begin{equation}
  |B_{k}| \leq \frac{C}{|k|^{s}}, \qquad |k| >m,\: k \in I.
  \label{eq:polbd}
\end{equation}
Then $W \oplus \Pi_{k \in I, |k|>m }B_k$ satisfies conditions
C2,C3.
\end{theorem}
An easy proof is left as an exercise for the reader.

\begin{theorem}
\label{thm:selfexists} Consider (\ref{eq:fugenpde}). Assume that
conditions (\ref{eq:lambdak}), (\ref{eq:lambdak2}) and
(\ref{eq:p>r}) hold. Let $s_0=p+d+1$

Let  $Z \oplus T_0$ form self-consistent bounds for
(\ref{eq:fugenpde}) for the time interval $[t_0,t_0+h_0]$, such that for some $C_0$ and $s \geq s_0$
it holds that
\begin{equation}
  |T_{0,k}| \leq \frac{C_0}{|k|^{s}}, \qquad |k| >m,\: k \in I, \:
  s> s_0.
\end{equation}

Then there exist $0 < h \leq h_0$, $W \oplus T_1$ -- self-consistent bounds
for (\ref{eq:fugenpde}) over time interval $[t_0,t_0+h]$ and $L >0$, such that for all $l >L$ and
$u \in  P_l(Z \oplus T_0)$
\begin{equation}
  \varphi^l(t_0,[0,h],u) \subset P_l(W \oplus T_1).
\end{equation}
and
\begin{equation}
  |T_{1,k}| \leq \frac{C_1}{|k|^{s}}, \qquad |k|>m, k \in I.
\end{equation}

Moreover, the set  $W \oplus T_1$ can be chosen to be convex.

There exists $M$, such that for $u \in W \oplus T_0$ holds
\begin{equation}
   |\varphi^l_k(t_0,(0,h],u)| < \frac{C_0}{|k|^s} , \quad |k| > M. \label{eq:entersfarmodes}
\end{equation}

\end{theorem}
\noindent
\textbf{Proof:}
 Let $W \subset X_m$, be  a compact set, such that
$Z \subset \inte_{X_m} W$.

By  increasing $C_0$, if necessary, we can assume that
\begin{equation}
  |u_k| \leq \frac{C_0}{|k|^{s}}, \qquad \mbox{for all  $u \in W \oplus T_0$ and $k \in
  I$}.
\end{equation}

We set $C_1=2C_0$ and define the tail $T_1$  by
\begin{equation}
   T_1=\Pi_{|k|>m, k \in I} \overline{B}\left(0,\frac{C_1}{|k|^s}\right).
\end{equation}

From Lemma~\ref{lem:Dgen} applied to the set $\{u \ | \ |u_k| \leq \frac{C_1}{|k|^s}  \}$ over the time interval $[t_0,t_0+h_0]$ it follows that there exists
$\D$, depending on $C_1$, $s$ and bounds of coefficients in $N$ such that
\begin{equation}
  |N_k(t,u)| < \frac{\D}{|k|^{s - r}}, \quad \mbox{for all $(t,u)$, such that $|u_k| \leq \frac{C_1}{|k|^s}$ and $t \in [t_0,t_0+h_0]$.}  \label{eq:proof-estm-Nk}
\end{equation}

From our assumption about smoothness of $V$ it follows that for any $s >0$ and any compact time interval $[t_1,t_2]$ there exists a constant $C=C(V,[t_1,t_2],s)$, such that
\begin{equation}
  |V_k(t)| \leq \frac{C}{|k|^s}, \quad k \in \mathbb{Z}^d, t \in [t_1,t_2].  \label{eq:Vk-gen-decay-estm}
\end{equation}

Let us take $(t,u) \in [t_0,t_0+h_0] \times W \oplus T_1$ and such that $|u_{k_0}|=\frac{C_1}{|k_0|^{s}}$ holds for
some $|k_0| > K_-$. Then from (\ref{eq:proof-estm-Nk}), (\ref{eq:Vk-gen-decay-estm}) and (\ref{eq:lambdak},\ref{eq:lambdak2}) it follows that
\begin{eqnarray*}
  \frac{1}{2}\frac{d}{dt}(u_{k_0}|u_{k_0})(t) < -\beta_0 |k_0|^p |u_{k_0}|^2 + |u_{k_0}|\cdot |N_{k_0}(t,u)| + |u_{k_0}| \cdot | V_{k_0}(t)| \leq   \\
   \left(- \beta_0 C_1 |k_0|^{p-s} + \D |k_0|^{r -s} + C(V,[t_0,t_0+h_0],s) |k_0|^{-s}\right) |u_{k_0}|,
 \end{eqnarray*}
hence
\begin{equation}
  \frac{d|u_{k_0}|^2}{dt}< 0, \qquad |k_0|> L,  \label{eq:farentry}
\end{equation}
for $L$ sufficiently large.

This means that any solution $u(t)$ of  (\ref{eq:fugenpde}) (or its Galerkin projection) with $u(t_0) \in W \oplus T_1$ can leave the set  $W \oplus T_1$ (in the time less than $t_0+h_0$) only through the boundary of $P_{L} (W \oplus T_1)$.

Consider now the differential inclusion
\begin{equation}
  u' \in P_LF(t,u) + \Delta, \qquad u \in X_L, \Delta \subset X_L  \label{eq:galediffincl}
\end{equation}
where the set $\Delta$ represents the Galerkin projection errors
on $W \oplus T_1$ for the time interval $[t_0,t_0+h_0]$ and is given by
\begin{equation}
  \Delta=\{ P_LF(t,u) - P_L F(t,P_Lu) \: | \: (t,u) \in  [t_0,t_0+h_0] \times W \oplus T_1  \}.
\end{equation}
As it was mentioned in the introduction, by a solution of
differential inclusion (\ref{eq:galediffincl}) we will understand
any $C^1$ function $u:[t_0,t_0+t_m] \to X_L$ satisfying condition
(\ref{eq:galediffincl}).

Now we will show that there exists $0 < h \leq h_0$, such that if
$u:[t_0,t_0+t_m]\to X_L$, where $t_m \leq h$, is a  solution of
(\ref{eq:galediffincl}) and $u(t_0) \in P_L(Z \oplus T(0)))$, then
\begin{equation}
u(t_0+t)  \in \inte_{X_L} P_L (W \oplus T_1), \qquad t \in [0,h].
\label{eq:dinclencl}
\end{equation}
Namely, it is enough to take $h>0$, $h \leq h_0$ satisfying the following
condition
\begin{equation}
   h \cdot \left(\max_{(t,u) \in [t_0,t_0+h_0] \times W \oplus T_1  } |P_LF(t,u)| + \max_{\delta \in \Delta}
   |\delta|\right) < \dist(P_L(Z \oplus T_0),\partial_{X_L} P_L(W\oplus
   T_1)).  \label{eq:ssb-h-encl}
\end{equation}
We  prove next that with such $h$ condition (\ref{eq:dinclencl}) is satisfied.

Let $l > L$ and let $u:[t_0,t_0+t_1) \to X_l$ be a solution of
\begin{equation}
  u'(t)=P_lF(t,u), \qquad u(t_0) = u_0 \in P_l(X \oplus T_0).
\end{equation}
By changing the vector field in the complement of  $P_l (W \oplus
T_1)$ we can assume that $t_1=\infty$.

Let
\begin{equation}
  t_m=\sup\{ t \geq 0 \: | \: t \leq h, \: u([t_0,t_0+t]) \subset P_l (W \oplus
T_1)\}.
\end{equation}
It is enough to prove that $t_m = h$.

Obviously $t_m \geq 0$.  We will do the proof by the contradiction.
Assume that $t_m < h$.

Observe that for $t \in [0,t_m]$ $P_Lu(t_0+t)$ is a solution of
(\ref{eq:galediffincl}), hence from (\ref{eq:ssb-h-encl}) we obtain
\begin{equation}
  P_L u([t_0,t_0+t_m])  \subset \inte_{X_L} P_L (W \oplus T_1).
\end{equation}
From (\ref{eq:farentry}) it follows immediately that
\begin{equation}
  Q_L u([t_0,t_0+t_m]) \subset \inte_{Y_l} P_lQ_L(W \oplus T_1).
\end{equation}
Hence
\begin{equation}
   u(t_0+t_m) \in \inte_{X_l} P_l(W \oplus T_1).
\end{equation}
From the above condition and the continuity of $u$ it follows that for
some $\eta >0$ holds
\begin{equation}
   u(t_0+t_m+t') \in \inte_{X_l} P_l(W \oplus T_1), \qquad t'\in [0,\eta].
\end{equation}
But this contradicts the definition of $t_m$. Therefore $t_m=t$.

To establish (\ref{eq:entersfarmodes}) observe that if $(t,u) \in [t_0,t_0+h_0] \times W \oplus T_1$ is such that $\frac{C_0}{|k_0|^{s}} \leq  |u_{k_0}|$ holds for some $|k_0| > K_-$ then from (\ref{eq:proof-estm-Nk}), (\ref{eq:Vk-gen-decay-estm}) and (\ref{eq:lambdak},\ref{eq:lambdak2}) we obtain
\begin{eqnarray*}
  \frac{1}{2}\frac{d}{dt}(u_{k_0}|u_{k_0})(t) < -\beta_0 |k_0|^p |u_{k_0}|^2 + |u_{k_0}|\cdot |N_{k_0}(t,u)|  + |u_{k_0}|\cdot |V_k(t)| \leq  \\
   \left(- \beta_0 C_0 |k_0|^{p-s} + \D |k_0|^{r -s} + C(V,[t_1,t_2],s) |k_0|^{-s}\right) |u_{k_0}|.
\end{eqnarray*}
Hence
\begin{equation}
  \frac{d|u_{k_0}|^2}{dt}< 0, \qquad |k_0|> M,
\end{equation}
for $M$ sufficiently large.

This means that  for any $u(t_0) \in W \oplus T_0$ and for $|k|>M$ if  $|u_k(t_0)| \leq \frac{C_0}{|k|^s}$, then we will have $|u_k(t_0+[0,h])| < \frac{C_0}{|k|^s}$.

 \qed

%% file: rapid-osc-infdim.tex
\section{Lemma on rapid oscillation in the context of self-consistent bounds}
\label{sec:lem-rapid-osc}
Our goal is to generalize Lemma~\ref{lem:prob1-estm-norm}  to dissipative PDEs in the context of self-consistent bounds.

Let us consider two problems of the class defined in Section~\ref{subsec:problem}
\begin{equation}
  \frac{dy}{dt} = F(y):= L y + N(y)  \label{eq:ssb-auto-prob1-loc}
\end{equation}
and its non-autonomous perturbation
\begin{equation}
  \frac{du}{dt}= \tilde{F}(t,u):=L u + N(u) + \tilde{N}(t,u) + V(t) \label{eq:ssb-prob1-loc}
\end{equation}
where
\begin{itemize}
\item $L$ satisfies assumptions from Section~\ref{subsec:problem}.
\item forcing term $V(t)$ is of the following form
\begin{equation}
V(t)=\sum_{k \in \mathbb{Z}^d} g_k(\omega_k t) v_k(t) e_k ,  \label{eq:Rosc-pde-v}
\end{equation}
where $g_k:\mathbb{R} \to \mathbb{R}$ are continuous and $v_k:\mathbb{R} \to \mathbb{R}$ are $C^1$ for every $k \in \mathbb{Z}^d$.
This choice of the forcing term is quite special, but it is enough for our purposes.

\item  $N(u)$ is a real polynomial of $u,Du,\dots,D^ru$
\item $\tilde{N}(t,u)=\tilde{N}(t,u,Du,\dots,D^ru)$ is a polynomial in variables $u,\dots,D^ru$ with
time dependent coefficients,
\begin{eqnarray}
  \tilde{N}(t,u)=\sum_{\sigma=1}^{\sigma_m} \tilde{g}_\sigma(\tilde{\omega}_\sigma t) \tilde{v}_\sigma(u)
\end{eqnarray}
with $\tilde{v}_\sigma(u)$ being polynomials in $u,\dots,D^ru$ and $ \tilde{g}_\sigma$ are bounded continuous functions. Observe that from the point of view of self-consistent bounds these are nice
functions.
\end{itemize}

Let us remind the reader that all considerations are in some Hilbert space $H$ and $\|\cdot\|$ denotes norm in that space, which in the case of coordinates expressed in some Hilbert base $\{e_k\}$ in $H$, is 
\begin{equation*}
  \left\|\sum a_k e_k\right\|=\sqrt{\sum |a_k|^2}.
\end{equation*}

Let  $G_{\eta}'(t)=g_{\eta}(t)$ for $\eta \in \mathbb{Z}^d$ and  $\tilde{G}'_\sigma=\tilde{g}_\sigma$ for $\sigma=1,\dots,\sigma_m$ .

We assume that there exist constants $A_V$, $B_V$ and $s_V$ such that for all $t \in \mathbb{R}$ and $k \in \mathbb{Z}^d$ holds
\begin{eqnarray}
  s_V &>& d+p+r+1, \\
  |v_k(t)| &\leq& \frac{A_V}{|k|^{s_V}},  \label{eq:v_k-bound} \\
  |v_k'(t)| &\leq& \frac{B_V}{|k|^{s_V}}.
\end{eqnarray}

We set (in order to apply later Lemma~\ref{lem:prob1-estm-norm})
\begin{eqnarray}
  C(v_k) &=& \frac{A_V}{|k|^{s_V}}, \label{eq:cv-k} \\
  C\left(\frac{\partial v_k}{\partial t}\right) &=& \frac{B_V}{|k|^{s_V}}. \label{eq:cvdt-k}
\end{eqnarray}

We assume that there exist constants $C(\cdot) \in \mathbb{R}$ such that
\begin{eqnarray}
  \sup_{t \in \mathbb{R}}|\tilde{g}_\sigma(t)| &=& C(\tilde{g}_\sigma),  \quad \sigma=1,\dots,\sigma_m, \\
  \sup_{t \in \mathbb{R}}|g_{\eta}(t)| &=& C(g_{\eta}), \quad \forall \eta \in \mathbb{Z}^d,  \\
  \sup_{t \in \mathbb{R}}|G_{\eta}(t)| &=& C(G_{\eta}), \quad \forall \eta \in \mathbb{Z}^d, \\
    \sup_{t \in \mathbb{R}}|\tilde{G}_{\sigma}(t)| &=& C(\tilde{G}_{\sigma}), \quad \sigma=1,\dots, \sigma_m.
\end{eqnarray}

Assume that
\begin{eqnarray}
     \sum_{\eta \in \mathbb{Z}^d} C(g_\eta) C(v_\eta) &<& \infty, \label{eq:Cgsmalletaveta} \\
      \sum_{\eta \in \mathbb{Z}^d} G(G_\eta)  \frac{1}{|\eta|^{s_V-p}}  &<& \infty. \label{eq:CGetaveta}
\end{eqnarray}

Observe that from (\ref{eq:Cgsmalletaveta}) it follows that
\begin{equation}
  \sup_{ t \in \mathbb{R}}\|V(t)\| \leq \sum_{\eta \in \mathbb{Z}^d} C(g_\eta) C(v_\eta) < \infty.
\end{equation}

All the notations and assumptions listed above are assumed for all the lemmas in this subsection.

The lemma below states that we can have common self-consistent bounds for both problems (\ref{eq:ssb-auto-prob1-loc}) and (\ref{eq:ssb-prob1-loc})
for the same time step.

\begin{lemma}
\label{lem:commonselfexists}
 Let $s_0=p+d+1$.
Let  $Z \oplus T_0$ form self-consistent bounds for (\ref{eq:ssb-auto-prob1-loc}) and
(\ref{eq:ssb-prob1-loc})  over time interval $[t_0,t_0+h_0]$, valid for arbitrary $\{\omega_\eta\}_{\eta \in \mathbb{Z}^d}$ and $\{\tilde{\omega}_\sigma\}_{\sigma=1,\dots,\sigma_m}$ , such that for some $C_0$ and $s \geq s_0$
\begin{equation}
  |T_{0,k}| \leq \frac{C_0}{|k|^{s}}, \qquad |k| >m,\: k \in I, \:
  s_V \geq s> s_0. \label{eq:T0-estm}
\end{equation}

Then there exist
\begin{itemize}
\item $0 < h \leq h_0$,  $d_0 >0$, $M>0$
\item  $W \oplus T_1$ - self-consistent bounds for (\ref{eq:ssb-auto-prob1-loc}) and
for (\ref{eq:ssb-prob1-loc})  over time interval $[t_0,t_0+h]$, valid for arbitrary $\{\omega_\eta\}_{\eta \in \mathbb{Z}^d}$ and $\{\tilde{\omega}_\sigma\}_{\sigma=1,\dots,\sigma_m}$,
\end{itemize}
such that for all $l >d_0$ and
$u \in  P_l(Z \oplus T_0)$   holds
\begin{eqnarray}
  |\varphi^l_k(t_0,(0,h],u)| &<& \frac{C_0}{|k|^s} , \quad |k| > M. \label{eq:twosys-entersfarmodes} \\
  \varphi^l(t_0,[0,h],u) &\subset& W \oplus T_1,
\end{eqnarray}
for both of the semiprocesses (\ref{eq:ssb-auto-prob1-loc}) and (\ref{eq:ssb-prob1-loc}) (denoted by the same symbol $\varphi^l$), where
\begin{equation}
  |T_{1,k}| \leq \frac{C_1}{|k|^{s}}, \qquad |k|>m, k \in I.
\end{equation}
Moreover, the set  $W \oplus T_1$ can be chosen to be convex.
\end{lemma}
\noindent
\textbf{Proof:}
The proof is the same as that of Theorem~\ref{thm:selfexists}. The only difference is that the constant $\D$, which was considered there
will now the following structure $\D=\D_1 + \D_2 + A_V$, where $\D_1$, $\D_2$ are obtained from Lemma~\ref{lem:Dgen} for polynomials $N$, $\tilde{N}$, respectively.  Such constant $\D$ is good for both equations (\ref{eq:ssb-auto-prob1-loc}) and (\ref{eq:ssb-prob1-loc}).

\qed

The next lemma shows that under our assumptions all the quantities required by Lemma~\ref{lem:prob1-estm-norm} are finite.
\begin{lemma}
\label{lem:ssb-estm-bnds}
  Let $s_0=d+p+r+1$. Let  $Z \oplus T_0$ be self-consistent bounds for
(\ref{eq:ssb-auto-prob1-loc}) and (\ref{eq:ssb-prob1-loc}) on the whole real line $\mathbb{R}$ valid for arbitrary $\{\omega_\eta\}_{\eta \in \mathbb{Z}^d}$ and $\{\tilde{\omega}_\sigma\}_{\sigma=1,\dots,\sigma_m}$,  such that for some $C_0$ and $s \geq s_0$, $s \leq s_V$
holds
\begin{equation}
  |T_{0,k}| \leq \frac{C_0}{|k|^{s}}, \qquad |k| >m,\: k \in I.
\end{equation}

Then there exist
\begin{itemize}
\item $0 < h \leq h_0$, $d_0 >0$, $C_1 >0$ and $s>0$
\item $W \oplus T_1$ -- self-consistent bounds
for problems (\ref{eq:ssb-auto-prob1-loc}) and (\ref{eq:ssb-prob1-loc})  over the time interval $[t_0,t_0+h]$ valid for arbitrary $t_0 \in \mathbb{R}$, $\{\omega_\eta\}_{\eta \in \mathbb{Z}^d}$ and $\{\tilde{\omega}_\sigma\}_{\sigma=1,\dots,\sigma_m}$,
\end{itemize}
 such that for all $n >d_0$ and
$u \in  P_{n}(Z \oplus T_0)$
\begin{eqnarray*}
  \varphi^{n}(t_0,[0,h],u) &\subset& W \oplus T_1, \\
  |T_{1,k}| &\leq& \frac{C_1}{|k|^{s}}, \qquad |k|>m, k \in I.
\end{eqnarray*}
where $\varphi^{n}$ stands for the semiprocess for $n$-th Galerkin projection of (\ref{eq:ssb-auto-prob1-loc}) or (\ref{eq:ssb-prob1-loc}).

The set  $W \oplus T_1$ can be chosen to be convex, $W \oplus T_1$ and $h$ do not depend on $t_0$, $\omega_k$'s and $\tilde{\omega}_j$'s.

Moreover, the following quantities are finite and do not depend on $t_0$ (we use symbol $\mathcal{P}_i$ to denote the Galerkin projection of the map, where
we restrict both the range and the domain, as defined in Def.~\ref{defn:Faddmissible})

\begin{eqnarray}
  \sup_{ i >0}\sup_{\substack{z \in P_i (W \oplus T_1)\\t \in \mathbb{R}}} \mu(D_z (\mathcal{P}_i F)(z))&=&l \nonumber  \\
    \sup_{ i >0} \sup_{\substack{z \in P_i (W \oplus T_1)\\t \in \mathbb{R}}} \|\mathcal{P}_i\tilde{v}_\sigma(z) \| &=& C\left(\tilde{v}_\sigma\right), \quad \sigma=1,\dots,\sigma_m, \nonumber \\
      \sup_{ i >0}\sup_{\substack{z \in P_i (W \oplus T_1)\\t \in \mathbb{R}}} \left\|(D_z (\mathcal{P}_i \tilde{v}_\sigma)(z)) (\mathcal{P}_i\tilde{F})(t,z) \right\| &=&  C(D_z \tilde{v}_\sigma \tilde{F}),
 \quad \sigma=1,\dots,\sigma_m \nonumber \\
  \sup_{ i >0} \sup_{\substack{z,z_1 \in P_i (W \oplus T_1)\\t \in \mathbb{R}}} \|(D_z \mathcal{P}_i F(z)) \mathcal{P}_i\tilde{v}_\sigma(z_1) \| &=& C\left(D_z F\tilde{v}_\sigma\right), \quad \sigma=1,\dots,\sigma_m, \nonumber \\
  \sup_{ i >0} \sup_{\substack{z \in P_i (W \oplus T_1)\\t \in \mathbb{R}}} \|(D_z \mathcal{P}_i F(z)) v_j(t) \| &=& C\left(D_z F v_j\right) \leq  \frac{A_V}{|j|^{s_V - p}}  \left(\beta_1 + \D  S_d (s-r) \right) \label{eq:DzFvjLem} \quad j \in \mathbb{Z}^d, \\
   \sum_{\sigma=1}^{\sigma_m} C(\tilde{G}_\sigma) C(\tilde{v}_\sigma) &<& \infty, \nonumber \\
    \sum_{\eta \in \mathbb{Z}^d} C(G_\eta) C(v_\eta) &<& \infty, \nonumber \\
     \sum_{\eta \in \mathbb{Z}^d} C(G_\eta) C\left(\frac{\partial v_\eta}{\partial t}\right) &<& \infty, \nonumber  \\
   \sum_{\sigma=1}^{\sigma_m} C(\tilde{G}_\sigma) C\left(D_z F\tilde{v}_\sigma\right)  &<& \infty,\nonumber \\
     \sum_{\eta \in \mathbb{Z}^d} C(G_\sigma) C\left(D_z F v_\eta\right)  &<& \infty,\nonumber \\
     \sum_{\sigma=1}^{\sigma_m} C(\tilde{G}_\sigma) C(D_z \tilde{v}_\sigma \tilde{F}) &<& \infty. \nonumber
\end{eqnarray}

\end{lemma}

\noindent
\textbf{Proof:}

In order to make our notation more transparent, we will drop the symbols $\mathcal{P}_i$, $P_i$ indicating that we are in fact working with Galerkin projections. This is justified by the fact that the estimates we are developing are independent of the projection.

The existence of $h>0$ and $W \oplus T_1 $ follows from Lemma~\ref{lem:commonselfexists}.

In order to estimate $l$ we use (\ref{eq:log-norm-eucl}), Lemma~\ref{lem:logNestm} and the Gershgorin  theorem \cite{G} to bound the eigenvalues of $Sym(D F(z))=(DF(z)+DF(z)^t)/2$ for $z \in W \oplus T_1 $. We have
\begin{equation*}
  \mbox{Sp}(Sym(DF(z))) \subset \bigcup_{k} \overline{B}\left(\lambda_k, G |k|^r \right)
\end{equation*}
for some constant $G$. From (\ref{eq:lambdak},\ref{eq:lambdak2},\ref{eq:p>r}) it follows that $ \mbox{Sp}(Sym(DF(z)))$ is bounded, hence
$l < \infty$.

To obtain $C(\tilde{v}_\sigma)$ observe that from Lemma~\ref{lem:Dgen} we have for $z \in W\oplus T_1$
\begin{equation}
  |\tilde{v}_{\sigma,k}(z)| \leq \frac{\D}{|k|^{s-r}}, \label{eq:estm-v-sigma}
\end{equation}
where $\D=\D(C,s,v_\sigma)$.

From this and since $s>s_0$ it follows that
\begin{equation}
  \|\tilde{v}_\sigma (z)\| \leq C(\tilde{v}_\sigma)\leq \sum_{k} \frac{\D}{|k|^{s-r}} = \D S_d(s-r) < \infty. \label{eq:estm-tildev-sigma}
\end{equation}

Now we estimate $ C\left(D_z \tilde{v}_\sigma \tilde{F}\right)$. From Lemma~\ref{lem:DNestm} it follows that
for $z \in W \oplus T_1$ holds
\begin{equation}
  \left|\frac{\partial \tilde{v}_{\sigma,k}}{\partial z_j} (z) \right| \leq \frac{\D_1 |j|^r}{|k-j|^{s-r}}, \label{eq:estmDtildev}
\end{equation}
where $\D_1$ depends on the polynomial $\tilde{v}_\sigma$, $d$, $r$, $W$, $C$ and $s$.
From our assumptions (\ref{eq:lambdak},\ref{eq:lambdak2},\ref{eq:p>r}), Lemma~\ref{lem:Dgen} and (\ref{eq:v_k-bound}) it follows that
\begin{equation}
  |\tilde{F}_k(t,z)| \leq \frac{\D_2}{|k|^{s-p}}, \quad t \in \mathbb{R}, z \in W \oplus T_1 \label{eq:estmtildeF}
\end{equation}
for some  $\D_2 \in \mathbb{R}$.

Therefore from (\ref{eq:estmDtildev}), (\ref{eq:estmtildeF}) and Lemma~\ref{lem:convolution} we obtain
\begin{eqnarray*}
   \left|(D_z \tilde{v}_{\sigma}(z) \tilde{F}(t,z))_k \right| \leq \D_1 \D_2 \sum_{j} \frac{|j|^r}{|k-j|^{s-r} |j|^{s-p}}  = \\
   \D_1 \D_2 \sum_{j} \frac{1}{|k-j|^{s-r} |j|^{s-p-r}} < \D_1 \D_2 \sum_{j} \frac{1}{|k-j|^{s-r-p} |j|^{s-p-r}} \leq\\
    \frac{\D_1 \D_2 C_2(d,s-r-p)}{|k|^{s-r-p}}.
\end{eqnarray*}
Hence there exists a constant $C(D_z \tilde{v}_\sigma \tilde{F})$ such that
\begin{equation*}
   \left\|D_z \tilde{v}_{\sigma}(z) \tilde{F}(t,z) \right\| \leq C(D_z \tilde{v}_\sigma \tilde{F}).
\end{equation*}

To obtain an upper bound for $C(D_z F \tilde{v}_\sigma)$ we proceed as follows. We have bound for $|\tilde{v}_\sigma(z)|$ for $z \in W \oplus T_1$ given by
(\ref{eq:estm-tildev-sigma}),
while for $D_z F(t,z)$ from our assumptions about $L$ (\ref{eq:lambdak},\ref{eq:lambdak2},\ref{eq:p>r}) and Lemma~\ref{lem:DNestm} applied to
$N$  it follows that
\begin{equation}
  |(D_z F(t,z))_{kj}| \leq |\lambda_k| \delta_{kj} + \frac{\D |j|^r}{|k-j|^{s-r}} \leq \beta_1 |k|^p  \delta_{kj} + \frac{\D |j|^r}{|k-j|^{s-r}}, \label{eq:DzF}
\end{equation}
where $\delta_{kj}=1$ when $k=j$ and $\delta_{kj}=0$ otherwise. We conclude as in the case of $C(D_z \tilde{v}_{\sigma} \tilde{F})$.

Now we will estimate $C(D_z F v_j)$.
Observe that
\begin{eqnarray*}
  (D_z F v_j)_k=  (D_z F)_{kj}  v_j(t).
\end{eqnarray*}
Therefore from (\ref{eq:DzF}) and (\ref{eq:v_k-bound}) we obtain
\begin{equation}
  |(D_z F v_j)_k| \leq \left(\delta_{kj} \beta_1 |k|^p +   \frac{\D |j|^r}{|k-j|^{s-r}}\right) \frac{A_V}{|j|^{s_V}}.
\end{equation}

Therefore  we obtain
\begin{eqnarray*}
  \|(D_z F v_j)\| \leq A_V \left( \beta_1 |j|^p  \frac{1}{|j|^{s_V}} + \sum_{k}\frac{\D |j|^r}{|k-j|^{s-r}|j|^{s_V}} \right) \leq \\
    A_V \left(\frac{\beta_1}{|j|^{s_V - p}} + \frac{\D }{|j|^{s_V - r}}  \sum_{k}\frac{1 }{|k-j|^{s-r}} \right) \leq \\
      \frac{A_V}{|j|^{s_V - p}}  \left(\beta_1 + \D  \sum_{k}\frac{1}{|k-j|^{s-r}} \right).
\end{eqnarray*}
Hence we have
\begin{equation}
   \|(D_z F v_j)\| \leq C(D_z F v_j)= \frac{A_V}{|j|^{s_V - p}}  \left(\beta_1 + \D  \sum_{k}\frac{1}{|k-j|^{s-r}} \right) = \frac{A_V}{|j|^{s_V - p}}  \left(\beta_1 + \D  S(d,s-r) \right).
 \label{eq:DzFvj}
\end{equation}

The assertions regarding the finiteness of sums $\sum_{\sigma}$ are true, because we add finite number of terms, only.

For the remaining sums over $\eta$ we have from (\ref{eq:cv-k}) and (\ref{eq:CGetaveta})
\begin{eqnarray*}
  \sum_{\eta \in \mathbb{Z}^d} C(G_\eta) C(v_\eta) \leq A_V \sum_{\eta \in \mathbb{Z}^d} C(G_\eta) \frac{1}{|\eta|^{s_V}} < \infty.
\end{eqnarray*}
Analogously we obtain from (\ref{eq:cvdt-k}) and (\ref{eq:CGetaveta})
\begin{eqnarray*}
  \sum_{\eta \in \mathbb{Z}^d} C(G_\eta) C\left(\frac{\partial v_\eta}{\partial t}\right) \leq B_V \sum_{\eta \in \mathbb{Z}^d} C(G_\eta) \frac{1}{|\eta|^{s_V}} < \infty.
\end{eqnarray*}
For the last sum we have from (\ref{eq:DzFvj}) for some constant $\D_3$
\begin{eqnarray*}
  \sum_{\eta \in \mathbb{Z}^d} C(G_\eta) C(D_z F v_\eta) \leq  \sum_{\eta \in \mathbb{Z}^d} \frac{C(G_\eta)  A_V \D_3}{|\eta|^{s_V - p}}
\end{eqnarray*}
which is finite due to (\ref{eq:CGetaveta}).

\qed

The lemma below is the variant of Lemma~\ref{lem:prob1-estm-norm} adapted to the context of dissipative PDEs and self-consistent bounds.

\begin{lemma}
\label{lem:ssb-estm-norm}
The same assumptions as in Lemma~\ref{lem:ssb-estm-bnds}.

Let $h$ and $W \oplus T_1$ be as in the assertion of Lemma~\ref{lem:ssb-estm-bnds}.

Let  $y:[t_0,t_0+h] \to \mathbb{R}^n$ and  $x:[t_0,t_0 + h] \to \mathbb{R}^n$  be solutions to (\ref{eq:ssb-auto-prob1-loc}) and (\ref{eq:ssb-prob1-loc}),
respectively, such that $x(t_0)=y(t_0)=z_0$, which are contained in  $W\oplus T_1$.

Then for $t \in [0,h]$ it holds that
\begin{equation}
 \|x(t_0+t)-y(t_0+t)\| \leq \sum_{j \in \mathbb{Z}^d} \frac{1}{|\omega_{j}|} b_{j}(t) + \sum_{\sigma=1}^{\sigma_m} \frac{1}{|\tilde{\omega}_{\sigma}|} \tilde{b}_{\sigma}(t)
 \end{equation}
where  continuous  functions $b_j,\tilde{b}_\sigma:[0,h] \to \mathbb{R}_+$ depend on the constants $C(\cdots)$ (see (\ref{eq:bk-expression})) defined in  Lemma~\ref{lem:prob1-estm-norm}
and the set $W \oplus T_1$, but do not depend on $t_0$.

The sum $\sum_{j \in \mathbb{Z}^d} \frac{1}{|\omega_{j}|} b_{j}(t)$ is convergent, if for all $j \in \mathbb{Z}^d$ $|\omega_{j}| > \epsilon >0$.

\end{lemma}
\noindent
\textbf{Proof:}
We apply Lemma~\ref{lem:prob1-estm-norm} to Galerkin projections, with uniform estimates of various constants $C(\cdot)$, which are finite due to Lemma~\ref{lem:ssb-estm-bnds}, and then we pass to the limit (we take
convergent subsequences).

 Observe that in the present case our rapidly oscillating part consists of two types of terms the ones parameterized by $\sigma$ (there is only a finite number of them) and the other parameterized by $j \in \mathbb{Z}^d$.
 Accordingly the expression for $b_j$'s  in  Lemma~\ref{lem:prob1-estm-norm} is split into two parts, parameterized by $\sigma$ and $j$.

\qed

\subsection{Some definitions}

In the next sections we will be considering the Burgers equation and the Navier-Stokes equations. In order to have the freedom to play with $s_V$ in the context
of the two above lemmas we introduce the following definition.

\begin{definition}
\label{def:uniform-C-infinity}
  Let $f: \mathbb{R} \times \mathbb{T}_d \to \mathbb{R}^{n}$. We will say that $f(t,\cdot) \in C^\infty$  uniformly  with respect to $t\in \mathbb{R}$, if $f(t,\cdot) \in C^\infty$ for all $t \in \mathbb{R}$ and for every $s \in \mathbb{N}$ there exists a constant $C(s) \in \mathbb{R}$, such that
  \begin{equation}
    \sup_{t \in \mathbb{R}, k \in \mathbb{Z}^d} |k|^s |f_k(t)| \leq C(s),
  \end{equation}
  where $f_k(t)$ is the $k$-th Fourier coefficient of $f(t,\cdot)$.
\end{definition}

If $u:\mathbb{T}_d \to \mathbb{R}^n$, then its Fourier coefficients satisfy
 for $k\in\mathbb{Z}$
\begin{equation}
    u_k=\overline{u}_{-k}.
\end{equation}
 This motivates the following definition.
\begin{definition}
\label{def:reality}
In the space of sequences $\{u_k\}_{k \in \mathbb{Z}^d}$,  where
$u_k\in \mathbb{C}^n$, we will say that the sequence $\{u_k\}$
satisfies the reality condition iff
\begin{equation}
   u_{k}=\overline{u}_{-k}, \quad k \in \mathbb{Z}^d.  \label{eq:reality-cond}
\end{equation}
We will denote the set of sequences  satisfying
(\ref{eq:reality-cond}) by $R$. It is  easy to see that $R$ is a
vector space over the field $\mathbb{R}$.
\end{definition}

Formally, space $R$ depend on $d$  and $n$, but it will be always clear from the context, what $R$ we are talking about.

\begin{definition}
\label{def:ZCs}
Let 
\begin{equation}
  Z(C,s)=\left\{a=\{a_k\}_{k \in \mathbb{Z}^d \setminus \{0\}}\ | \   |a_k| \leq \frac{C}{|k|^s}\right\}.
\end{equation}
\end{definition} 

%% file: burgers.tex
\section{Viscous Burgers equation with periodic boundary conditions on the line}
\label{sec:def-burgers}

The Burgers equation was proposed in \cite{B} as a mathematical
model of turbulence. There is a significant number of applications
of the Burgers equation, see e.g. \cite{Wh}. We consider the
initial value problem for viscous Burgers equation on the real line
with periodic boundary conditions and \emph{a non-autonomous forcing $F$}, i.e.
\providecommand{\timeInterval}{[t_0, \infty)}
\providecommand{\forcingTimeInterval}{\mathbb{R}}
\begin{subequations}
    \label{eq:burgers}
    \begin{align}
        &u_t(t,x)+u(t,x)\cdot u_x(t,x)-\nu u_{xx}(t,x)=f(t,x),\quad t\in\timeInterval,\ x\in\mathbb{R},\label{eq:burgers1}\\
        &u(t, x)=u(t, x+2\pi),\quad t\in\timeInterval,\ x\in\mathbb{R},\label{eq:burgers3}\\
        &f(t, x)=f(t, x+2\pi),\quad t\in\forcingTimeInterval,\ x\in\mathbb{R},\label{eq:burgers4}\\
        &u(t_0,x)=\bar{u}(x),\quad t_0\in\mathbb{R},\ x\in\mathbb{R},\label{eq:burgers2}
    \end{align}
\end{subequations}
where $\nu >0$.

We will use the Fourier series to study \eqref{eq:burgers}. Let
\begin{equation}
u(t,x)=\sum_{k \in \mathbb{Z}}u_k(t) \exp(ikx).
\end{equation}

It is straightforward to write the problem \eqref{eq:burgers} in
the Fourier basis. We obtain the following infinite ladder
of equations
\begin{equation}
\label{eq:burgers_infinite1}
 \frac{d u_k}{d t}=-i\frac{k}{2}\sum_{k_1\in\mathbb{Z}}{u_{k_1}\cdot u_{k-k_1}}+
        \lambda_k u_k+f_k(t),\quad t\in\timeInterval,\ k\in\mathbb{Z},
\end{equation}
where
\begin{subequations}
    \label{eq:burgers_infinite2}
    \begin{align}
        &u_k(t_0)=\frac{1}{2\pi}\int_0^{2\pi}{\bar{u}(x)e^{-ikx}}\,dx,\quad k\in\mathbb{Z},\\ 
        &f_k(t)=\frac{1}{2\pi}\int_0^{2\pi}{f(t,x)e^{-ikx}}\,dx,\quad k\in\mathbb{Z},\label{eq:burgers_infinite3}\\
        &\lambda_k=-\nu k^2. \label{eq:burgers_infinite5}
    \end{align}
\end{subequations}

 The reality of $u$ and $f$ implies
that for $k\in\mathbb{Z}$
\begin{equation}
    u_k=\overline{u}_{-k}, \quad  f_k=\overline{f}_{-k},\quad \text{ for }t\in\forcingTimeInterval.
\end{equation}
We see that variables $\{u_k\}_{k\in\mathbb{Z}}$ are not
independent.
Observe that (see Def.~\ref{def:reality})
\begin{equation}
  \{u_k(t)\} \in R, \quad \{f_k(t)\} \in R, \qquad \mbox{for $t \in \mathbb{R}$}.
\end{equation}

We  assume that the initial condition for
\eqref{eq:burgers_infinite1} satisfies
\begin{equation}
  \label{eq:fixedInt}
  \frac{1}{2\pi}\int_{0}^{2\pi}{\bar{u}(x)\,dx}=\alpha,\quad\text{for a fixed }\alpha\in\mathbb{R}.
\end{equation}
We  require additionally that $f_0(t)=0$  for $t\in\forcingTimeInterval$,
and then \eqref{eq:fixedInt} implies that
\begin{equation}
  \label{eq:fixedA0}
  u_0(t)=\alpha, \quad \forall t \geq t_0.
\end{equation}

\begin{definition}
    \label{def:symmetricGalerkinProjection}
    For any given number $m>0$ {\rm the $m$-th Galerkin projection} of \eqref{eq:burgers_infinite1} is
    \begin{equation}
        \label{eq:symmetricGalerkinProjection}
        \frac{d u_k}{d t}=-i\frac{k}{2}\sum_{\substack{|k-k_1|\leq m\\|k_1|\leq m}}{u_{k_1}\cdot u_{k-k_1}}+\lambda_k u_k+f_k(t),\quad t\in\timeInterval,\ |k|\leq m.
    \end{equation}
\end{definition}
Note that the condition \eqref{eq:fixedA0} holds also for all
Galerkin projections \eqref{eq:symmetricGalerkinProjection} as
long as $f_0(t)=0$ for all
$t\in\forcingTimeInterval$. Also observe that   the reality condition
\eqref{eq:reality-cond} is invariant under all Galerkin
projections \eqref{eq:symmetricGalerkinProjection}, i.e. if
$u_k(t_0)=\overline{u}_{-k}(t_0)$, then $u_k(t)=\overline{u}_{-k}(t)$
for all $t > t_0$ if the solution of
(\ref{eq:symmetricGalerkinProjection}) exists up to that time.

\begin{definition}
\label{def:tildeH}
  Let $H$ be the space $l_2(\mathbb{Z},\mathbb{C})$, i.e. $u \in H$ is a sequence $u:\mathbb{Z}\to \mathbb{C}$ such that
  $\sum_{k \in \mathbb{Z}}|u_k|^2 < \infty$  over the coefficient field $\mathbb{R}$.
  The subspace $\widetilde{H}\subset H$ is defined by
  \begin{equation*}
    \widetilde{H}:=\left\{\{u_k\}\in H\colon\text{there exists }0\leq C<\infty\text{ such that }|u_k|\leq\frac{C}{\nmid k\nmid^4}\text{ for }k\in\mathbb{Z}\right\}.
  \end{equation*}

  Let the space $\subspaceH$ be given by
  \begin{equation*}
    \subspaceH:=H\cap R.
  \end{equation*}

  For the subspace with $a_0 = 0$ we introduce
  \begin{equation}
     \subspaceH_0 :=  \subspaceH \cap \{ a_0 = 0\}.
  \end{equation}
\end{definition}

Let us comment on Definition~\ref{def:tildeH}. Despite the fact that we are dealing with complex sequences we use as the coefficient field
the set of real numbers, because the reality condition is not compatible with the complex multiplication.

  The choice of the particular subspace $\subspaceH$ is motivated by the fact that the order
  of decay of coefficients $\{u_k\}\in\subspaceH$ is sufficient for the uniform convergence of $\sum{u_ke^{ikx}}$ and every term appearing in
  \eqref{eq:burgers1}.

\subsection{The effect of the moving coordinate frame}

Let us transform the Burgers equation to a coordinate frame, which is moving with the velocity $c$.
Since the function $u(t,x)$ has the meaning of velocity, this transformation on the function level works as
follows: $u(t,x)$ is transformed into  $a(t,x)$
\begin{equation}
  a(t,x)=u(t,x+ct) -c.
\end{equation}

We have the following easy lemma.
\begin{lemma}
Assume that for $t \in \mathbb{R}$ holds $\int_0^{2\pi}f(t,x)dx=0$ and $\frac{1}{2\pi}\int_0^{2\pi} \bar{u}(x)dx= c$.

Let $u(t,x)$ be a solution  of (\ref{eq:burgers}). Then function
 \begin{equation}
  a(t,x)=u(t,x+ct) -c
\end{equation}
is a solution of (\ref{eq:burgers}) with the forcing term $\tilde{f}(t,x)=f(t,x+ct)$. Moreover,
\begin{eqnarray*}
\int_0^{2\pi}\tilde{f}(t,x)dx&=&0, \quad \forall t \in \mathbb{R}, \\
\int_{0}^{2\pi}a(t,x)dx&=&0, \quad \forall t \in \mathbb{R}.
\end{eqnarray*}
\end{lemma}

\comment{
\noindent
\textbf{Proof:}
We have
\begin{equation}
  u(t,x)=a(t,x-ct)+c.
\end{equation}
Hence
\begin{eqnarray*}
  u_t(t,x)&=&a_t(t,x-ct) - c a_x(t,x-ct) \\
  u_x(t,x)&=&a_x(t,x-ct) \\
   u_{xx}(t,x)&=&a_{xx}(t,x-ct)
\end{eqnarray*}
Therefore from
\begin{equation*}
  u_t(t,x) + u(t,x)u_x(t,x) = \nu u_{xx}(t,x) + f(t,x)
\end{equation*}
we obtain
\begin{eqnarray*}
  a_t(t,x-ct) - c a_x(t,x-ct) + (a(t,x-ct)+c)a_x(t,x-ct) = \nu a_{xx}(t,x-ct) + f(t,x) \\
  a_t(t,x-ct)  + a(t,x-ct)a_x(t,x-ct) = \nu a_{xx}(t,x-ct) + f(t,x) \\
  a_t(t,x)  + a(t,x)a_x(t,x) = \nu a_{xx}(t,x) + f(t,x+ct).
\end{eqnarray*}
Observe that from the space periodicity of $f$ it follows that for any $t$ holds
\begin{equation}
  \int_0^{2\pi}  f(t,x+ct)dx=  \int_0^{2\pi}  f(t,x)dx=0
\end{equation}
We also have for all $t \in \mathbb{R}$
\begin{equation}
  \int_0^{2\pi} a(t,x)dx=\int_0^{2\pi} (u(t,x-ct) - c)dx=\int_0^{2\pi} u(t,x)dx - 2 \pi c= 0.
\end{equation}
\qed
}

\subsection{The action of movement of coordinate frame on the Fourier modes}

Assume that
\begin{equation*}
  f(t,x)=\sum_{k \in \mathbb{Z}} f_k(t) \exp(ikx)
\end{equation*}

Observe that the Fourier expansion of $\tilde{f}(t,x)$ is
\begin{eqnarray*}
  \tilde{f}(t,x)=f(t,x+ct)=\sum_{k \in \mathbb{Z}} f_k(t) \exp(ik(x+ct))= \sum_{k \in \mathbb{Z}} f_k(t)\exp(ikct) \exp(ikx)
\end{eqnarray*}
Hence
\begin{equation}
  \tilde{f}_k(t)= f_k(t)\exp(ikct).
\end{equation}

If $u(t,x)=\sum_{k \in \mathbb{Z}} u_k(t) \exp(ikx)$ with $u_0(t)=c$, then
\begin{eqnarray*}
  a(t,x)= u(t,x+ct) - c= \sum_{k \in \mathbb{Z} \setminus \{0\}} u_k(t) \exp(ik(x+ct))=\\
    \sum_{k \in \mathbb{Z} \setminus \{0\}} (u_k(t)\exp(ikct)) \exp(ikx)
\end{eqnarray*}
Hence
\begin{eqnarray}
  a_0(t)=0, \qquad a_k(t)= u_k(t) \exp(ikct), \quad k \neq 0.  \label{eq:fourier-modes-trans}
\end{eqnarray}

%% file: fast-mov-burgers.tex
\section{ Burgers equation with large average speed}
\label{sec:burgers}
We consider Burgers equation (\ref{eq:burgers}).

Recall that if $f_0 \equiv 0$, then  $u_0(t)=\alpha$ for $t \in \mathbb{R}$. Therefore we can write (\ref{eq:burgers_infinite1}) as follows
\begin{equation}
 u_k'=  \lambda_k u_k -ik \alpha u_k + N_k(u) + f_k(t), \quad k \in \mathbb{Z} \setminus \{0\},
\end{equation}
where
\begin{equation}
 N_k(u)= \frac{-ik}{2} \sum_{k_1 \in \mathbb{Z}\setminus \{0,k\}}u_{k_1}u_{k-k_1}.
\end{equation}

Observe that the transformation of function $u(t,x)$ when passing to the moving coordinate frame (compare (\ref{eq:fourier-modes-trans})) given by
\begin{equation}
a_0=0, \qquad a_k=u_k \exp(-i k\alpha t), \ k \neq 0,
\end{equation}
preserves the reality condition, i.e. if $\overline{u_{-k}}=u_k$, then $\overline{a_{-k}}=a_k$
and
\begin{equation}
 N_k(a(u))=N_k(u).
\end{equation}

Therefore in the moving coordinate frame we obtain the system
\begin{equation}
 a_k'=  \lambda_k a_k + N_k(a) + f_k(t)\exp(i k \alpha t)=\widetilde{F}_k(t,a), \quad k \in \mathbb{Z} \setminus \{0\}. \label{eq:osc-burgers}
\end{equation}
In order to have a reference, we will also write the autonomous problem
\begin{equation}
 a_k'=  \lambda_k a_k + N_k(a)=F_k(a), \quad k \in \mathbb{Z} \setminus \{0\}, \label{eq:osc-nog-burgers}
\end{equation}
and, to be consistent  with the previous sections, we denote the forcing term in new coordinates by
\begin{equation}
  V_k(t)=f_k(t)\exp(i k \alpha t),
\end{equation}
so that in the notation introduced in Sections~\ref{sec:basic-estm} and Section~\ref{sec:lem-rapid-osc} we have
\begin{equation*}
  v_k(t,x)=f_k(t),\qquad g_k(t)=\exp(it),\qquad \omega_k=k\alpha.
\end{equation*}
and we do not have terms $\tilde{N}$, $\tilde{v}_\sigma$ and $\tilde{\omega}$'s.

Let us recall some results regarding the existence of forward invariant absorbing sets for \eqref{eq:osc-burgers}.

\begin{definition}{\rm\cite[Def. 3.1]{Cy}}
\label{def:energy}
\emph{Energy}  of $a \in H$ is given by the formula
\begin{equation}
    \label{eq:energy}
    E(\{a_k\})=\sum_{k\in\mathbb{Z}}{|a_k|^2}.
\end{equation}
Energy of $a \in H_0$ is given by the formula
\begin{equation}
  \label{eq:energyWithoutZero}
  \mathcal{E}(\{a_k\})=\sum_{k\in\mathbb{Z}\setminus\{0\}}{|a_k|^2}.
\end{equation}
\end{definition}

The following result is well known
\begin{lemma}{ \cite[Lemma 3.2]{Cy} }
\label{lem:burgers-energy-decay}
Consider  Galerkin projection of (\ref{eq:osc-burgers}). Then
\begin{equation}
  \frac{d \mathcal{E}(a(t))}{dt} \leq - 2 \nu \mathcal{E}(a(t)) + 2\sqrt{\mathcal{E}(a(t))} \sqrt{\mathcal{E}(f(t))}.
\end{equation}
\end{lemma}

\providecommand{\logarithmicTrapReg}{\mathcal{R}}
\providecommand{\logTrapReg}{\logarithmicTrapReg}



\providecommand{\jednadruga}{\frac{1}{2}}
\providecommand{\Ezero}{\mathcal{E}}
\providecommand{\Fassumption}{
$f_k(t)=\overline{f_{-k}(t)}$, $f_k(t)=0$ for $|k|>J$, and $f_0(t)=0$}

\providecommand{\satisfiesFassumption}{the external force $f$ satisfies: $f(t,\cdot)\in C^\infty$
uniformly  with respect to $t\in \mathbb{R}$, $f_0=0$, $f_k(t)=\overline{f_{-k}(t)}$, and  $\Ef<\infty$}

\providecommand{\energyBound}{\widetilde{\mathcal{E}}}
\providecommand{\Ef}{\sup_{t\in\forcingTimeInterval}\mathcal{E}\left(\left\{f_k(t)\right\}\right)}
\providecommand{\invariantSubspace}{restricted to the invariant
subspace given by $a_k=\overline{a_{-k}}$}
\providecommand{\Dd}{2^{s-\jednadruga}+\frac{2^{s-1}}{\sqrt{2s-1}}}

The two next theorems are about the existence of a trapping region and an absorbing set  for (\ref{eq:osc-burgers}) (see Definitions~\ref{def:trappingReg} and~\ref{def:absorbingSet}). Both these results
are well known, we recall them here in the self-consistent bounds context as given in \cite{Cy}. The results from \cite{Cy} require a little
adaptation, because there  $f$ had only a finite number of nonzero Fourier coefficients. Moreover, now we assume additionally that $f_0(t)\equiv 0$
and $a_0=0$. This is the reason why we can replace $E(a)$ and $E(f)$ by $\mathcal{E}(a)$ and $\mathcal{E}(f)$ in the present paper.

\begin{theorem}{\rm \cite[Thm. 2.8]{CyZ}}
  \label{thm:analyticTrappingRegion}
  Assume that \satisfiesFassumption.
  Let $\{a_k\}_{k\in\mathbb{Z}}\in H$,  $s>0.5$, $E_0=\frac{\Ef}{\nu^2} < \infty$,
  $\energyBound>E_0 $, $D=\Dd$,
  $N = \left(\frac{\sqrt{\energyBound}D+1}{\nu}\right)^{2}$, $C>\sqrt{\energyBound}N^s$ and
  $C>\sup_{\substack{k\in\mathbb{Z}\\t\in\forcingTimeInterval}}{|k|^{s-3/2}|f_k(t)|}$.
  Then
    \begin{equation*}
         W(\energyBound, C, s)=\left\{ \{a_k\} \in \subspaceH_0 \ |\ \Ezero(\{a_k\})\leq \energyBound,\ |a_k|\leq\frac{C}{|k|^s}\right\}
    \end{equation*}
    is a \emph{trapping region} for each \gp.
\end{theorem}
\paragraph{Proof:} The proof follows the same lines as the proof of \cite[Thm. 3.4]{Cy}. What is new here, is that we consider
an arbitrary smooth forcing, and in the previous work we considered only forcing having a finite number of nonzero Fourier coefficients,
therefore we present only new elements of the required proof.

First, the value $\sup_{\substack{k\in\mathbb{Z}\\t\in\forcingTimeInterval}}{|k|^{s-3/2}|f_k|}$ is finite, because by assumption $f(t,\cdot) \in C^\infty$ uniformly for $t \in \mathbb{R}$.

The arguments below are valid for all Galerkin projections, then we can pass to the limit using a standard technique.

Now, let us check if for $a \in\partial W(\tilde{\mathcal{E}},C,s)$, such that $|a_k|=\frac{C}{|k|^s}$ for a $|k|>N$,
the condition $\frac{d |a_k|}{d t}<0$
holds, as this is the element of the proof where the forcing enters, and has to be revised.

From \cite[Lemma 3.3]{Cy} it follows that
\begin{equation}
  |N_k(a)| \leq \frac{D\sqrt{\energyBound}C}{|k|^{s-3/2}}.
\end{equation}
We have
\begin{equation*}
  \frac{d |a_k|}{d t}<-\nu|k|^2\frac{C}{|k|^s}+\frac{D\sqrt{\energyBound}C}{|k|^{s-3/2}}+\frac{\sup_{t\in\forcingTimeInterval}{|k|^{s-3/2}|f_k(t)|}}{|k|^{s-3/2}}.
\end{equation*}
Hence to have  $\frac{d |a_k|}{d t}<0$ it is enough to have
\begin{equation}
  \nu\sqrt{|k|}>D\sqrt{\energyBound}+\frac{\sup_{t\in\forcingTimeInterval}{|k|^{s-3/2}|f_k(t)|}}{C}.
\end{equation}
Therefore for any $|k|>\left(\frac{D\sqrt{\energyBound}}{\nu}+\frac{\sup_{\substack{t\in\forcingTimeInterval\\k\in\mathbb{Z}}}{|k|^{s-3/2}|f_k(t)|}}{C\nu}\right)^2$ we have
$\frac{d |a_k|}{d t}<0$. Because $|k|$ was assumed
$|k|>N=\left(\frac{\sqrt{\energyBound}D+1}{\nu}\right)^{2}>\left(\frac{D\sqrt{\energyBound}}{\nu}+\frac{\sup_{t\in\forcingTimeInterval}{|k|^{s-3/2}|f_k(t)|}}{C\nu}\right)^2$ we are done. \qed

\begin{theorem}{\rm \cite[Thm 2.9]{CyZ}}
    \label{thm:absorbingSet}
     Assume that \satisfiesFassumption.

    Let $\varepsilon>0$,  $E_0=\frac{\Ef}{\nu^2} < \infty$,
     $\energyBound>E_0$.
    Put
    \begin{eqnarray}
      s_i&=&i/2\text{ for }i\geq 2,\nonumber\\
      D_i&=&2^{s_i-\frac{1}{2}}+\frac{2^{s_i-1}}{\sqrt{2s_i-1}}\text{ for }i\geq 2,\nonumber\\
      C_2&\geq&\varepsilon+\left.\left(\jednadruga\energyBound+\sup_{\substack{k\in\mathbb{Z}/\{0\}\\t\in\forcingTimeInterval}}{\frac{|f_k(t)|}{|k|}}\right)\right/\nu,\label{eq:depOnF1}\\
      C_i&\geq&\varepsilon+\left.\left(C_{i-1}\sqrt{\energyBound}D_{i-1}+\sup_{\substack{k\in\mathbb{Z}/\{0\}\\t\in\forcingTimeInterval}}{|k|^{s_i-2}|f_k(t)|}\right)\right/\nu\text{ for }i>2,\label{eq:depOnF2}.
    \end{eqnarray}\\
    Then, there exists a sequence of \emph{absorbing sets} for large \gps, $i\geq 2$, such that
    \begin{equation*}
      \mathcal{A}_i\bigl(\energyBound,C_i, s_i\bigr)\subset \left\{\{a_k\}_{k\in\mathbb{Z}}\in \subspaceH_0 \ |\ \Ezero(\{a_k\}_{k\in\mathbb{Z}})\leq\energyBound,\ |a_k|\leq\frac{C_i}{|k|^{s_i}}\right\}.
    \end{equation*}

    For any $s \in \mathbb{N}$, $s\geq 2$ there exists $C_s$, such that $W(\energyBound,C_s,s)$ is an absorbing set  for large \gps.

\end{theorem}

\paragraph{Proof:}
The proof is essentially the same as the proof of \cite[Thm. 4.7]{Cy}, where the forcing is assumed to have a finite number of nonzero modes.
Here we consider a general forcing which is $C^\infty$ uniformly with respect to $t \in \mathbb{R}$. Therefore the values of
$\sup_{\substack{k\in\mathbb{Z}/\{0\}\\t\in\forcingTimeInterval}}{|k|^{s_i-2}|V_k(t)|}$ are finite.
To show that $\mathcal{A}_i$ are forward invariant, the argument is essentially the same, as in the proof of Theorem~\ref{thm:ns-absorbing-set}.

To show that there exists $C_s$, such that $W(V_0,s,C_s)$ is a forward invariant absorbing set, observe that it is enough to take
\begin{equation}
  C_s=\max\left\{C_i, \widetilde{C}_i\right\},
\end{equation}
where  $\widetilde{C}_i>\sqrt{\energyBound}N^{s_i}$, and $N$ is defined in Theorem~\ref{thm:analyticTrappingRegion}.

\qed

\providecommand{\omegaS}{\alpha}

\subsection{Main theorem about the attracting orbit}

\begin{theorem}
\label{thm:burgers-main}
Consider (\ref{eq:osc-burgers}).
   Assume that \satisfiesFassumption. Assume also that $\frac{\partial f}{\partial t}(t,\cdot) \in C^\infty$  uniformly with respect to $t \in \mathbb{R}$.

Let $d=1, p=2, r=1$.
Assume that $A_V,B_V  \in \mathbb{R}$ and $s_V$ are such that
\begin{eqnarray}
  s_V &>& d+p+r+1, \\
  A_V&=&\sup_{t \in \mathbb{R}, k \in \mathbb{Z}\setminus \{0\}} |k|^{s_V} |f_k(t)|, \\
  B_V&=&\sup_{t \in \mathbb{R}, k \in \mathbb{Z}\setminus \{0\}} |k|^{s_V} |f'_k(t)|.
\end{eqnarray}

There exists $\hat{\omegaS}>0$ such that for all $|\omegaS|>\hat{\omegaS}$

    \begin{itemize}
      \item there exists an eternal bounded solution  of \eqref{eq:osc-burgers} $\overline{a}\colon\mathbb{R}\to\subspaceH_0 \cap \tilde{H}$,
      \item $\overline{a}(t)=\mathcal{O}\left(\frac{A_V}{|\omegaS|}\right)+\mathcal{O}\left(\frac{A_V + B_V}{ \nu |\omegaS|}\right)$ for all $t\in\mathbb{R}$
      \item $\overline{a}$ attracts \emph{exponentially} all orbits in $\subspaceH_0 \cap \tilde{H}$.
    \end{itemize}

  \end{theorem}

\paragraph{Proof:}

In this proof we will use $\varphi$ to denote the semiprocess induced by the Galerkin projection of (\ref{eq:osc-burgers}). We will produce bounds and
use arguments that are valid for all Galerkin projections. Then we pass to the limit.

From our assumptions it follows that for $t \in \mathbb{R}$ and $k \in \mathbb{Z}\setminus \{0\}$ holds
\begin{eqnarray}
  |f_k(t)| &\leq& \frac{A_V}{|k|^{s_V}},  \label{eq:burgers-v_k-bound} \\
  |f_k'(t)| &\leq& \frac{B_V}{|k|^{s_V}}.
\end{eqnarray}

First we  find a forward invariant absorbing set $\mathcal{A}$, which forms self-consistent bounds for arbitrary $\omegaS$ for the system (\ref{eq:osc-burgers}).
We take $\mathcal{A}=W(\energyBound,C_s,s)$ for any $s\geq 2$ and sufficiently large $C_s$, i.e.  the forward invariant absorbing set obtained in Theorem~\ref{thm:absorbingSet}.
Observe that $\mathcal{A}$ is also an forward invariant absorbing set for the autonomous system (\ref{eq:osc-nog-burgers}).

We want to use Lemma~\ref{lem:ssb-estm-norm} on $\mathcal{A}$ on the time interval $[t_0,t_0+h]$ for arbitrary $t_0$ with $h$ to be specified later.
 Let us see first that its assumptions are satisfied. First of all since $\mathcal{A}$ is forward invariant we have a priori bounds  valid on any interval
 $[t_0,t_0+h]$ and in the notation of Lemma~\ref{lem:ssb-estm-norm} we can set $Z \oplus T_0 = W \oplus T_1=\mathcal{A}$. Formally this is not correct,
 but the lemma is valid also when we consider sets contained in self-consistent bounds. This is our present situation.

We do not have the term $\tilde{N}(t,a)$ and our oscillating part has the following  form
\begin{equation*}
  V(t) = \sum_{k \neq 0} \exp(ik \omegaS t) f_k(t) \exp(ikx).
\end{equation*}
Therefore in the notation used in Section~\ref{sec:lem-rapid-osc} we see that
 $g_k(t)=\exp(i t)$, $\omega_k=k \omegaS $, $v_k(t)=f_k(t)$, $G_k(t)=i \exp(i t)$, $\sigma_m=0$ (i.e. we do not have the terms $\tilde{g}_\sigma$ and $\tilde{v}_\sigma$). Observe that
 \begin{equation}
   |\omega_k| \geq \omegaS.
 \end{equation}

We have
\begin{eqnarray}
  C(g_k)&=&C(G_k)=1,   \label{eq:burgers-Cgk}\\
  C(v_k) &=& \frac{A_V}{|k|^{s_V}}, \label{eq:burgers-cv-k} \\
  C\left(\frac{\partial v_k}{\partial t}\right) &=& \frac{B_V}{|k|^{s_V}}. \label{eq:burgers-cvdt-k}
\end{eqnarray}

It is easy to see that
\begin{eqnarray}
     \sum_{\eta \in \mathbb{Z} \setminus \{0\}} C(g_\eta) C(v_\eta) &<& \infty, \\
      \sum_{\eta \in \mathbb{Z} \setminus \{0\}} G(G_\eta)  \frac{1}{|\eta|^{s_V-p}}  &<& \infty.
\end{eqnarray}
This means that all assumptions from Lemma~\ref{lem:ssb-estm-norm} are satisfied.

It follows from Lemma~\ref{lem:logNormNegCloseToZero}  that
there exists $E_->0$, $E_-$ depends on $C_i$ and $s_i$ used in the definition of $\mathcal{A}$, such that the set
\begin{equation}
\mathcal{R}=\mathcal{A} \cap \left\{a \in \subspaceH_0 \ | \  \|a\| \leq \sqrt{E_-}\right\}
\end{equation}
 satisfies
\begin{equation}
 \mu(D_a\tilde{F}(t,a),\mathbb{R} \times \mathcal{R})<0.  \label{eq:lognorm-on-R}
 \end{equation}
 It is important to notice that the estimate for the tail in $\mathcal{R}$ and $\mathcal{A}$ is the same.  Observe that (\ref{eq:lognorm-on-R}) implies that any two orbits approach each other as long as they stay in $\mathcal{R}$. However $\mathcal{R}$ might not be forward invariant for (\ref{eq:osc-burgers}).

Let
\begin{equation}
  E_1 \leq \frac{E_-}{4}.
\end{equation}

 Let us set
 \begin{equation}
   \mathcal{R}_1 = \mathcal{A} \cap \left\{a \in \subspaceH_0 \ | \  \|a\| \leq \sqrt{E_1}\right\} \subset \mathcal{R}.
 \end{equation}

From Lemma~\ref{lem:burgers-energy-decay} it follows that for  autonomous system (\ref{eq:osc-nog-burgers}) (and any  Galerkin projection) holds
\begin{equation*}
  \frac{d \|a(t)\|}{dt} \leq -\nu \|a(t)\|,
\end{equation*}
hence
\begin{equation}
  \|a(t_0+t)\| \leq \|a(t_0)\| \exp(-\lambda_1 t), \quad \mbox{for $t >0$},  \label{eq:burgers-en-decay}
\end{equation}
where
\begin{equation*}
   \lambda_1=\nu.
\end{equation*}
Let us fix $h_0>0$. Let us set
\begin{equation}
  \Delta=\frac{2\sqrt{E_1}}{3} (1-\exp(-\lambda_1 h_0)).  \label{eq:burgers-delta-E1}
\end{equation}

From Lemma~\ref{lem:ssb-estm-norm} it follows that there exits $\hat{\omegaS}$, such that for $|\omegaS| > \hat{\omegaS}$ holds
\begin{equation}
  \|y(t_0+t) - a(t_0+t)\| < \Delta, \quad t \in [0,h_0]  \label{eq:diff-y-x}
\end{equation}
where $a(t)$ is a solution of (\ref{eq:osc-nog-burgers}) and $y(t)$ is a solution of (\ref{eq:osc-burgers}), such that $a(t_0)=y(t_0) \in \mathcal{A}$.

We will prove that  for any $t_0 \in \mathbb{R}$
\begin{eqnarray}
  \varphi(t_0,[0,h_0],\mathcal{R}_1) &\subset& \mathcal{R},  \label{eq:varphi0hR1} \\
  \varphi(t_0,h_0,\mathcal{R}_1) &\subset&  \mathcal{R}_1.  \label{eq:varphihR1}
\end{eqnarray}
Since $\mathcal{A}$ is forward invariant
and contains both $\mathcal{R}$ and $\mathcal{R}_1$, so to establish (\ref{eq:varphi0hR1}) and (\ref{eq:varphihR1}) we just need to worry with the value of the norm
(or energy) of the solution starting from $\mathcal{R}_1$, only.

We  have from (\ref{eq:burgers-en-decay},\ref{eq:diff-y-x}) for $a \in \mathcal{R}_1$, $t\in [0,h_0]$ and $|\omegaS| > \hat{\omegaS}$
\begin{eqnarray}
  \|\varphi(t_0,t,a)\| &\leq& \|a\| \exp(-\lambda_1 t) + \Delta \leq \sqrt{E_1}\exp(-\lambda_1 t) + \frac{2\sqrt{E_1}}{3} (1-\exp(-\lambda_1 h_0))
  \leq \nonumber \\
   &\leq&  \frac{5\sqrt{E_1}}{3} \leq \frac{5\sqrt{E_-}}{6} < \sqrt{E_-}. \label{eq:burgers-size-orbit-inter}
\end{eqnarray}
This establishes (\ref{eq:varphi0hR1}).

To establish (\ref{eq:varphihR1}) we compute as above to obtain
\begin{eqnarray*}
   \|\varphi(t_0,h_0,a)\|  \leq \sqrt{E_1}\exp(-\lambda_1 h_0) + \frac{2\sqrt{E_1}}{3} (1-\exp(-\lambda_1 h_0)) = \\
    \sqrt{E_1} \left(\exp(-\lambda_1 h_0) + \frac{2}{3}(1-\exp(-\lambda_1 h_0)) \right)= \sqrt{E_1} \left(\frac{2}{3} + \frac{1}{3}\exp(-\lambda_1 h_0)  \right) <  \sqrt{E_1}.
\end{eqnarray*}

 Therefore from (\ref{eq:lognorm-on-R},\ref{eq:varphi0hR1}) and the standard estimates for Lipschitz constants in terms of logarithmic norm (see for example Lemma 4.1 in \cite{KZ}) it follows that
\begin{equation}
  \|\varphi(t_0,h_0,a) - \varphi(t_0,h_0,y) \| \leq e^{\mu(D_a \tilde{F},\mathbb{R} \times \mathcal{R}_1)h_0} \|a- y\|, \quad a,y \in \mathcal{R}_1.
\end{equation}
We consider now a family of time shifts by $h_0$: $\varphi(k h_0,h_0,\cdot)$ for $k \in \mathbb{Z}$. In the terminology used in \cite{CyZ} this is a discrete semiprocess. From \cite[Thm. 5.2, 6.16]{CyZ}  it follows that in $\mathcal{R}_1$ there exists $\overline{a}$ a unique orbit defined for $t \in \mathbb{R}$ (an eternal bounded solution), which attracts all other forward orbits with initial condition for $t_0=0$ in $\mathcal{R}_1$. From (\ref{eq:varphi0hR1}) it follows that this attracting orbit is contained in $\mathcal{R}$.

We will show that all orbits  with an initial condition in $\mathcal{A}$ at time $t=0$ or $t=k h_0$, $k \in \mathbb{N}$,  enter $\mathcal{R}_1$. This implies that all orbits are exponentially attracted by $\bar{a}$.
From (\ref{eq:burgers-en-decay}) and (\ref{eq:diff-y-x}) it follows that for $a \in \mathcal{A}$ we have
\begin{eqnarray*}
  \|\varphi(t_0,h_0,a) \| < \|a\| \exp(-\lambda_1 h_0) + \Delta = \\
  \|a\|  \exp(-\lambda_1 h_0) + \frac{2\sqrt{E_1}}{3} (1-\exp(-\lambda_1 h_0)) = \\
  \exp(-\lambda_1 h_0) \left( \|a\| - \frac{2\sqrt{E_1}}{3}  \right) + \frac{2\sqrt{E_1}}{3}.
\end{eqnarray*}
Hence
\begin{equation*}
  \|\varphi(t_0,h_0,a)\| < \|a\|, \quad \mbox{if $\|a\| > \frac{2\sqrt{E_1}}{3}$}.
\end{equation*}
Therefore there exists $K \in \mathbb{N}$, such that for any $t_0 \in \mathbb{R}$ holds
\begin{equation}
  \varphi(t_0,K h_0, \mathcal{A}) \subset \mathcal{R}_1.
\end{equation}

To finish the proof it remains to establish the bound for the attracting orbit. From (\ref{eq:burgers-size-orbit-inter}) it follows that this bound is given by $\|a\| \leq \frac{5\sqrt{E_1}}{3}$, where
$E_1$ depends on the size of the forcing term through the following relation (compare (\ref{eq:burgers-delta-E1}))
\begin{equation}
  \frac{3 \Delta}{2(1- \exp(-\lambda_1 h_0))} \leq \sqrt{E}_1, \qquad E_1 \leq E_-/4.
\end{equation}

  Observe that the expression for $\Delta$ has the following form (see Lemma~\ref{lem:ssb-estm-norm})
\begin{equation}
  \Delta \leq \sum_{k \in \mathbb{Z} \setminus \{0\}} \frac{\max_{t \in [0,h_0]}b_k(t)}{|k \omegaS|}= \frac{1}{|\omegaS|} \sum_{k \in \mathbb{Z} \setminus \{0\}} \frac{\max_{t \in [0,h_0]} b_k(t)}{|k|},
\end{equation}
where the expression for $b_k$ is given by (\ref{eq:bk-expression}). It should be noted that in order to estimate the size of the periodic orbit the values of $b_k(t)$ can be computed on $\mathcal{R}$. From (\ref{eq:lognorm-on-R}) since $D_a \tilde{F}=D F$ we have
\begin{equation}
l=\mu(D F,\mathcal{R})) <0. \label{eq:burgers-logNorm-neg}
\end{equation}
This implies that
\begin{eqnarray}
  1+e^{lt} \leq 2, \quad (e^{lt}-1)/l < t, \qquad \mbox{for $t \geq 0$}. \label{eq:burgers-estm-neg-logorm}
\end{eqnarray}

 After inserting into (\ref{eq:bk-expression}) our estimates (\ref{eq:burgers-Cgk},\ref{eq:burgers-cv-k},\ref{eq:burgers-cvdt-k}) and (\ref{eq:burgers-estm-neg-logorm}) we obtain (observe that  $D_z v_k(t,z)=0$, hence $C(D_z v_k \tilde{F})=0$)
\begin{eqnarray*}
  b_k(t) &\leq& \frac{2 A_V}{|k|^{s_V}} + C\left(D_z Fv_k\right)h_0 +  \frac{B_V}{|k|^{s_V}}h_0, \quad t \in [0,h_0].
\end{eqnarray*}

For $C\left(D_z Fv_k\right)$ we use (\ref{eq:DzFvjLem})  in Lemma~\ref{lem:ssb-estm-bnds} and $\beta_1=\nu$ to obtain for $t \in [0,h_0]$
\begin{eqnarray*}
 b_k(t) \leq
   \frac{2 A_V}{|k|^{s_V}} +  \left( \left(\frac{A_V}{|k|^{s_V - p}}  \left(\nu + \D  S_1(s_V-r) \right) \right) +  \frac{B_V}{|k|^{s_V}}\right) h_0.
\end{eqnarray*}
Constant $\D$ depends on $C$, $s$ (which has been fixed by the choice of $\mathcal{A}$) and the polynomial $N$. Since we construct $\mathcal{R}$, which 'isolates' the periodic orbit, as a subset of $\mathcal{A}$, we will treat $\D$ as the constant depending on $\mathcal{A}$, which is fixed as we change $\alpha$.

After summing up over $k$ we obtain
\begin{eqnarray*}
  |\omegaS|\Delta \leq  2 A_V S_1(s_V+1) +  \left( A_V S_1(s_V-p+1)  \left(\nu + \D  S_1(s_V-r) \right)  +  B_V S_1(s_V+1)\right) h_0
\end{eqnarray*}

The size of $E_1$ is given by
\begin{eqnarray*}
  \sqrt{E_1} \leq \frac{3 \cdot  A_V S_1(s_V+1) +  \frac{3}{2} \cdot \left( A_V S_1(s_V-p+1)  \left(\nu + \D  S_1(s_V-r) \right)  +  B_V S_1(s_V+1)\right) h_0}{|\omegaS| (1- \exp(-\nu h_0))}.
\end{eqnarray*}
For the optimal $E_1$ we should minimize the above expression with respect to $h_0$.

For example, if we set $h_0=\frac{1}{\nu}$ (we have the freedom of increasing of $h_0$ at the cost of increasing $\hat{\omegaS}$), then we obtain
\begin{eqnarray*}
  \sqrt{E_1}&\leq& \frac{e}{(e-1)  |\omegaS|}\left(A_V \left(3 \cdot  S_1(s_V+1) +  \frac{3}{2} \cdot  S_1(s_V-p+1)\right)  + \right. \\
    &+& \left. \frac{3}{2\nu} \cdot \left( A_V S_1(s_V-p+1)   \D  S_1(s_V-r)  +  B_V S_1(s_V+1)\right)\right)
\end{eqnarray*}

\qed

After returning to a moving coordinate frame we obtain from Theorem~\ref{thm:burgers-main} the following result.
\begin{theorem}
\label{thm:burgers-attracting-moving}
The same assumptions about $f$, $A_V$, $B_V$, $s_V$ are in Theorem~\ref{thm:burgers-main}. Consider problem (\ref{eq:burgers}) together with condition
(\ref{eq:fixedInt}).

There exists $\hat{\omegaS}>0$ such that for all $|\omegaS|>\hat{\omegaS}$
    \begin{itemize}
      \item there exists an eternal bounded solution  of \eqref{eq:osc-burgers} $\overline{a}\colon\mathbb{R}\to\subspaceH \cap \tilde{H}$,
      \item $\overline{a}(t)=\omegaS + \mathcal{O}\left(\frac{A_V}{|\omegaS|}\right)+\mathcal{O}\left(\frac{A_V + B_V}{ \nu |\omegaS|}\right)$ for all $t\in\mathbb{R}$
      \item $\overline{a}$ attracts \emph{exponentially} all orbits in $\subspaceH \cap \{a_0=\omegaS\} \cap \tilde{H}$.
    \end{itemize}
\end{theorem}

\subsection{Comments on the scaling of the attracting eternal solution in Theorem~\ref{thm:burgers-main}}
\label{subsec:comments-size}
In the context of the Theorem~\ref{thm:burgers-main}, the scaling of the attracting eternal bounded solution 
$\overline{a} = \mathcal{O}\left(\frac{A_V}{|\omegaS|}\right)+\mathcal{O}\left(\frac{A_V + B_V}{ \nu |\omegaS|}\right)$
is a bit counter-intuitive, because the part $\mathcal{O}\left(\frac{A_V}{|\omegaS|}\right)$ does not go to zero with $\nu$, which 
can easily be proved for the eternal solution $\overline{a}$. 
In fact the size of the eternal bounded solution should converge to zero when $\nu \to \infty$  or $|\alpha|\to \infty$.
This means that the bounds for the eternal solution provided by theorems \ref{thm:burgers-main}, \ref{thm:burgers-attracting-moving}
are not optimal.

This is well illustrated for an  ODE on the complex plane given by
\begin{equation}
  z'=-\nu z + e^{i \alpha t}. \label{eq:ode-planar-dissipative}
\end{equation}
The solution of (\ref{eq:ode-planar-dissipative}) is
\begin{equation}
  z(t)=z(0) e^{-\nu t} + \frac{e^{i \alpha t} }{\nu + i \alpha}. \label{eq:ode-sol}
\end{equation}
The  attracting solution is periodic with period $T=\frac{2 \pi}{\alpha}$ and is obtained from (\ref{eq:ode-sol}) with $z(0)=\frac{1}{\nu + i \alpha}$. We have
\begin{equation}
 \bar{z}(t)= z_{attr}(t) = \frac{e^{i \alpha t}}{\nu + i \alpha}.
\end{equation}

Hence we obtain in this example $\bar{z}(t)=\mathcal{O}\left(\frac{1}{\nu + |\alpha|}\right)$.  Observe that also $\bar{z}(t)=\mathcal{O}\left(\frac{1}{|\alpha|}\right)$ holds, which corresponds to the estimate established in Theorem~\ref{thm:burgers-main}.
The question is wether the scaling obtained in Theorem~\ref{thm:burgers-main} can be improved to obtain the desired dependence on $\nu$ in all terms.

Obtaining such scaling using the approach from Lemma~\ref{lem:prob1-estm-norm} appears difficult, if possible at all. Below we reproduce the computations from
this lemma for equation (\ref{eq:ode-planar-dissipative}). Observe that we have (we use notation from the proof of Lemma~\ref{lem:prob1-estm-norm}) and we have only one frequency
\begin{equation}
  M(t_0+t,t_0)=e^{-\nu t}, \quad g(t)=e^{it}, \quad G(t)=-i e^{it}, \quad v(t)=1.
\end{equation}
Therefore we have (we do the integration by parts as performed in Lemma~\ref{lem:prob1-estm-norm})
\begin{eqnarray*}
  I(t+t_0)=\int_0^t e^{-\nu (t-s)} e^{i \alpha (t_0 + s)}ds = \left. \frac{1}{i \alpha}e^{-\nu (t-s)} e^{i \alpha (t_0 + s)}\right|_{s=0}^{s=t} -
     \frac{1}{i \alpha} \int_0^s \nu e^{-\nu (t-s)} e^{i \alpha (t_0 + s)} ds = \\
     \frac{e^{i \alpha t_0}}{i \alpha} \left(e^{i \alpha t} - e^{-\nu t} \right) - \frac{e^{i \alpha t_0} \nu e^{-\nu t}}{i \alpha (\nu + i \alpha)}\left( e^{(\nu + i \alpha)t} -1 \right)= \\
     e^{i \alpha t_0} \left( \frac{1}{i \alpha}\left(e^{i \alpha t} - e^{-\nu t} \right) - \frac{\nu }{i \alpha (\nu + i \alpha)}\left( e^{i \alpha t} - e^{-\nu t} \right) \right)
\end{eqnarray*}
In the proof of Lemma~\ref{lem:prob1-estm-norm} we estimate both terms in parentheses by $\mathcal{O}\left(\frac{1}{|\alpha|} \right)$, however in the above equation they partly cancel out to give  $e^{i \alpha t_0} \frac{1}{\nu + i \alpha}=\mathcal{O}\left(\frac{1}{\nu + |\alpha|}\right)$.

%% file: ns.tex
\section{Navier-Stokes equation}
\label{sec:NSE}

\subsection{The Navier-Stokes equation in the Fourier domain}

We will use the following notation.  For $z \in \mathbb{C}$,  by
$\overline{z}$  we denote the conjugate of $z$. For any two
vectors $u=(u_1,\dots,u_n)$ and $v=(v_1,\dots,v_n)$ from $\mathbb{C}^n$ or $\mathbb{C}^\infty$ we set (if it makes sense)
\begin{eqnarray*}
  (u|v)&=&\sum_{i} u_i \overline{v_i} \\
  (u\cdot v) &=&\sum_{i} u_i v_i.
\end{eqnarray*}

Below we give a derivation of Navier-Stokes equations in the Fourier domain. This is a standard material (see \cite{ES,MS} and references given there) included here
just for the sake of making the paper reasonably self-contained.

 The general $d$-dimensional Navier-Stokes equation
(NSE) is written for $d$ unknown functions
$u(t,x)=(u_1(t,x),\dots,u_d(t,x))$ of $d$ variables
$x=(x_1,\dots,x_d)$ and time $t$, and the pressure $p(t,x)$.
\begin{subequations}
\label{eq:NSsystem}
\begin{align}
  \frac{\partial u_i}{\partial t} + \sum_{k=1}^d u_k \frac{\partial u_i}{\partial
  x_k}&=\nu \triangle u_i - \frac{\partial p}{\partial x_i} +
  f^{(i)}(t,x) \label{eq:NS} \\
  \mbox{div}\ u&=\sum_{i=1}^d \frac{\partial u_i}{\partial x_i}=0 \label{eq:div}
\end{align}
\end{subequations}
The functions $f^{(i)}(t,x)$ are the components of the external
forcing, $\nu >0$ is the viscosity.

We consider (\ref{eq:NS},\ref{eq:div}) on the torus $\mathbb{T
}_d$ (i.e. we consider the periodic boundary conditions).

An easy computation shows that
\begin{equation}
  \frac{d}{dt} \int_{\mathbb{T}_d} u(t,x)dx= \int_{\mathbb{T}_d} f(t,x)dx.  \label{eq:duo}
\end{equation}

The periodic boundary conditions enable us to use the Fourier
series. We write
\begin{equation}
  u(t,x)=\sum_{k \in \mathbb{Z}^d} u_k(t)e^{i(k,x)}, \qquad
  p(t,x)=\sum_{k \in \mathbb{Z}^d} p_k(t)e^{i(k,x)}
\end{equation}
Observe that $u_k(t) \in \mathbb{C}^d$, i.e. they are
$d$-dimensional vectors and $p_k(t) \in \mathbb{C}$. We will always
assume that
\begin{equation}
f_0=0.  \label{eq:f00}
\end{equation}

From (\ref{eq:duo}) and (\ref{eq:f00}) it follows that
\begin{equation}
  u_0(t)= u_0(0), \quad \forall t > 0, \mbox{ where the solution is defined}.
\end{equation}

Observe that (\ref{eq:div}) is reduced to the requirement $u_k
\bot k$. Namely
\begin{eqnarray*}
  \mbox{div}\ u&=& \sum_{k \in \mathbb{Z}^d} i (u_k(t),k)e^{i(k,x)} = 0 \\
    (u_k,k)&=&0 \quad k \in  \mathbb{Z}^d
\end{eqnarray*}

We obtain the following infinite ladder of differential equations
for $u_k$ (see for example \cite{ZNS})
\begin{equation}
  \frac{d u_k}{d t}=-i \sum_{k_1}(u_{k_1}|k)u_{k-k_1} - \nu k^2u_k
  -i p_k k + f_k, \qquad k \in \mathbb{Z}^d. \label{eq:NSgal}
\end{equation}
Here $f_k$ are components of the external forcing. Let $\sqcap_k$
denote the operator of orthogonal projection onto the
$(d-1)$-dimensional plane orthogonal to $k$. Observe that since
$(u_k,k)=0$, we have $\sqcap_k u_k=u_k$. We apply the projection
$\sqcap_k$ to (\ref{eq:NSgal}). The term $p_k k$ disappears and we
obtain
\begin{equation}
  \frac{d u_k}{d t}=-i \sum_{k_1}(u_{k_1}|k)\sqcap_k u_{k-k_1} - \nu k^2u_k + \sqcap_k f_k
    \label{eq:NSgal1}
\end{equation}
The pressure is given by the following formula
\begin{equation}
  -i \sum_{k_1}(u_{k_1}|k)(I - \sqcap_k)u_{k-k_1} -i p_k k + (I -
  \sqcap_k)f_k=0
\end{equation}

Observe that solutions of  (\ref{eq:NSgal1}) satisfy
the incompressibility condition $(u_k,k)=0$. The subspace of real
functions, $R$ (see Def.~\ref{def:reality}), is invariant under (\ref{eq:NSgal1}). In the
sequel, we will investigate the equation (\ref{eq:NSgal1})
restricted to $R$.

Observe that since $u_0(t)=u_0$ the system (\ref{eq:NSgal1}) becomes
\begin{equation}
 \frac{d u_k}{d t}=-i \sum_{k_1 \neq 0}(u_{k_1}|k)\sqcap_k u_{k-k_1} - (\nu k^2 + i (u_0|k))u_k + \sqcap_k f_k, \quad  k \in \mathbb{Z}^d \setminus \{0\}.
   \label{eq:NSgal2}
\end{equation}
For $u_0=0$ we obtain the following equation
\begin{equation}
 \frac{d u_k}{d t}=-i \sum_{k_1 \neq 0}(u_{k_1}|k)\sqcap_k u_{k-k_1} - \nu k^2 u_k + \sqcap_k f_k, \quad  k \in \mathbb{Z}^d \setminus \{0\}.
   \label{eq:NSgal3}
\end{equation}

According with the notation used in the previous sections we define
\begin{eqnarray*}
  (Lu)_k &=& -\nu k^2 u_k, \\
  N_k(u) &=& -i \sum_{k_1 \neq 0}(u_{k_1}|k)\sqcap_k u_{k-k_1}.
\end{eqnarray*}

\begin{definition}
Energy of $\{u_k, \quad k \in  \mathbb{Z}^d\}$  is
\begin{displaymath}
  E(\{u_k, \quad k \in  \mathbb{Z}^d\})= \sum_{k \in  \mathbb{Z}^d}|u_k|^2.
\end{displaymath}

Enstrophy of $\{u_k, \quad k \in  \mathbb{Z}^d\}$  is
\begin{displaymath}
  V(\{u_k, \quad k \in  \mathbb{Z}^d\})= \sum_{k \in  \mathbb{Z}^d}
  |k|^2|u_k|^2.
\end{displaymath}

Let
\begin{equation}
  \|y\|_1 := \sqrt{V(y)}.
\end{equation}

\end{definition}

\begin{definition}
\label{def:NStildeH}
  Let $H$ be the space $l_2(\mathbb{Z}^d,\mathbb{C}^d)$, i.e. $u \in H$ is a sequence $u:\mathbb{Z}^d\to \mathbb{C}^d$ such that
  $\sum_{k \in \mathbb{Z}^d}|u_k|^2 < \infty$  over the coefficient field $\mathbb{R}$.
  The subspace $\widetilde{H}\subset H$ is defined by
  \begin{equation*}
    \widetilde{H}:=\left\{\{u_k\}\in H\colon\text{there exists }0\leq C<\infty\text{ such that }|u_k|\leq\frac{C}{\nmid k\nmid^{d+3}}\text{ for }k\in\mathbb{Z}^d \right\}.
  \end{equation*}

  Let the space $\subspaceH$ be given by
  \begin{equation*}
    \subspaceH:=H\cap R \cap \{ (u_k,k)=0 \quad k \in \mathbb{Z}^d\}.
  \end{equation*}
  We set
   \begin{equation*}
    \subspaceH_0:=\subspaceH \cap \{ u_0=0\}.
  \end{equation*}
\end{definition}

In the sequel we consider the Galerkin projections given by the following definition.
\begin{definition}
Let us consider Galerkin projections of (\ref{eq:NSgal2}) or (\ref{eq:NSgal3}) parameterized by $n \in \mathbb{R}_+$ and denoted by $P_n$ given by
\begin{equation}
  (P_n u)_k = \left\{
                \begin{array}{ll}
                  u_k, & \hbox{if $|k| \leq n$;} \\
                  0, & \hbox{otherwise.}
                \end{array}
              \right.
\end{equation}
We will call such projections \emph{symmetric}.
\end{definition}
Observe that the reality condition (i.e. the subspace $\mathcal{R}$) is preserved by symmetric Galerkin projections, therefore also $\subspaceH$ and $\subspaceH_0$ are preserved.

\subsection{Effect of the introduction of uniformly moving coordinate system}

As in the case of the Burgers equation we consider NSE in the moving coordinate frame. Under the assumption of $f_0$ and in coordinate frame moving with the velocity $\alpha$, such that
\begin{equation}
  \alpha = \frac{1}{(2\pi)^d} \int_{\mathbb{T}_d} u(0,x)dx
  \end{equation}
  the function
  \begin{equation}
    v(t,x)=u(t,x+\alpha t) - \alpha
  \end{equation}
will satisfy NSE equation with the forcing term $\tilde{f}(t,x)=f(t,x+\alpha t)$ and
\begin{equation}
  \int_{\mathbb{T}_d}v(0,t)=0.
\end{equation}

On the Fourier series level we obtain
\begin{equation}
  \tilde{f}_k(t)=f_k(t) e^{i(k,\alpha)t}.  \label{eq:ns-f-moving}
\end{equation}
Hence, if $(k,\alpha)\neq 0$ for all $k$, such that $f_k(t) \neq 0$ for a $t \in \mathbb{R}$, then we can expect an averaging effect if $|\alpha| \to \infty$.

 If  $J$ is  finite, then  we can get the result about NSE in 2D repeating the arguments used for Burgers equation (with the trapping region  based on the enstrophy function - see Thm.~\ref{thm:ns-trap1}) under the assumption
\begin{equation}
   \forall k \in J \quad (k,\alpha) \neq 0,
\end{equation}
because then
\begin{equation}
  \inf_{k \in J} |(k,\alpha)| >0.  \label{eq:cond-non-res}
\end{equation}

We  show below this we cannot  hope for (\ref{eq:cond-non-res}) for a typical function $f$ with all nonzero modes.
\begin{lemma}
 Assume that $d > 1$ and $J=\mathbb{Z}^d \setminus \{0\}$, then for any $\alpha=(\alpha_1,\dots,\alpha_d) \in \mathbb{R}^d$ holds
 \begin{equation}
     \inf_{k \in J} |(k,\alpha)|=0.  \label{eq:inf-is-zero}
 \end{equation}
 \end{lemma}
\noindent
 \textbf{Proof:}
 We will use the Dirichlet Approximation Theorem \cite{N}, which says that for any irrational number $\beta$ there exist infinitely many fractions $p/q$, with $p,q \in \mathbb{Z}$, such that
 \begin{equation}
   \left| \beta - \frac{p}{q} \right| \leq \frac{1}{q^2}.  \label{eq:dir-approx-thm}
 \end{equation}

 It is enough we prove the result for $d=2$.   We can assume that  $\alpha_1 \neq 0$. We apply (\ref{eq:dir-approx-thm}) to $\beta=\frac{\alpha_2}{\alpha_1}$ to obtain
 \begin{equation}
   \left| \frac{\alpha_2}{\alpha_1} - \frac{p}{q} \right| \leq \frac{1}{q^2}.
 \end{equation}
 Hence when we multiply the above inequality by $|\alpha_1 q| $ we obtain
 \begin{eqnarray*}
    |(\alpha,(p,-q))|  = \left| q \alpha_2 - p \alpha_1 \right|\leq \frac{|\alpha_1|}{|q|}
 \end{eqnarray*}
 Therefore if we take $|q|$ large enough we can obtain $|(\alpha,(p,-q))|$ to be arbitrary small. \qed

\comment{
If we do not assume that $J$ is finite, then we see  two ways to proceed. In both cases we will assume that $J=\mathbb{Z}^d \setminus \{0\}$ and the following non-resonance condition
\begin{equation}
   \forall k \in J \quad (k,\alpha) \neq 0.
\end{equation}
\begin{itemize}
\item[1.] $f$ analytic and $\alpha$ satisfies the diophantine condition
  \begin{equation}
    \exists \gamma, \tau >0   \quad (k,\alpha) \geq \frac{\gamma}{\|k\|^\tau}, \ k \in \mathbb{Z}^d \setminus \{0\}. \label{eq:diof-cond}
  \end{equation}
\item[2.] $f$ is $C^\infty$ \textbf{JC ROZWINAC}
\end{itemize}
Obviously case 2 is more general, but the final result will be less quantitative than in the first case. In case 1, we should be able to follow the pattern of the proof for the Burgers
equation, the additional element is the analysis of the infinite sum of $b_k$'s  in Lemma~\ref{lem:ssb-estm-norm}, which should be convergent and decaying to $0$ with $\|\alpha\| \to \infty$ thanks to (\ref{eq:diof-cond}).
}

\comment{
\textbf{Case 2, is the one we are going to investigate. We split the the forcing $f$ into two parts: $f=f_D + f_T$ the finite (dominant part) $f_D$ and the tail $f_T$. We will
treat NSE with the tail part as the equation perturbed by $f_D$. For NSE with  forcing  $f_T$ we will establish the existence of the attracting orbit
and then we add $f_D$, whose effect will be treated as in the case of the Burgers equation, with the only difference that we will be perturbing the non-autonomous system.}
}

\subsection{Trapping regions, absorbing sets in dimension two}

\begin{lemma}
\label{lem:NS-derEnstrophy}
Let $d=2$ and consider (\ref{eq:NSgal3})  on $\subspaceH_0$. For any solution of (\ref{eq:NSgal3}) (such that all
necessary Fourier  series converge) or the symmetric Galerkin
projection of (\ref{eq:NSgal3}) we have
\begin{equation}
  \frac{dV\{u_k(t)\}}{dt} \leq - 2\nu V(\{u_k(t) \}) + 2 \sqrt{V(F)}\sqrt{ V(\{u_k(t) \})}.
   \label{eq:enstrineq}
\end{equation}
\end{lemma}

The proof can be found in many text-books, see also \cite{Si}.

Inequality (\ref{eq:enstrineq}) shows that
\begin{equation}
   \frac{dV\{u_k(t)\}}{dt} <0 , \qquad \mbox{when} \qquad
    V > V^*=\frac{V(F)}{ \nu^2}  \label{eq:enstrineq2}
\end{equation}

The following lemma is basically contained in \cite{MS} (there it is expressed using the vorticity)
\begin{lemma}
\label{lem:ns-estmLin}
Assume $d=2$ and  $\{u_k, k \in \mathbb{Z}^d\}$ is such that for some $D <
 \infty$, $\gamma > 1$
 \begin{equation}
   |u_k| \leq \frac{D}{|k|^\gamma}, \quad \mbox{and} \quad V(\{u_k\}) \leq
   V_0.
 \end{equation}
 Then  for any $\epsilon > 0$  there exists $C(\epsilon, \gamma)$, such that
 \begin{equation}
   | N_k(u)| \leq
   \frac{C(\epsilon,\gamma) \sqrt{V_0} D}{|k|^{\gamma - 1 - \epsilon}}, k \in \mathbb{Z}^d \setminus \{0\}.
 \end{equation}
\end{lemma}
\noindent
\textbf{Proof:}
We will use  the following inequality
\begin{eqnarray*}
    | (u_{k_1}|k)\sqcap_k u_{k-k_1}|=| (u_{k_1}|k - k_1)\sqcap_k u_{k-k_1}|
      \leq  |u_{k_1}| \ |k-k_1| \ |u_{k-k_1}|
\end{eqnarray*}

We will split $N_k$ into two sums $N_k = N^I_k + N^{II}_k$ and bound each sum separately.

\noindent \textbf{Sum I.} $|k_1| \leq \frac{1}{2}|k|$.

Here $|k - k_1| \geq \frac{1}{2}|k|$ and therefore $|u_{k-k_1}|\
|k- k_1| \leq \frac{D}{|k-k_1|^{\gamma-1}} \leq \frac{2^{\gamma-1}
D}{|k|^{\gamma-1}}$. Now observe that
\begin{equation}
  \sum_{|k_1| \leq \frac{1}{2}|k|} |u_{k_1}| =
  \sum_{|k_1| \leq \frac{1}{2}|k|} |k_1| \ |u_{k_1}|
  \frac{1}{|k_1|} \leq \sqrt{\sum |k_1|^2 |u_{k_1}|^2} \cdot
  \sqrt{\sum_{|k_1| < \frac{1}{2}|k|} \frac{1}{|k_1|^2}}
\end{equation}

The sum $\sum_{|k_1| < \frac{1}{2}|k|} \frac{1}{|k_1|^2}$ can be
estimated from above by a constant times an integral  of
$\frac{1}{r^2}$ over the ball of radius $\frac{1}{2}|k|$ with the
ball around the origin removed. Therefore we obtain
\begin{equation}
  \sum_{|k_1| \leq \frac{1}{2}|k|} \frac{1}{|k_1|^2} \leq
  C\int_1^{|k|/2} \frac{r dr}{r^2} \leq C \ln|k|.
\end{equation}
and
\begin{equation}
 \left| N_k^I \right|= \left| \sum_{|k_1| \leq \frac{|k|}{2}}(u_{k_1}|k)\sqcap_k u_{k-k_1}\right|
  \leq \frac{2^{\gamma -1} D}{|k|^{\gamma-1}} \sqrt{V_0} \sqrt{C}
  \sqrt{ \ln |k|} < \frac{C \sqrt{V_0} D}{|k|^{\gamma -1 -
  \epsilon}}.
\end{equation}

\noindent \textbf{ Sum II. } $\frac{1}{2}|k| < |k_1|$. We have
\begin{equation}
 |N_k^{II}| = \sum_{\frac{1}{2}|k|  < |k_1|}
     |u_{k_1}| \cdot |u_{k-k_1}|\cdot|k-k_1| \leq  D
  \sum_{\frac{1}{2}|k| < |k_1| }
  \frac{1}{|k_1|^\gamma}|u_{k-k_1}|\cdot|k-k_1|.
\end{equation}
We interpret $\sum_{\frac{1}{2}|k| < |k_1|}
 \frac{1}{|k_1|^\gamma}|u_{k-k_1}|\cdot|k-k_1|$ as a scalar product of
$|u_{k-k_1}|\cdot|k-k_1|$ and $\frac{1}{|k_1|^\gamma}$, hence, by the Schwarz
inequality, for some constant $\tilde{C}$ we obtain
\begin{equation}
 |N_k^{II}|
\leq D \sqrt{V_0}  \sqrt{\sum_{\frac{1}{2}|k| < |k_1|} \frac{1}{|k_1|^{2\gamma}}} \leq D \tilde{C} \sqrt{V_0} \frac{2^{\gamma -1}}{\sqrt{2(\gamma -1)} |k|^{\gamma-1}}.
\end{equation}
In the above computations we used the following estimate
\begin{eqnarray*}
\sum_{\frac{1}{2}|k| < |k_1|} \frac{1}{|k_1|^{2\gamma}} \leq \tilde{C_1} \int_{\frac{|k|}{2}< |k_1|} \frac{1}{|k_1|^{2\gamma}}d^2 k_1 = \tilde{C}_1 2 \pi \int_{\frac{|k|}{2}}^\infty \frac{dr}{r^{2\gamma-1}}= \frac{\tilde{C}_1 2 \pi 2^{2(\gamma-1)}}{2(\gamma-1) |k|^{2(\gamma-1)}}
\end{eqnarray*}

Observe that we used here the assumption $\gamma > 1$, which guarantees that the sum and the integral above converges.

 We obtain
\begin{equation}
|N_k| \leq  \frac{C \sqrt{V_0}D}{ |k|^{ \gamma - 1 - \epsilon } }.
\end{equation}
\qed

\begin{theorem}
\label{thm:ns-trap1}
Let $d=2$.
Assume that $V_0 > V^*$, $\gamma > 1$  and $D \in \mathbb{R}_+$.  We set
\begin{equation}
  \mathcal{W}(V_0,\gamma,D)=\left\{  \{u_k\}\in \subspaceH_0 \ | \ V(\{u_k\}) \leq V_0, \qquad
     |u_k| \leq \frac{D}{|k|^\gamma}, \quad k \in \mathbb{Z}^2 \setminus \{0\} \right\}
\end{equation}

Assume that $f(t,\cdot) \in C^\infty$ uniformly with respect $t \in \mathbb{R}$.

Then for sufficiently large $D$ set  $\mathcal{W}(V_0,\gamma,D)$ is  a forward invariant set for each symmetric Galerkin projection of (\ref{eq:NSgal3}).

\end{theorem}
\noindent
\textbf{Proof:}
Let $C=C(\epsilon=\frac{1}{2},\gamma)$
be a constant from Lemma \ref{lem:ns-estmLin}.

From our assumption about $f$ it follows that there exists $A_\gamma$, such that
\begin{equation}
  |f_k(t)| \leq \frac{A_\gamma}{|k|^\gamma}, \quad k \in \mathbb{Z}^2 \setminus \{0\}, \ t \in \mathbb{R}. \label{eq:ns-fxak}
\end{equation}

Let us take  point $u$ from the boundary of $\mathcal{W}(V_0,\gamma,D)$. We need to check if the vector field points in the direction of the interior of $\mathcal{W}(\cdot)$. We have either $V(u)=V_0$ or for some $k$ holds $|u_k| = \frac{D}{|k|^\gamma}$. In the first
case from Lemma~\ref{lem:NS-derEnstrophy} we know that $\frac{d V}{dt}(u) <0$.

To handle the second case using Lemma~\ref{lem:ns-estmLin} and (\ref{eq:NSgal3},\ref{eq:ns-fxak}) we compute as follows
\begin{eqnarray*}
  \frac{d |u_k|}{dt} \leq - \nu |k|^2 \frac{D}{|k|^\gamma} + |N_k| + |f_k| \leq - \nu |k|^2 \frac{D}{|k|^\gamma} + \frac{C \sqrt{V_0} D}{|k|^{\gamma - 3/2}} +
   \frac{A_\gamma}{|k|^\gamma}.
\end{eqnarray*}
Since we want  $\frac{d |u_k|}{dt}$ to be negative we will require that the following inequalities hold
\begin{eqnarray}
  - \frac{1}{2} \nu |k|^2 \frac{D}{|k|^\gamma} + \frac{C \sqrt{V_0} D}{|k|^{\gamma - 3/2}} &<& 0, \label{eq:ns-trap1-cond} \\
   - \frac{1}{2} \nu |k|^2 \frac{D}{|k|^\gamma} +  \frac{A_\gamma}{|k|^\gamma} &<& 0. \label{eq:ns-trap2-cond}
\end{eqnarray}

To have (\ref{eq:ns-trap1-cond},\ref{eq:ns-trap2-cond}) we need
\begin{eqnarray}
     |k| &>& \frac{4 C^2 V_0 }{\nu^2}, \label{eq:ns-tc1} \\
     D &>& \frac{2 A_\gamma}{\nu |k|^2}.  \label{eq:ns-tc2}
\end{eqnarray}
Observe that (\ref{eq:ns-tc2}) is satisfied for all $k \in \mathbb{Z}^2 \setminus \{0\}$ for $D$ large enough. Condition (\ref{eq:ns-tc1}) does not contain
$D$, but it usually it does not hold for all $k \in \mathbb{Z}^2 \setminus \{0\}$. We need it only on the boundary of our set.  For this if $u_k = \frac{D}{|k|^\gamma}$, then since $V(u) \leq V_0$ we have
\begin{equation}
  |k|^2 \frac{D^2}{|k|^{2\gamma}} \leq  V_0,
\end{equation}
hence
\begin{equation}
 |k| \geq \frac{D^{\frac{1}{\gamma-1}}}{V_0^{\frac{1}{2\gamma -2}}}. \label{eq:ns-k-on-bnd}
\end{equation}
Therefore we need (\ref{eq:ns-tc1}) for $k$'s satisfying (\ref{eq:ns-k-on-bnd}), which leads to
\begin{eqnarray*}
  \frac{D^{\frac{1}{\gamma-1}}}{V_0^{\frac{1}{2\gamma -2}}} >  \frac{4 C^2 V_0 }{\nu^2},
\end{eqnarray*}
hence
\begin{eqnarray*}
    D > \left(\frac{4 C^2 V_0^{1 + \frac{1}{2\gamma -2} } }{\nu^2}\right)^{\gamma-1}.
\end{eqnarray*}

\qed

Our next theorem shows that there exists a sequence of absorbing sets with good compactness properties. Moreover, the set $\mathcal{W}(V_0,\gamma,D)$
defined above are also absorbing sets for $D$ large enough.
\begin{theorem}
\label{thm:ns-absorbing-set}
Assume that $f(t,\cdot) \in C^\infty$ uniformly with respect to $t \in \mathbb{R}$.

Let $V_0 > V^*$.

Then there exists a sequence of  absorbing sets $\mathcal{A}_i$, $i \geq 2$, such that
\begin{eqnarray}
  \mathcal{A}_i \subset \{u\in \subspaceH_0 \ | \ V(u) \leq V_0, \quad |u_k| \leq \frac{C_i}{|k|^{\frac{i}{2}+1}}\}
\end{eqnarray}

For any $s \in \mathbb{N}$, $s\geq 2$ there exists $C_s$, such that $\mathcal{W}(V_0,s,C_s)$ (as defined in Theorem~\ref{thm:ns-trap1}) is an absorbing set
for all symmetric Galerkin projections of (\ref{eq:NSgal3})
\end{theorem}
\noindent
\textbf{Proof:}
From Lemma~\ref{lem:NS-derEnstrophy} it follows   that any solution of a symmetric Galerkin projection of (\ref{eq:NSgal3}) enters
set $\mathcal{A}_1$ defined by
\begin{equation}
  \mathcal{A}_1=\{u \in \subspaceH_0 \ | \ V(u) \leq V_0 \}
\end{equation}

On $ \mathcal{A}_1$ we have the following estimate of $N_k$ (we use  $(u_{k_1}|k)=(u_{k_1}| k-k_1)$)
\begin{eqnarray*}
  |N_k| \leq \sum_{k_1} |u_{k_1}| \cdot |k-k_1| \cdot |u_{k-k_1}| \leq \sqrt{E(u)} \cdot \sqrt{V(u)} \leq V(u) \leq V_0.
\end{eqnarray*}
Therefore on $ \mathcal{A}_1$ we have
\begin{eqnarray*}
  \frac{d |u_k|}{dt} \leq -\nu k^2 |u_k| + V_0 + \sup_{t \in \mathbb{R}}|f_k(t)| \leq -\nu k^2 \left( |u_k| - \frac{V_0 + \sup_{t \in \mathbb{R}, k }k^2 |f_k(t)|}{\nu |k|^2} \right),
\end{eqnarray*}
hence every point from $\mathcal{A}_1$ enters after finite time into set $\mathcal{A}_2$ (see \cite[Lemma 4.4.]{Cy})
\begin{equation}
 \mathcal{A}_2=\left\{ u \in \subspaceH_0 \ | \ V(u) \leq V_0, \quad |u_k| \leq \frac{C_2}{|k|^2} \right\}, \quad
   C_2= \frac{2 \left(V_0 + \sup_{t \in \mathbb{R}, k }k^2 |f_k(t)|\right)}{\nu}.
\end{equation}

Now we can setup an inductive argument.

Assume that $\mathcal{A}_i$ is given and for $u \in \mathcal{A}_i$ holds, $V(u) \leq V_0$ and $|u_k| \leq \frac{C_i}{|k|^{s_i}}$ where $s_i >1$.

Then from Lemma~\ref{lem:ns-estmLin} it follows that for $u \in \mathcal{A}_i$ holds
\begin{multline*}
   \frac{d |u_k|}{dt} \leq -\nu k^2 |u_k| + \frac{C\left(\frac{1}{2},s_i\right) \sqrt{V_0} C_i}{|k|^{s_i-3/2}}  + \sup_{t \in \mathbb{R}}|f_k(t)|  \\
     \leq -\nu k^2 \left(|u_k| - \frac{C\left(\frac{1}{2},s_i\right) \sqrt{V_0} C_i}{\nu |k|^{s_i+1/2}} - \frac{1}{\nu |k|^{s_i+1/2}}\sup_{\substack{t \in \mathbb{R}\\k\in\mathbb{Z}}}|k|^{s_i-3/2}|f_k(t)|  \right)
\end{multline*}
Hence $ \frac{d |u_k|}{dt} <0 $ if the following condition holds
\begin{equation}
  |u_k| >  \frac{C\left(\frac{1}{2},s_i\right) \sqrt{V_0} C_i+\sup_{\substack{t \in \mathbb{R}\\k\in\mathbb{Z}}}|k|^{s_i-3/2}|f_k(t)|}  {\nu |k|^{s_i+1/2}}. \label{eq:ns-entry-absorbing-set}
\end{equation}
Observe that if $|u_k| = \frac{C_i}{|k|^{s_i}}$ for $|k|$ large enough the above inequality is satisfied.

Let $\epsilon_i >0$ be arbitrary, we set
\begin{eqnarray*}
  s_{i+1}&=&s_i+ \frac{1}{2}, \\
  C_{i+1}&=& (1+\epsilon_i)\frac{C\left(\frac{1}{2},s_i\right) \sqrt{V_0} C_i+\sup_{t \in \mathbb{R}, k}|k|^{s_i-3/2}|f_k(t)|}  {\nu } \\
  \mathcal{A}_{i+1}&=& \mathcal{A}_i \cap \left\{u \in \subspaceH_0 \ | \  |u_k| \leq \frac{C_{i+1}}{|k|^{s_{i+1}}}\right\}.
\end{eqnarray*}
Observe that when $u \in \partial \mathcal{A}_{i+1}$ and $|u_k| = C_{i+1} / |k|^{ s_{i+1} }$, then from (\ref{eq:ns-entry-absorbing-set}) it follows that
\begin{equation}
  \frac{d |u_k|}{dt} < 0.
\end{equation}
hence $\mathcal{A}_{i+1}$ is forward invariant.

Since $\epsilon_i >0$, then every point from $\mathcal{A}_i$ enters $\mathcal{A}_{i+1}$ in finite time (see \cite[Lemma 4.4.]{Cy}).

To prove the assertion about $\mathcal{W}(\cdot)$ observe that
\begin{equation}
  \mathcal{A}_i \subset \mathcal{W}(V_0,s_i,C_i).
\end{equation}
From Theorem~\ref{thm:ns-trap1} it follows that we can find $C \geq C_i$, such  $\mathcal{W}(V_0,s_i,C)$ is forward invariant. But $\mathcal{A}_i \subset \mathcal{W}(V_0,s_i,C)$, hence every trajectory enters $\mathcal{W}(V_0,s_i,C)$.
\qed

\begin{rem}
\label{rem:ns-abs}
In the context of the above theorem if the forcing term $f$ is replaced by its Galerkin projection, then  absorbing sets $\mathcal{A}_i$ are also absorbing sets for equation with such modified forcing.
\end{rem}

Later we will need also  the following lemma.
\begin{lemma}
\label{lem:ns-trap-proj}
The same assumptions as in Theorem~\ref{thm:ns-trap1}. Assume that $V_0,\gamma,D$ are such that $\mathcal{W}(V_0,\gamma,D)$ is a trapping region for each symmetric Galerkin projection of (\ref{eq:NSgal3}).

Let  $V_1 >0$ and
$K>0$ is such that
\begin{equation*}
V_1 > \frac{1}{\nu^2} \sup_{t \in \mathbb{R}} V((I-P_K)f(t)).
\end{equation*}

Then  $\mathcal{W}(V_1,\gamma,D)$ is a trapping region for each symmetric Galerkin projection of (\ref{eq:NSgal3}) with the forcing equal to $(I - P_K)f(t)$.

\end{lemma}
\textbf{Proof:}
Set $\{ u \ | \ V(u) \leq V_1\}$ and $\mathcal{W}(V_0,\gamma,D)$  are forward invariant for (\ref{eq:NSgal3}) with the forcing equal to $(I - P_K)f(t)$, hence also
their intersection which is equal to  $\mathcal{W}(V_1,\gamma,D)$ is.
\qed

\subsection{Main theorem about the existence of attracting orbit}

\begin{theorem}
\label{thm:ns-main}
Let $J \subset \mathbb{Z}^d \setminus \{0\}$, $d\in\{2,3\}$.
Consider (\ref{eq:NSgal2}) with forcing term $f(t,x)=\sum_{k \in J} f_k(t) \exp(i kx)$, such that $f_0(t)\equiv 0$, and the constraint
\begin{equation}
  \int_{\mathbb{T}_d} u(t,x)dx = \alpha.
\end{equation}
Assume that
\begin{itemize}
\item $f(t,\cdot) \in C^\infty$ and $\frac{\partial f}{\partial t}(t,\cdot) \in  C^\infty$ uniformly with respect to $t \in \mathbb{R}$.
\item the following \emph{non-resonance condition} is satisfied
\begin{equation}
  (\alpha,k) \neq 0, \quad k \in J.  \label{eq:ns-nonres}
\end{equation}
\end{itemize}

 Then there exists a bounded eternal locally attracting solution $\bar{u}\in \subspaceH \cap \{u_0=\alpha \} \cap \tilde{H}$, $\bar{u}(t)=\alpha + \tilde{u}(t)$,
where $\sup_{t \in \mathbb{R}}\|\tilde{u}(t)\|\to 0$ when $|\alpha| \to \infty$.

If $d=2$, then $\bar{u}$ attracts exponentially all orbits in $\subspaceH \cap \{u_0=\alpha \}  \cap \tilde{H}$.

\end{theorem}

Before the proof we establish the following lemma.
\begin{lemma}
\label{lem:enstropy-energy}
Let $s > d+1$ and $C>0$.
Then for every $\epsilon >0$ there exists $\delta=\delta(C,s,d)$, such that if $a,b \in Z(C,s)$ (see Def.~\ref{def:ZCs}) and $\|a-b\| < \delta$, then $\|a-b\|_1 < \epsilon$.
\end{lemma}
\noindent
\textbf{Proof:}
Let us fix $\epsilon >0$. There exists $n>0$, such that
\begin{equation}
  \|(I - P_n)a\|_1 < \epsilon/2.
\end{equation}

Let $a,b \in Z(C,s)$, such that $\|a-b\|< \delta$. We have
\begin{eqnarray*}
  \|a-b\|_1^2 &\leq& \|P_n(a-b)\|^2_1 +   \|(I - P_n)a\|^2_1 +  \|(I - P_n)b\|^2_1 \leq L^2 \|P_n(a-b)\|^2 + \frac{ 2\epsilon^2}{4} \leq \\
  & & n^2\|a-b\|^2 + \frac{ 2\epsilon^2}{4} \leq   n^2 \delta^2 + \frac{ \epsilon^2}{2}.
\end{eqnarray*}
Therefore for $\delta < \frac{\epsilon}{\sqrt{2}  n}$ we obtain our assertion.
\qed

\noindent
\textbf{Proof of Theorem~\ref{thm:ns-main} for $d=2$:}

First we pass to  the coordinate frame moving with the velocity $\alpha$. From now on we will consider (\ref{eq:NSgal3}) on $\subspaceH_0$ with forcing given by
(\ref{eq:ns-f-moving}).

 The main idea of the proof is as follows.
We split the the forcing $f$ into two parts: $f=f_D + f_T$ the finite (dominant part) $f_D$ and the tail $f_T$. We will
treat NSE with the tail part as the equation perturbed by $f_D$. For NSE with  forcing  $f_T$ we will establish the existence of the attracting orbit
and then we add $f_D$, whose effect will be treated as in the case of the Burgers equation, with the only difference that we will be
perturbing the non-autonomous system. Observe that in this case the number of $\omega$'s is finite.

The decision on how to split $f$ into $f_D$ and $f_T$ is based on the following rule.
From Lemma~\ref{lem:logNormNegCloseToZero} it follows that the logarithmic norm is negative sufficiently close to zero. We will split $f$ into $f_D+f_T$,
so that for the perturbation given by $f_T$  we will have the absorbing set containing the  attracting eternal orbit.

Now we will realize the above idea.

Let us fix $V_\mathcal{A}> V^*$ and $s \in \mathbb{N}$. To apply Lemma~\ref{lem:ssb-estm-norm} observe that we have $d=2$, $p=2$, $r=1$, hence we need $s\geq 6$.

From Theorem~\ref{thm:ns-absorbing-set} it follows that there exists $C$, such that
 $\mathcal{A}=\mathcal{W}(V_\mathcal{A},s,C)$ is an absorbing set for (\ref{eq:NSgal3}) with  forcing $\tilde{f}$ or Galerkin projections of $\tilde{f}$ (see Remark~\ref{rem:ns-abs}).

From Lemma~\ref{lem:logNormNegCloseToZero}  it follows that
$\mathcal{R}=\mathcal{W}(V_-,s,C) \subset \mathcal{A}$  for $V_-$ small enough,   is a trapping region for equation (\ref{eq:NSgal1}) with $f\equiv 0$ with a negative logarithmic norm.
From Lemma~\ref{lem:ns-trap-proj} it follows that taking $K>0$ big enough $\mathcal{R}$ will be a trapping region with a negative logarithmic norm for   (\ref{eq:NSgal1}) with
the forcing $f_T=(I-P_K)f$.   For this the following must hold
\begin{equation}
   V^*_T=\frac{1}{\nu^2}\sup_{t \in \mathbb{R}} V(f_T(t))< V_-.  \label{eq:ns-VfT-small}
\end{equation}

Observe that (because the same parameters $C$ and $s$ are used to define $\mathcal{R}$ and $\mathcal{A}$)
\begin{equation}
  \mathcal{R}=\{a \in \subspaceH_0 \ | \  \|a\|_1 \leq V_-\} \cap \mathcal{A}. \label{eq:NS-R-in-A}
\end{equation}

Now we consider two problems
\begin{equation}
  \frac{du}{dt} = \tilde{F}(t,a):= \nu La + N(a) + f_D(t) + f_T(t), \label{eq:NS+osc}
\end{equation}
and
\begin{equation}
  \frac{dy}{dt} = F(t,y):= \nu Ly + N(y) +  f_T(t). \label{eq:NS-fartail}
\end{equation}

In the notation introduced above we have
\begin{equation}
 \mu(D_u\tilde{F}(t,u),\mathbb{R} \times \mathcal{R})<0.  \label{eq:NS-lognorm-on-R}
 \end{equation}
It should be stressed that the logarithmic norm in (\ref{eq:NS-lognorm-on-R}) is derived from the $\|\cdot\|$, i.e. the $l_2$ norm.

We want to use Lemma~\ref{lem:ssb-estm-norm}, with vector fields $F$ and $\tilde{F}$ as above, on $\mathcal{A}$ on the time interval $[t_0,t_0+h]$ for arbitrary $t_0$ with $h$ to be specified later.
 Let us see first that its assumptions are satisfied. First of all since $\mathcal{A}$ is forward invariant we have a priori bounds  valid on any interval
 $[t_0,t_0+h]$ and in the notation of Lemma~\ref{lem:ssb-estm-norm} we can set $Z \oplus T_0 = W \oplus T_1=\mathcal{A}$. Formally this is not correct,
 but the lemma is valid also when we consider sets contained in self-consistent bounds. This is our present situation.

We assume that $s$ is chosen big enough and we set $s_V=s+1$.

Let us denote set
\begin{equation}
  J=\{k \in \mathbb{Z}^2 \ | \  0 < |k| \leq K\}.
\end{equation}

We do not have the term $\tilde{N}(t,a)$, and our oscillating part (denoted in  Section~\ref{sec:lem-rapid-osc} and Lemma~\ref{lem:ssb-estm-norm} by $V(t)$) has the following  form
\begin{equation*}
   \sum_{k \in J} \exp(i (k \cdot \omegaS) t) f_k(t) \exp(ikx).
\end{equation*}

Therefore in the notation used in Section~\ref{sec:lem-rapid-osc} we see that for $k \in J$
 $g_k(t)=\exp(i t)$, $\omega_k=(k \cdot \omegaS) $, $v_k(t)=f_k(t)$, $G_k(t)=i \exp(i t)$ and all these quantities are zero for $k \notin J$. Observe that $\sigma_m=0$ (i.e. we do not have the terms $\tilde{g}_\sigma$ and $\tilde{v}_\sigma$).  From  non-resonance condition (\ref{eq:ns-nonres}) we have
 \begin{equation}
  \inf_{k \neq 0, |k| \leq L } |\omega_k| >0,
 \end{equation}
and from our assumption about $f$ and $\frac{\partial f}{\partial t}$ it follows that
\begin{eqnarray*}
  C(g_k)&=&C(G_k)=1,  \quad k \in J \label{eq:ns-Cgk}\\
   C(g_k)&=&C(G_k)=0,  \quad k \notin J \label{eq:ns-Cgk-2}\\
  A_V= &=& \sup_{t \in \mathbb{R}, k \in J} |k|^{s_V} |f_k(t)|,  \label{eq:ns-cv-k} \\
   C(v_k)&=& \frac{A_V}{|k|^{s_V}} \\
   B_V= &=& \sup_{t \in \mathbb{R}, k \in J} |k|^{s_V} \left|\frac{\partial f_k}{\partial t}(t)\right|, \\
  C\left(\frac{\partial v_k}{\partial t}\right) &=& \frac{B_V}{|k|^{s_V}} . \label{eq:ns-cvdt-k}
\end{eqnarray*}

It is easy to see that (because the sums are in fact finite)
\begin{eqnarray*}
     \sum_{\eta \in J} C(g_\eta) C(v_\eta) &<& \infty, \\
      \sum_{\eta \in J} G(G_\eta)  \frac{1}{|\eta|^{s_V-p}}  &<& \infty.
\end{eqnarray*}
This means that all assumptions from Lemma~\ref{lem:ssb-estm-norm} are satisfied.

Observe that (\ref{eq:NS-lognorm-on-R}) implies that any two orbits approach each other as long as they stay in $\mathcal{R}$. However, $\mathcal{R}$ might not be forward invariant for (\ref{eq:NS+osc}).

Let $V_1$ and $\eta \in (0,1)$ be such that
\begin{equation}
  V^*_T < \eta V_1 <  V_1 < V_-.  \label{eq:NS-V1}
\end{equation}

 Let us set
 \begin{equation}
   \mathcal{R}_1 = \mathcal{R} \cap \left\{a \in \subspaceH_0 \ | \  \|a\|_1 \leq \sqrt{V_1}\right\}.  \label{eq:ns-defR1}
 \end{equation}

From Lemma~\ref{lem:NS-derEnstrophy} for solutions of (\ref{eq:NS-fartail}) ( and for Galerkin projections)  with $V(y) \neq 0$ we have
\begin{equation}
  \frac{d}{dt}\sqrt{V(y(t))} \leq - \nu \sqrt{V(y(t))} + \nu \sqrt{V^*_T}
\end{equation}
Hence for $t>0$ holds (as long as $V(y(t)) \neq 0$)
\begin{equation}
  \| y(t_0+t)\|_1 \leq \left(\|y(t_0)\|_1 -\sqrt{V^*_T}\right) e^{-\lambda_1 t} + \sqrt{V^*_T}, \quad \mbox{for $t >0$},  \label{eq:NS-en-decay}
\end{equation}
where
\begin{equation*}
  \lambda_1=\nu.
\end{equation*}

Let us fix $h_0>0$ and $\Delta_1 >0$, such that
\begin{eqnarray}
  \Delta_1 &<& \sqrt{V_-} - \sqrt{V_1}\text{, and }  \label{eq:ns-delta1-1} \\
  \Delta_1 &<& \left( \sqrt{\eta V_1} - \sqrt{V_T^*}\right) \left(1-e^{-\lambda_1 h_0}\right). \label{eq:ns-delta1-2}
\end{eqnarray}

From Lemma~\ref{lem:enstropy-energy} it follows that there exists $\Delta>0$  such that for any $a,b \in \mathcal{A}$ we have the following implication
\begin{equation}
 \mbox{if} \quad \|a-b\| < \Delta, \quad \mbox{then} \quad \|a-b\|_1 \leq \Delta_1.  \label{eq:ns-impl-norm}
\end{equation}

Let $a(t_0+t)$ be a solution of (\ref{eq:NS+osc}) and $y(t_0+t)$ be a solution of (\ref{eq:NS-fartail}) with the same initial condition, $a(t_0) = y(t_0) \in \mathcal{A}$.
Since in $f_D$ we have only a finite number of frequencies and by our non-resonance condition (\ref{eq:ns-nonres}) all are non-zero, then from Lemma~\ref{lem:ssb-estm-norm}  it follows that there exists $\hat{\omegaS}$, such that for $|\omegaS| > \hat{\omegaS}$ holds   for all $a(t_0) \in \mathcal{A}$ and $t_0 \in \mathbb{R}$
\begin{equation}
  \|y(t_0+t) - a(t_0+t)\| < \Delta, \quad t \in [0,h_0].
\end{equation}
From (\ref{eq:ns-impl-norm}) we obtain for $|\omegaS| > \hat{\omegaS}$, for all $a(t_0) \in \mathcal{A}$, and $t_0 \in \mathbb{R}$ it holds that
\begin{equation}
  \|y(t_0+t) - a(t_0+t)\|_1 < \Delta_1, \quad t \in [0,h_0].  \label{eq:ns-diff-y-x}
\end{equation}

We want to prove that
\begin{eqnarray}
   \varphi(t_0,[0,h_0],\mathcal{R}_1) &\subset& \mathcal{R},  \label{eq:NS-R1-R-finv} \\
  \varphi(t_0,h_0,\mathcal{R}_1) &\subset& \mathcal{R}_1.  \label{eq:NS-R1-inv}
\end{eqnarray}
For the proof observe first that since $\mathcal{A}$ is forward invariant, so to establish (\ref{eq:NS-R1-R-finv}) and (\ref{eq:NS-R1-inv}) we just need to worry with the value of the enstrophy (or $\|\cdot\|_1$ norm)  of the solution starting from $\mathcal{R}_1$.

From (\ref{eq:ns-defR1},\ref{eq:NS-en-decay},\ref{eq:ns-delta1-1}) we have for $a \in \mathcal{R}_1$, $t \in [0,h_0]$ and $|\omegaS| > \hat{\omegaS}$
\begin{eqnarray*}
\|\varphi(t_0,t,a)\|_1 \leq \left(\sqrt{V_1} - \sqrt{V_T^*} \right) e^{-\lambda_1 t} + \sqrt{V_T^*} + \Delta_1 \leq \sqrt{V_1} + \Delta_1 < \sqrt{V_-}.
\end{eqnarray*}
This proves (\ref{eq:NS-R1-R-finv}).

To establish (\ref{eq:NS-R1-inv}) we compute as above using (\ref{eq:ns-delta1-2}) to obtain
\begin{eqnarray*}
\|\varphi(t_0,h_0,a)\|_1 \leq \left(\sqrt{V_1} - \sqrt{V_T^*} \right) e^{-\lambda_1 h_0} + \sqrt{V_T^*} + \Delta_1 < \sqrt{V_1}.
\end{eqnarray*}

We consider now a family of time shifts by $h_0$: $\varphi(k h_0,h_0,\cdot)$ for $k \in \mathbb{Z}$. In the terminology used in \cite{CyZ} this is a discrete semiprocess. From (\ref{eq:NS-R1-R-finv}),(\ref{eq:NS-R1-inv}) and \cite[Thm. 5.2, 6.16]{CyZ}  it follows that in $\mathcal{R}_1$ there exists $\overline{a}$ a unique orbit defined for $t \in \mathbb{R}$ (an eternal solution), contained in $\mathcal{R}$, which attracts all other forward orbits with initial condition for $t_0=0$ in $\mathcal{R}_1$.

We will show that all orbits  with an initial condition in $\mathcal{A}$ at time $t=0$ or $t=k h_0$, $k \in \mathbb{N}$,  enter $\mathcal{R}_1$. This implies that all orbits are exponentially attracted by $\bar{a}$.
From (\ref{eq:NS-en-decay}) and (\ref{eq:ns-diff-y-x}) it follows that for $a \in \mathcal{A}$ we have
\begin{eqnarray*}
  \|\varphi(t_0,h_0,a) \|_1 < \left(\|a\|_1 -\sqrt{V^*_T}\right) e^{-\lambda_1 h_0} + \sqrt{V^*_T} + \Delta_1 \leq \\
 \|a\|_1 e^{-\lambda_1 h_0} + \sqrt{V^*_T}(1 - e^{-\lambda_1 h_0} ) + \left(\sqrt{\eta V_1} - \sqrt{V_T^*}\right) \left(1-e^{-\lambda_1 h_0}\right)=\\
 \|a\|_1 e^{-\lambda_1 h_0}  + \sqrt{\eta V_1}  \left(1-e^{-\lambda_1 h_0}\right)= \left(\|a\|_1 -\sqrt{\eta} \sqrt{V_1}\right) e^{-\lambda_1 h_0} + \sqrt{\eta} \sqrt{V_1}
\end{eqnarray*}
Hence
\begin{equation*}
  \|\varphi(t_0,h_0,a)\|_1 < \|a\|, \quad \mbox{if $\|a\|_1 > \sqrt{\eta} \sqrt{V_1}$}.
\end{equation*}
Therefore there exists $K \in \mathbb{N}$, such that for any $t_0 \in \mathbb{R}$ holds
\begin{equation}
\label{eq:absorbingIsMapped}
  \varphi(t_0,K h_0, \mathcal{A}) \subset \mathcal{R}_1.
\end{equation}

\qed

\subsection{The proof of Theorem~\ref{thm:ns-main} for $d=3$ }

We will need several lemmas about the trapping regions for NSE in 3D.
\begin{lemma}
 \label{lem:NS3D-cs-trapping}
Assume that $s > 3$. Then set $Z(C,s)$ (see Def.~\ref{def:ZCs}) for
\begin{equation}
C < \frac{\nu}{C_2(3,s)} \label{eq:ns3d-c}
\end{equation}
(see Lemma~\ref{lem:convolution} for definition of $C_2(d,s)$) is
a trapping region for (\ref{eq:NSgal3}) \emph{with zero external force} $(f\equiv 0)$.
\end{lemma}
\noindent
\textbf{Proof:}
From Lemma~\ref{lem:convolution} it follows that for $u \in Z(C,s)$ holds
\begin{equation*}
  |N_k(u)| \leq \frac{C_2(3,s) C^2}{|k|^{s-1}}.
\end{equation*}
Therefore if $|u_k|=\frac{C}{|k|^s}$ for some $k$, then from the above and (\ref{eq:NSgal3}) it follows immediately that
\begin{equation*}
  \frac{d |u_k|}{dt} \leq -\nu |k|^2 \frac{C}{|k|^s} + \frac{C_2(3,s) C^2}{|k|^{s-1}} <0,
\end{equation*}
if
\begin{eqnarray*}
  \nu |k|  > C_2(3,s) C.
\end{eqnarray*}
Since we want this for any $k$, hence we obtain the following requirement
 $C < \frac{\nu}{C_2(3,s)}$.
\qed

\begin{lemma}
\label{lem:ns3d-cs-trapp-ft}
Assume $s>3$ and $C>0$  satisfies (\ref{eq:ns3d-c}). Assume that there exist constants $s_V \geq s$ and $A_V$, such that
\begin{equation}
  |f_k(t)| \leq \frac{A_V}{|k|^s}, \quad t \in \mathbb{R}, k \in \mathbb{Z}^3 \setminus \{0\}.
\end{equation}

Then there exists $K_0$, such that
\begin{equation}
  \frac{d|u_k|}{dt} <0, \quad \mbox{for $|k| > K_0$ and $u \in Z(C,s)$ such that $|u_k|=\frac{C}{|k|^s}$}
\end{equation}
and
$Z(C,s)$ is a trapping region for (\ref{eq:NSgal3}) with the forcing term $(I - P_n)f$ for all $n > K_0$.
\end{lemma}
\noindent
\textbf{Proof:}
On the point on the boundary of $Z(C,s)$ for some $k$ holds $|u_k|=\frac{C}{|k|^s}$. For such point we have
\begin{eqnarray*}
  \frac{d |u_k|}{dt} \leq - \nu |k|^2 \frac{C}{|k|^s}  + \frac{C_2(3,s) C^2}{|k|^{s-1}} + \frac{A_V}{|k|^{s_V}}.
\end{eqnarray*}
We want $\frac{d |u_k|}{dt} < 0$ for $|k|$ big enough.  This is implied by the following inequality
\begin{eqnarray*}
  C |k| (\nu |k| - C_2(3,s) C) > \frac{A_V}{|k|^{s_V-s}},
\end{eqnarray*}
from (\ref{eq:ns3d-c}) it follows that
\begin{eqnarray*}
  \nu |k| - C_2(3,s) C > \nu |k| - \nu.
\end{eqnarray*}
Hence it is enough to have
\begin{equation*}
  C \nu |k| (|k|-1) > \frac{A_V}{|k|^{s_V-s}}.
\end{equation*}
which leads to
\begin{equation}
  |k|^{s_V - s + 1} (|k|-1) > \frac{A_V}{C \nu}.
\end{equation}
This holds for $|k| > K_0$, so some $K_0$. In fact, since $C$ is bounded from above by (\ref{eq:ns3d-c}) we can have bound on $K_0$, which will depend on $s$
\begin{equation}
   |k|^{s_V - s + 1} (|k|-1) > \frac{A_V C_2(3,s)}{ \nu^2}.
\end{equation}
\qed

The proof of the following lemma can be found in many text-books.
\begin{lemma}
\label{lem:ns-derEnergy}
Let $u(t)$ be a solution of Galerkin projection of (\ref{eq:NSgal3}). Then
\begin{equation}
  \frac{d}{dt} E(u) \leq - 2\nu E(u) +2 \sqrt{E(u)} \cdot \sqrt{E(f)}
\end{equation}
\end{lemma}

\noindent
\textbf{Proof of Theorem~\ref{thm:ns-main} for $d=3$:}

We fix $s \in \mathbb{Z}_+$ large enough. To apply Lemma~\ref{lem:ssb-estm-norm} observe that we have $d=3$, $p=2$, $r=1$, hence we need $s\geq 7$.

  From Lemma~\ref{lem:NS3D-cs-trapping} it follows that
   for $C$ small enough,  $Z(C,s)$ is a trapping region for (\ref{eq:NSgal3}) when $f\equiv 0$. Let us fix such value of $C$.

In Lemma~\ref{lem:ns3d-cs-trapp-ft} we proved that if there exists $K_0$, such that if $n > K_0$, then $ Z(C,s)$ is a trapping region for (\ref{eq:NSgal3}) with the perturbation
   $f_T=(I-P_n)f$. We have the freedom of further increasing of $n$.

From Lemma~\ref{lem:logNormNegCloseToZero} it follows that there exists $E_- > 0$, such that the logarithmic norm for (\ref{eq:NSgal3})  is negative
  on the set $\mathbb{R} \times \mathcal{R}$, where
\begin{equation}
  \mathcal{R}=\left\{ u \in \subspaceH_0 \cap Z(C,s)   \ | \  E(u) \leq E_- \right\}.
\end{equation}

For further construction we need the following lemma.
\begin{lemma}
\label{lem:ns3d-r1r2}
Assume that $f,A,s_V,K_0$ are as in Lemma~\ref{lem:ns3d-cs-trapp-ft}.

Assume $0 < E_2 < E_1 \leq E_-$, and let us define
\begin{equation}
\mathcal{R}_i=Z(C,s) \cap \{E(u) \leq E_i\}, \quad i=1,2.
 \end{equation}

Assume that
\begin{equation}
  E_1 < \frac{C^2}{K_0^{2s}}.  \label{eq:ns3d-E1}
\end{equation}
 and  $n > K_0$ is such that
\begin{equation}
  E_1 > E_2 > \frac{1}{\nu^2} \sup_{t \in \mathbb{R}} |(I - P_n)f(t)|,
\end{equation}
then $\mathcal{R}_1$ and $\mathcal{R}_2$ are trapping regions for (\ref{eq:NSgal3}) with the forcing given by $(I-P_n)f$.

Moreover, if the forward trajectory $u(t)$ for (\ref{eq:NSgal3}) with the (full) forcing term $f$ starting in $\mathcal{R}_i$ leaves this set at time $t_e$, then $E(u(t_e))=E_i$.
\end{lemma}
\noindent
\textbf{Proof:}
Consider first (\ref{eq:NSgal3}) with the forcing given by $(I-P_n)f$.
From Lemma~\ref{lem:ns-derEnergy} it follows immediately that for $i=1,2$ set $\{u \ | \ E(u) \leq E_i\}$ is forward invariant and from Lemma~\ref{lem:ns3d-cs-trapp-ft} it follows that $Z(C,s)$ is  forward invariant. Therefore
$\mathcal{R}_i$ is also forward invariant, as the intersection two forward invariant sets. This proves the first assertion.

To establish the second assertion it is enough to show that: if $u \in \partial \mathcal{R}_i$ and $|u_k| = \frac{C}{|k|^s}$, then the vector field with full perturbation is pointing inwards.

This is achieved by demanding that any point $u \in Z(C,s)$  for which for some $k$ holds $|u_k|=\frac{C}{|k|^s}$ and  $\frac{d |u_k|}{dt} \geq 0$ has the energy larger than $E_1$. From Lemma~\ref{lem:ns3d-cs-trapp-ft} it follows that it is enough to require that
\begin{equation*}
  \frac{C}{|k|^s} > \sqrt{E_1}, \quad \mbox{for $|k| \leq K_0$}.
\end{equation*}
Observe that this holds due to assumption (\ref{eq:ns3d-E1}).

Let $u(t)$ be the solution for (\ref{eq:NSgal3}) with the (full) forcing term $f$ starting in $\mathcal{R}_i$ leaving this set at time $t_e$. At the exit moment, we cannot have  $|u_k(t_e)|= \frac{C}{|k|^s}$ for some $|k| > K_0$, because $\frac{d |u_k|}{dt}(t_e+[-\epsilon,\epsilon]) <0$ for some $\epsilon>0$. Hence  $|u_k(t_e)|= \frac{C}{|k|^s}$ cannot hold for the exit point .
\qed

We continue with the proof of Theorem~\ref{thm:ns-main} for $d=3$. Let $E_1,E_2$, $\mathcal{R}_1$ and $\mathcal{R}_2$ be as in the above lemma.

\begin{lemma}
\label{lem:ns3d-small-time-step}
Let $\mathcal{R}_1$, $\mathcal{R}_2$ be as in Lemma~\ref{lem:ns3d-r1r2}.
There exists $h>0$, such that for (\ref{eq:NSgal3})  (perturbation is  $f$) we have for any $t_0 \in \mathbb{R}$
\begin{equation}
  \varphi(t_0,[0,h],\mathcal{R}_2) \subset \mathcal{R}_1. \label{eq:ns-3d-enlo-h}
\end{equation}
\end{lemma}
\noindent
\textbf{Proof:}
From Lemma~\ref{lem:NS3D-cs-trapping} we know  that if $u(t_0) \in \mathcal{R}_1$ and $E(u(t_0+t)) < E_1$ for $t \in [0,h]$, then $u(t_0+[0,h]) \subset \mathcal{R}_1$. From Lemma~\ref{lem:ns-derEnergy} it follows that for $u \in \mathcal{R}_1$ and any $t\in \mathbb{R}$ holds
\begin{equation*}
  \frac{dE(u(t)}{dt} \leq - 2\nu E(u(t)) + 2 \sqrt{E(u)} \sqrt{E(f(t))} \leq 2 \sqrt{E_1} \sup_{t \in \mathbb{R}} \sqrt{E(f(t))}.
\end{equation*}
Therefore it is enough to take $h=(E_1-E_2)\left(\sqrt{E_1} \sup_{t \in \mathbb{R}} \sqrt{E(f(t))}\right)^{-1}$

\qed

We continue with the proof of Theorem~\ref{thm:ns-main} for $d=3$.

Now we are ready for the application of the  averaging lemma.

Let us fix $n>K_0$ and set
\begin{equation}
  f_T=(I-P_n)f, \quad f_D=P_n f, \quad E_T^*=\frac{1}{\nu^2}\sup_{t \in \mathbb{R}} E(f_T(t)).
\end{equation}
We consider two problems
\begin{equation}
  \frac{du}{dt} = \tilde{F}(t,a):= \nu La + N(a) + f_D(t) + f_T(t), \label{eq:NS3d+osc}
\end{equation}
and
\begin{equation}
  \frac{dy}{dt} = F(t,y):= \nu Ly + N(y) +  f_T(t). \label{eq:NS3d-fartail}
\end{equation}

In the notation introduced above we have
\begin{equation}
 \mu(D_u\tilde{F}(t,u),\mathbb{R} \times \mathcal{R})<0.  \label{eq:NS3d-lognorm-on-R}
 \end{equation}

We want to use Lemma~\ref{lem:ssb-estm-norm}, with vector fields $F$ and $\tilde{F}$ as above, on $\mathcal{R}_1$ on the time interval $[t_0,t_0+h]$ for arbitrary $t_0$ with $h$ obtained in Lemma~\ref{lem:ns3d-small-time-step}.
 Let us see first that its assumptions are satisfied.

We assume that $s$ is already chosen big enough, and we set $s_V=s+1$.

Let us denote set
\begin{equation}
  J=\{k \in \mathbb{Z}^3 \ | \  0 < |k| \leq n\}.
\end{equation}

We do not have the term $\tilde{N}(t,a)$ and our oscillating part has the following  form
\begin{equation*}
  V(t) = \sum_{k \in J} \exp(i (k \cdot \omegaS) t) f_k(t) \exp(ikx).
\end{equation*}

Therefore in the notation used in Section~\ref{sec:lem-rapid-osc} we see that for $k \in J$
 $g_k(t)=\exp(i t)$, $\omega_k=(k \cdot \omegaS) $, $v_k(t)=f_k(t)$, $G_k(t)=i \exp(i t)$ and all these quantities are zero for $k \notin J$. Observe that $\sigma_m=0$ (i.e. we do not have the terms $\tilde{g}_\sigma$ and $\tilde{v}_\sigma$).  From  non-resonance condition (\ref{eq:ns-nonres}) we have
 \begin{equation}
  \inf_{k \neq 0, |k| \leq n } |\omega_k| >0,
 \end{equation}
and from our assumption about $f$ and $\frac{\partial f}{\partial t}$ it follows that
\begin{eqnarray*}
  C(g_k)&=&C(G_k)=1,  \quad k \in J \label{eq:ns3d-Cgk}\\
   C(g_k)&=&C(G_k)=0,  \quad k \notin J \label{eq:ns3d-Cgk-2}\\
   A_V&=& \sup_{t \in \mathbb{R}, k \in J} |k|^{s_V} |f_k(t)|,  \label{eq:ns3d-cv-k} \\
   C(v_k)&=& \frac{A_V}{|k|^{s_V}} \\
   B_V&=& \sup_{t \in \mathbb{R}, k \in J} |k|^{s_V} \left|\frac{\partial f_k}{\partial t}(t)\right|, \\
  C\left(\frac{\partial v_k}{\partial t}\right) &=& \frac{B_V}{|k|^{s_V}} . \label{eq:ns3d-cvdt-k}
\end{eqnarray*}

It is easy to see that (because the sums are in fact finite)
\begin{eqnarray*}
     \sum_{\eta \in J} C(g_\eta) C(v_\eta) &<& \infty, \\
      \sum_{\eta \in J} G(G_\eta)  \frac{1}{|\eta|^{s_V-p}}  &<& \infty.
\end{eqnarray*}
This means that all assumptions from Lemma~\ref{lem:ssb-estm-norm} are satisfied.

Observe that (\ref{eq:NS3d-lognorm-on-R}) implies that any two orbits approach each other as long as they stay in $\mathcal{R}$. However $\mathcal{R}$ might not be forward invariant for (\ref{eq:NS3d+osc}).

From Lemma~\ref{lem:ns-derEnergy} for solutions of (\ref{eq:NS3d-fartail}) ( and for Galerkin projections)  with $E(y) \neq 0$ we have
\begin{equation}
  \frac{d}{dt}\sqrt{E(y(t))} \leq - \nu \sqrt{E(y(t))} + \nu \sqrt{E^*_T}
\end{equation}
Hence for $t>0$ holds (as long as $E(y(t)) \neq 0$)
\begin{equation}
  \| y(t_0+t)\| \leq \left(\|y(t_0)\| -\sqrt{E^*_T}\right) e^{-\lambda_1 t} + \sqrt{E^*_T}, \quad \mbox{for $t >0$},  \label{eq:NS3d-en-decay}
\end{equation}
where
\begin{equation*}
  \lambda_1=\nu.
\end{equation*}

Let us fix  $\Delta >0$, such that
\begin{eqnarray}
  \Delta &<& \left( \sqrt{E_2} - \sqrt{E_T^*}\right) \left(1-e^{-\lambda_1 h}\right). \label{eq:ns3d-delta-2}
\end{eqnarray}

Let $a(t_0+t)$ be a solution of (\ref{eq:NS3d+osc}) and $y(t_0+t)$ be a solution of (\ref{eq:NS3d-fartail}) with the same initial condition, $a(t_0) = y(t_0) \in \mathcal{A}$.
Since in $f_D$ we have only a finite number of frequencies and by our non-resonance condition (\ref{eq:ns-nonres}) all are non-zero, then from Lemma~\ref{lem:ssb-estm-norm}  it follows that there exists $\hat{\omegaS}$, such that for $|\omegaS| > \hat{\omegaS}$ holds   for all $a(t_0) \in \mathcal{R}_2$ and $t_0 \in \mathbb{R}$
\begin{equation}
  \|y(t_0+t) - a(t_0+t)\| < \Delta, \quad t \in [0,h].
\end{equation}

We want to prove that
\begin{eqnarray}
  \varphi(t_0,h_0,\mathcal{R}_1) &\subset& \mathcal{R}_1  \label{eq:NS3d-R2-inv}
\end{eqnarray}
 For the proof it is enough the compute the energy (or a norm)  of the solution starting from $\mathcal{R}_2$.

From (\ref{eq:NS3d-en-decay},\ref{eq:ns3d-delta-2}) we have for $a \in \mathcal{R}_2$, $t \in [0,h]$ and $|\omegaS| > \hat{\omegaS}$
\begin{eqnarray*}
\|\varphi(t_0,h,a)\| \leq \left(\sqrt{E_2} - \sqrt{E_T^*} \right) e^{-\lambda_1 h} + \sqrt{E_T^*} + \Delta < \sqrt{E_2}.
\end{eqnarray*}

We consider now a family of time shifts by $h_0$: $\varphi(k h_0,h_0,\cdot)$ for $k \in \mathbb{Z}$. In the terminology used in \cite{CyZ} this is a discrete semiprocess. From (\ref{eq:NS3d-R2-inv}), (\ref{lem:ns3d-small-time-step}) and \cite[Thm. 5.2, 6.16]{CyZ}  it follows that in $\mathcal{R}_2$ there exists $\overline{a}$ a unique orbit defined for $t \in \mathbb{R}$ (an eternal solution), which attracts all other forward orbits with initial condition for $t_0=0$ in $\mathcal{R}_2$.
\qed

\subsection{Alternative approach to the proof of Theorem~\ref{thm:ns-main} }

In this section we present a lemma, which gives a different proof  for the part of Theorem~\ref{thm:ns-main} related to the global attraction
of all solutions to a small eternal orbit. The argument is of global nature and  it replaces entirely the arguments of negativeness of the logarithmic norm
 in a local neighborhood \eqref{eq:NS-lognorm-on-R}, and then showing that the absorbing set is eventually mapped into this neighborhood \eqref{eq:absorbingIsMapped}.

We remark that the presented approach is valid for the 2D Navier-Stokes equations exclusively.

\paragraph{Notation} For a differentiable function $u\colon \mathbb{T}^2\to \mathbb{R}^2$ we are going to use the following notation to denote components of its partial derivatives
$ u_{x_k}^j:=\frac{\partial u^j}{\partial x_k}$.

\begin{definition}
Let $u\in H^\prime_0$. We define the norm $\|\nabla u\|_{\infty}$ as follows
\begin{equation*}
  \|\nabla u\|_{\infty}  =\|u^1_{x_1}\|_{L^\infty(\mathbb{T}^2)}+\|u^1_{x_2}\|_{L^\infty(\mathbb{T}^2)}+
  \|u^2_{x_1}\|_{L^\infty(\mathbb{T}^2)}+\|u^2_{x_2}\|_{L^\infty(\mathbb{T}^2)}.
\end{equation*}
\end{definition}

\begin{lemma}
\label{lem:ns-small-attracts-all}
  Let $d=2$, $J \subset \mathbb{Z}^2 \setminus \{0\}$.
Consider (\ref{eq:NSgal2}) with forcing term $f(t,x)=\sum_{k \in J} f_k(t) \exp(i kx)$, such that $f_0(t)\equiv 0$, and the constraint
\begin{equation}
  \int_{\mathbb{T}_d} u(t,x)dx = \alpha.
\end{equation}

Assume that there exists a bounded eternal orbit $\overline{u}\colon \mathbb{R}\to \subspaceH_0 \cap \tilde{H}$ satisfying the following bound
\begin{equation}
\|\overline{u}(t) \|_\infty < \nu,\text{ for all }t\in\mathbb{R},  \label{eq:ns-small-orbit}
\end{equation}

then $\overline{u}$ attracts exponentially all orbits in $\subspaceH_0 \cap \tilde{H}$.
\end{lemma}

\noindent
\textbf{Proof:}
First, let $u$, $v\colon\mathbb{R}_+\to \subspaceH_0 \cap \tilde{H}$ be  solutions of \eqref{eq:NSsystem}
, we subtract from
\begin{equation*}
  u_t+u\cdot\nabla u - \nu\Delta u +\nabla p_1 = f,
\end{equation*}
the equation
\begin{equation*}
  v_t+v\cdot\nabla v - \nu\Delta v +\nabla p_2 = f,
\end{equation*}
and obtain the following equation for the difference $w=u-v$, with $\widetilde{p}=p_1-p_2$,
\begin{equation*}
  w_t + u\cdot\nabla w + w\cdot\nabla u - w\cdot\nabla w - \nu\Delta w +\nabla\widetilde{p} = 0.
\end{equation*}
We multiply the above equation by $w$ and then integrate it over the torus
\providecommand{\TT}{\mathbb{T}^2}
\begin{equation*}
  \int_{\TT}{w_t\cdot w} + \int_{\TT}{(u\cdot\nabla w)\cdot w} + \int_{\TT}{(w\cdot\nabla u)\cdot w} -
\int_{\TT}{(w\cdot\nabla w)\cdot w} - \int_{\TT}{\nu\Delta w\cdot w} +\int_{\TT}{\nabla\widetilde{p}\cdot w} = 0.
\end{equation*}
Obviously $\int_{\TT}{(w\cdot\nabla w)\cdot w}=0$,
$\int_{\TT}{\nabla\widetilde{p}\cdot w}=-\int_{\TT}{\widetilde{p}\cdot {\rm div\,} w}=0$, and it follows that
(with all integrals taken over $\mathbb{T}^d$), using integration by parts, and the zero divergence condition for $w$,
that $\int_{\TT}{(u\cdot\nabla w)\cdot w}=0$, namely
\begin{multline*}
  \int{(u\cdot\nabla w)\cdot w}=\int{(u^1w^1_{x_1}+u^2w^1_{x_2})w^1+(u^1w^2_{x_1}+u^2w^2_{x_2})w^2}=\\
-\int{w^1(u^1w^1)_{x_1}+w^1(u^2w^1)_{x_2}+w^2(u^1w^2)_{x_1}+w^2(u^2w^2)_{x_2}}=\\
-\int{w^1u^1w^1_{x_1}+w^1u^2w^1_{x_2}+w^2u^1w^2_{x_1}+w^2u^2w^2_{x_2}}=
\int{(w^1)^2u^1_{x_1}+(w^1)^2u^2_{x_2}+(w^2)^2u^1_{x_1}+(w^2)^2u^2_{x_2}}=0.
\end{multline*}

We obtain
\begin{equation*}
  \frac12\frac{\partial}{\partial t}\int_{\TT}{|w|^2} = - \int_{\TT}{(w\cdot\nabla u)\cdot w}
  - \nu\int_{\TT}{|\nabla w|^2},
\end{equation*}
We estimate
\begin{equation*}
  \frac12\frac{\partial}{\partial t}\int_{\TT}{|w|^2} \leq \left|\int_{\TT}{(w\cdot\nabla u)\cdot w}\right|
  - \nu\int_{\TT}{|w|^2}\leq \|\nabla u\|_{\infty}\int_{\TT}{|w|^2} - \nu\int_{\TT}{|w|^2},
\end{equation*}
Now, we set $u=\overline{u}$, and from Gronwalls' lemma it follows
\begin{equation*}
  \|w(t)\|_{L^2(\mathbb{T}^d)}^2\leq \exp\left(2(\|\nabla \overline{u}\|_{\infty}-\nu)t\right)\|w(0)\|_{L^2(\mathbb{T}^d)}^2.
\end{equation*}
Thus, as we assume that there exists $\overline{u}$ -- an eternal solution to \eqref{eq:NSsystem}, such that
$\|\nabla \overline{u}(t)\|_{\infty}<\nu$ for all $t\in\mathbb{R}$, $\overline{u}$ exponentially attracts
forward in time any other solution of \eqref{eq:NSsystem}.
\qed

Using Lemma~\ref{lem:ns-small-attracts-all} the proof of Theorem~\ref{thm:ns-main} for $d=2,3$ can be done as follows.
\begin{description}
 \item[1.] Just as in the 3D-case we construct small trapping regions for $|\alpha| \to \infty$, so that any eternal orbit in them will satisfy condition
     (\ref{eq:ns-small-orbit}).
  \item[2.] From Lemma~\ref{lem:ns-small-attracts-all}  we know that an eternal orbit satisfying (\ref{eq:ns-small-orbit})  attracts all orbits.
\end{description}
